\documentclass[11pt,english,reqno]{amsart}
\usepackage[dvipsnames]{xcolor}
\usepackage{comment}
\usepackage{bm}
\usepackage{color}
\usepackage{babel}
\usepackage{verbatim}
\usepackage{amsthm}
\usepackage{amstext}
\usepackage{amssymb}
\usepackage{amsbsy}
\usepackage{ifpdf}
\usepackage[nobysame,abbrev,alphabetic]{amsrefs}
\usepackage{enumerate}
\usepackage{amsmath} 
\usepackage{amscd}
\usepackage{amsfonts}
\usepackage{graphicx}
\usepackage[all]{xy}
\usepackage{verbatim}
\usepackage{gensymb}
\usepackage{mathrsfs}
\usepackage[right,displaymath,mathlines]{lineno}

\usepackage[original]{imakeidx}
\usepackage{tikz}
\usepackage{pgfplots}
\pgfplotsset{compat = newest}
\usepackage{tikz-3dplot}
\usetikzlibrary{shapes.geometric}

\makeindex 

\usepackage[colorlinks]{hyperref}
\hypersetup{
linkcolor=blue,          
citecolor=green,        
} 

\newtheorem{theorem}{Theorem}[section]
\newtheorem{claim}[theorem]{Claim}
\newtheorem{lemma}[theorem]{Lemma}
\newtheorem{corollary}[theorem]{Corollary}
\theoremstyle{definition}\newtheorem{definition}[theorem]{Definition}

\theoremstyle{theorem}\newtheorem{proposition}[theorem]{Proposition}

\theoremstyle{definition}

\theoremstyle{definition}\newtheorem{remarks}[theorem]{Remarks}
\theoremstyle{definition}\newtheorem{remark}[theorem]{Remark}
\theoremstyle{definition}
\newcommand{\al}{\alpha}
\newcommand{\be}{\beta}
\newcommand{\ga}{\gamma}
\newcommand{\Ga}{\Gamma}
\newcommand{\del}{\delta}
\newcommand{\Del}{\Delta}
\newcommand{\lam}{\lambda}
\newcommand{\Lam}{\Lambda}
\newcommand{\eps}{\epsilon}

\newcommand{\vre}{\varepsilon}

\newcommand{\sig}{\sigma}

\newcommand{\Om}{\Omega}
\newcommand{\vphi}{\varphi}

\newcommand{\cB}{\mathcal{B}}
\newcommand{\cC}{\mathcal{C}}
\newcommand{\cD}{\mathcal{D}}

\newcommand{\cF}{\mathcal{F}}

\newcommand{\cM}{\mathcal{M}}

\newcommand{\cO}{\mathcal{O}}
\newcommand{\cP}{\mathcal{P}}

\newcommand{\cS}{\mathcal{S}}

\newcommand{\cU}{\mathcal{U}}
\newcommand{\cV}{\mathcal{V}}
\newcommand{\cW}{\mathcal{W}}

\newcommand{\cY}{\mathcal{Y}}

\newcommand{\XX}{{\mathscr{X}}}

\newcommand{\XXnA}{{\mathscr{X}_n^{\bA}}}
\newcommand{\SXX}{{\mathcal{S}\mathscr{X}}}

\newcommand{\bA}{\mathbb{A}}

\newcommand{\bG}{\mathbb{G}}
\newcommand{\bL}{\mathbb{L}}

\newcommand{\bR}{\mathbb{R}}
\newcommand{\R}{\mathbb{R}}
\newcommand{\bZ}{\mathbb{Z}}
\newcommand{\Z}{\mathbb{Z}}
\newcommand{\bQ}{\mathbb{Q}}
\newcommand{\Q}{\mathbb{Q}}

\newcommand{\bK}{\mathbb{K}}
\newcommand{\bN}{\mathbb{N}}
\newcommand{\N}{\mathbb{N}}

\newcommand{\bT}{\mathbb{T}}
\newcommand{\bS}{\mathbb{S}}

\newcommand{\bx}{\mathbf{x}}
\newcommand{\by}{\mathbf{y}}
\newcommand{\bp}{\mathbf{p}}
\newcommand{\bv}{\mathbf{v}}
\newcommand{\bu}{\mathbf{u}}
\newcommand{\bw}{\mathbf{w}}
\newcommand{\bz}{\mathbf{z}}

\newcommand{\SL}{\operatorname{SL}}
\newcommand{\SO}{\operatorname{SO}}

\newcommand{\PGL}{\operatorname{PGL}}

\newcommand{\GL}{\operatorname{GL}}

\newcommand{\defi}{\overset{\on{def}}{=}}
\newcommand{\df}{\overset{\on{def}}{=}}

\newcommand\norm[1]{||#1||}
\newcommand\wt[1]{\widetilde{#1}}

\newcommand\set[1]{\left\{#1\right\}}
\newcommand\pa[1]{\left(#1\right)}
\newcommand\idist[1]{\langle#1\rangle}

\newcommand\av[1]{|#1|}
\newcommand\on[1]{\operatorname{#1}}
\newcommand\diag[1]{\operatorname{diag}\left(#1\right)}

\newcommand{\disp}{\mathrm{disp}}

\newcommand{\prim}{\operatorname{prim}}

\newcommand\mb[1]{\mathbf{#1}}

\newcommand\tb[1]{\textbf{#1}}
\newcommand\mat[1]{\pa{\begin{matrix}#1\end{matrix}}}
\newcommand\br[1]{\left[#1\right]}
\newcommand\smallmat[1]{\pa{\begin{smallmatrix}#1\end{smallmatrix}}}
\newcommand\crly[1]{\mathscr{#1}}

\newcommand{\spa}{\on{span}}
\newcommand{\supp}{\on{supp}}
\newcommand{\sm}{\smallsetminus}

\newcommand{\sro}{\cS_{r_0}}
\newcommand{\srosharp}{\cS^{\sharp}_{r_0}}

\newcommand{\sr}{\cS_r}
\newcommand{\srotilde}{\wt{\cS}_{r_0}}

\newcommand{\srtilde}{\tilde{\cS}_r}
\newcommand{\srogeps}{\cS_{r_0,\ge\vre}}
\newcommand{\sroleps}{\cS_{r_0,<\vre}}

\newcommand{\srleps}{\cS_{r,<\vre}}
\newcommand{\uro}{\cU_{r_0}}

\newcommand{\xna}{X_n^{\bA}}
\newcommand{\xn}{X_n}
\newcommand{\Gaa}{\Gamma_{\bA}}

\newcommand{\onto}{\xymatrix{\ar@{>>}[r]&}}

\newcommand{\usnote}[1]{\marginpar{\color{cyan}\tiny [US] #1}}

\newcommand{\red}[1]{\textcolor{red}{#1}}

\newcommand {\ignore}[1]  {}

\begin{document}
\title{Geometric and arithmetic aspects of approximation vectors}
\author{Uri Shapira}
\address{Dept. of Mathematics, Technion, Haifa, Israel
{\tt ushapira@tx.technion.ac.il}
}
\author{Barak Weiss}
\address{Dept. of Mathematics, Tel Aviv University, Tel Aviv, Israel
{\tt barakw@post.tau.ac.il}}

\maketitle
\begin{abstract}
  Let $\theta\in \R^d$. We associate three objects to each
  approximation $(\bp, q) \in \Z^d \times \N$ of $\theta$: the  
  projection  of the lattice  $\Z^{d+1}$ to $\R^d$, along the
  approximating vector $(\bp, q)$; the displacement
  vector $(\bp -q\theta )$; and the residue classes of the
  components of the  $(d+1)$-tuple $(\bp, q)$ modulo all primes. All
  of these   have been studied in
  connection with Diophantine approximation problems. We
  consider the asymptotic distribution of all of these quantities,
  properly rescaled, as $(\bp,
  q)$ ranges over the best approximants and $\vre$-approximants of
  $\theta$, and describe limiting measures on the relevant
  spaces, which hold for Lebesgue a.e.\;$\theta$. We also
  consider a similar problem for vectors $\theta$ whose components,
  together with $1$, span a totally real number field of degree
  $d+1$. Our technique involve recasting the problem as an
  equidistribution problem for a cross-section of a one-parameter flow
  on an adelic space, which is a fibration over the space of
  $(d+1)$-dimensional lattices. Our results generalize results of many
  previous authors, to higher
  dimensions and to joint equidistribution.
\end{abstract}
\tableofcontents

\section{Introduction}\label{sec: intro}
Our results concern the asymptotic statistics of certain geometric and
arithmetic quantities associated with approximation vectors. We begin by
introducing the notions and notations necessary for
formulating the results (more details will be given in subsequent
sections of the paper).

Throughout this paper, $d$ and $n$ are positive integers with $n =
d+1$. Let $\mb{e}_1, \ldots, \mb{e}_n$ be the standard basis of
$\R^n$, and let $\mb{v} = \sum v_i \mb{e}_i  \in \R^n$ with $v_n \neq 0$.
We will continuously use the direct sum decomposition \index{R@$\R^d$
  -- horizontal space}  
\begin{equation} \label{eq: direct sum decomposition}
\R^n = \R^d \oplus \spa (\mb{v}), \ \ \text{ where } \R^d \df
\{\mb{u} \in \R^n : u_n=0\},
\end{equation}
and the first summand in this decomposition will be called the {\em
  horizontal space}. The projection $\R^n \to 
\R^d$ will be denoted by
\index{P@$\pi^{\mb{v}}_{\R^d}$ -- projection to $\R^d$ along $\bv$} 
$\pi^{\mb{v}}_{\R^d}$, and $\pi_{\R^d}$ will be an abbreviation for
$\pi_{\R^d}^{\mathbf{e}_n}$.

Choose some norm on $\R^d$, and for $\by\in \bR^d$, denote
\index{D@$\idist{\by}$ -- distance of $\by$ to $\Z^d$} 
$\idist{\by} = \min_{\bp \in \Z^{d}} \|\by- \bp\|.$
For $\theta \in \R^d$, we say that $(\bp, q) \in \Z^{d} \times \N$
is a {\em best approximation} \index{best approximation} of $\theta$
if for any $q' < q$, $\idist{q\theta}<\idist{q'\theta}$, and
$\idist{q\theta} =\|q\theta -\bp\|$.
For all $\theta \in
\R^d \sm \Q^d$, the set of best
approximations of $\theta$ is an infinite sequence $\left(\bp_k,
  q_k\right)_{k \in \N}$, where
the order is chosen so that $q_1 \leq q_2 \leq \cdots$. In fact there
will be $k_0$ such that $q_{k+1} > q_k$ for all $k > k_0$, and the $\bp_k$ will be
  uniquely determined for $k>k_0$. The potential ambiguity of choices
  for $k \leq k_0$ will have no effect on the
objects we will consider in this paper.  For $\theta \in
\R^d$ and $\bv = (\bp, q) \in 
\Z^d \times \N,$ we will write 
\begin{equation}\label{eq: def displacement vector}
\on{disp}(\theta, \bv) \df q^{1/d}(\bp - q\theta ) \in \bR^d,
\end{equation}
\index{D@$\on{disp}(\theta, \bv) $ -- displacement}
and refer to this vector as the {\em displacement}.

All measure spaces in this paper will be standard Borel spaces and all
measures will be Borel measures. The collection of probability
measures on a measure space $X$ is denoted \index{P@$\cP(X)$
  -- probability measures on $X$} by $\cP(X)$. Two measures
on $X$ are said to 
be in \textit{the same measure class} if they are mutually 
absolutely continuous (i.e.\ have the same
null-sets). If $X$ is a locally compact second countable 
Hausdorff (lcsc) space, $\mu$ is a regular measure on $X$, and $(x_k)_{k \in \N}$ is
a sequence in $X$, we say that {\em $(x_k)$ equidistributes with respect to
$\mu$} \index{equidistributes}  if the measures $\frac{1}{N}\sum_{k=1}^N \delta_{x_k}$ converge
weak-* to $\mu$.  

A lattice in $\R^n$ is a discrete
subgroup of full rank. Its covolume is the volume of a fundamental
domain in $\R^n$, and the space of all lattices of covolume 1 is
denoted \index{X@$\XX_n = \SL_n(\R)/\SL_n(\Z)$ -- space of lattices in $\R^n$} by
$\XX_n$. It can be identified with the quotient 
$\SL_n(\R)/\SL_n(\Z)$ via the map $g\Z^n \mapsto g \SL_n(\Z)$, and
thus acquires a natural $\SL_n(\R)$-invariant probability measure,
which we denote by $m_{\XX_n}$. \index{M@$m_{\XX_n}$ -- Haar-Siegel
  measure on $\XX_n$} This measure is sometimes called the
{\em Haar-Siegel measure}. Two lattices in $\R^n$ are 
{\em homothetic} if one can be obtained from the other by
multiplication by a nonzero scalar, and the homothety class of each
$\Lambda$ contains
a unique representative in $\XX_n$, which we denote by
\index{H@$[\Lam]$ -- homothety class of $\Lam$} $[\Lambda]$. If $\Lambda$ is
a lattice in $\R^n$ and 
$\mb{v} \in \Lambda \sm \R^d$ then $\pi_{\R^d}^{\mb{v}} (\Lambda)$ is
a lattice in $\R^d$ 
and $\left[\pi_{\R^d}^{\mb{v}} (\Lambda) \right] \in \XX_d.$

We denote by $\widehat \Z$ \index{P@$\widehat \Z$  -- profinite
  completion of $\Z$} the profinite completion of $\Z$, that is,
the inverse limit of groups $\Z/m\Z$, with respect to the natural maps
$\Z/m_1\Z \to \Z /m_2\Z$ whenever $m_2 | m_1$. This is a compact topological
ring which is isomorphic to $\prod_{p \text{ prime}} \Z_p$, where
$\Z_p$ is the ring of $p$-adic integers. We denote by
\index{P@$\widehat{ \Z}^n$ -- profinite completion of $\Z^n$}
$\widehat{ \Z}^n$ the additive group which is the $n$-fold Cartesian
product of $\widehat \Z$; it is also the inverse
limit of quotient groups 
$\Z^n/\Lambda$, where $\Lambda$ ranges over finite index subgroups of $\Z^n$. We
denote by $m_{\widehat{ \Z}^n}$ \index{H@ $m_{\widehat{ \Z}^n}$ --
  Haar measure on $\widehat{\Z}^n$} the Haar probability measure on $\widehat{
\Z}^n$. For $\Lambda \in \XX_n$, a vector $v \in \Lambda$ is said to
be {\em primitive} \index{primitive vector} if it is not a 
multiple of a  vector in $\Lambda$ by an integer different from $\pm
1$, and we denote the primitive elements of $\Lam$ by
$\Lam_{\mathrm{prim}}$.\index{P@$\Lam_{\mathrm{prim}}$ -- primitive
  vectors in $\Lam$} The natural diagonal embedding 
   $\Z^n \hookrightarrow
\widehat{ \Z}^n$ has a dense image, and we 
denote the closure of the image of $\Z^n_{\mathrm{prim}}$ by $\widehat{
\Z}^n_{\mathrm{prim}}$. \index{P@$\widehat{
\Z}^n_{\mathrm{prim}}$ -- closure of $\Z^n_{\prim}$ in $\widehat{Z}^n$} Note that $\widehat{
\Z}^n_{\mathrm{prim}}$ is one orbit for the natural action of the
group $\SL_n(\widehat \Z)$, and  
$$\widehat{
\Z}^n_{\mathrm{prim}} = \{(v_p)_{p \ \mathrm{prime} }  : \forall p, \,
\| v_p\|_p=1\}.$$
Since $\SL_n(\widehat \Z)$ acts
transitively on $\widehat{
\Z}^n_{\mathrm{prim}}$,  there is a unique invariant probability
measure for this action, which we denote by $m_{ \widehat{
\Z}^n_{\mathrm{prim}}}.$ \index{M@$m_{ \widehat{
\Z}^n_{\mathrm{prim}}}$ -- the $\SL_n(\hat{\Z})$-invariant probability
measure on $\widehat{
\Z}^n_{\mathrm{prim}} $}

\begin{theorem}\label{thm: main best}
For any norm $\norm{\cdot}$ on $\bR^d$ there is a probability measure
$\mu = \mu_{\on{best},\norm{\cdot}}$  
on  $\XX_d \times \bR^d \times \widehat{\bZ}^n$ such that
for Lebesgue almost any $\theta\in\bR^d$, the following holds. Let
$\mb{v}_k\in\bZ^n$ be the 
sequence of best approximations to $\theta$ with  
respect to the norm $\norm{\cdot}$. Then the sequence
\begin{equation}\label{eq: seq best}
  \left(\br{\pi^{\mb{v}_k}_{\bR^d}(\bZ^n)}, \on{disp}(\theta, \mb{v}_k),
 \mb{v}_k\right)_{k \in \N} \in \XX_d \times \bR^d \times
\widehat{\bZ}^n
\end{equation}
 equidistributes with respect to $\mu$. The measure $\mu$ has the
 following properties: 
\begin{enumerate}
\item It is a product $\mu = \mu^{(\infty)} \times \mu^{(f)}$
  where $\mu^{(\infty)}\in \cP(\XX_d \times \bR^d),$ 
  $\mu_f\in\cP(\widehat{\bZ}^n)$.
  \item The measure $\mu^{(f)}$ is $m_{ \widehat{
\Z^n}_{\mathrm{prim}}}$ (and in particular, does not depend on the
choice of the norm). 
\item The projection $\mu^{(\XX_d)}$ of $\mu^{(\infty)}$ to $\XX_d$ is equivalent to
  $m_{\XX_d}$ 
  (i.e., satisfies $\mu^{(\XX_d)} \ll
    m_{\XX_d} \ll \mu^{(\XX_d)} $), but is equal to it only in case $d=1$. 
\item \label{item: 5} The projection $\mu^{(\R^d)}$ of $\mu^{(\infty)}$ to
  $\bR^d$ is boundedly supported, absolutely continuous  
w.r.t.\,Lebesgue with a nontrivial density (i.e., is not the
restriction of Lebesgue measure to a 
subset of $\bR^d$). If $\norm{\cdot}$ is the Euclidean norm, then 
$\mu^{(\R^d)}$ and $\mu^{(\XX_d)}$ are 
$\SO_d(\R)$-invariant.

\item For $d>1$, $\mu^{(\infty)}
  \neq \mu^{(\XX_d)} \times \mu^{(\R^d)}$. 
\end{enumerate}
In particular, each of the
coordinate sequences
\begin{equation} \label{eq: coordinate sequences}
  \left(\br{\pi^{\mb{v}_k}_{\bR^d}(\bZ^n)}
  \right)
  \subset
 \XX_d, \ \ \ \left(\on{disp}(\theta, 
   \mb{v}_k)\right)
 \subset \R^d, \ \ \ (
 \mb{v}_k)
 \subset
 \widehat{\bZ}^n
 \end{equation}
 equidistributes in its respective space, with respect to the
 pushforward of $\mu^{(\crly{X}_d)}, \mu^{(\bR^d)}, \mu^{(f)}$
 respectively.
 \end{theorem}

Let $\vre>0$ and fix a norm $\| \cdot \|$ on $\bR^d$. Given $\theta \in \R^d$, we
say that $\bw = (\bp, q) \in \Z^d 
\times \N$ is an \index{E@$\vre$-approximation} {\em $\vre$-approximation of $\theta$ (with respect to
  $\| \cdot \|$)} if $\|\disp(\theta, \bw)\| \leq \vre$ and $\gcd(p_1,
\ldots, p_d, q)=1$. Standard results in Diophantine 
approximation imply that for a.e.\;$\theta$, for all $\vre$, there are
infinitely many $\vre$-approximations.
When we refer to the
sequence $(\bp_k, q_k)$ of $\vre$-approximations, we will always
assume that they are ordered so that $q_1 \leq q_2 \leq \cdots. $  
\begin{theorem}\label{thm: main epsilon}
For any norm $\norm{\cdot}$ on $\bR^d$ and any $\vre>0$ there is a
probability measure 
$\nu = \nu_{\vre-\on{approx},\norm{\cdot}}$  
on  $\XX_d \times \bR^d \times \widehat{\bZ}^n$ such that
for Lebesgue almost any $\theta\in\bR^d$, the following holds. Let
$\bw_k\in\bZ^n$ be the 
sequence of $\vre$-approximations of $\theta$ with  
respect to the norm $\norm{\cdot}$. Then the sequence
\begin{equation}\label{eq: seq eps}
  \left(\br{\pi^{\bw_k}_{\bR^d}(\bZ^n)}, \on{disp}(\theta, \bw_k),
    \bw_k\right)_{k \in \N} \in \XX_d \times \bR^d \times \widehat{\bZ}^n
  \end{equation}
 equidistributes with respect to $\nu$. The measure $\nu$ has the
 following properties: 
\begin{enumerate}
\item The measure $\nu$ is a product $\nu = \nu^{(\XX_d)}\times
  \nu^{(\R^d)} \times \nu^{(f)}$ 
  where $\nu^{(\XX_d)}\in \cP(\XX_d), \ \nu^{(\R^d)}\in
  \cP(\R^d)$ and  $\nu^{(f)}\in\cP(\widehat{\bZ}^n)$.
  \item The measure $\nu^{(\XX_d)}$ is $m_{\XX_d}$ and the measure
    $\nu^{(f)}$ is $m_{ \widehat{ 
\Z}^n_{\mathrm{prim}}}$ (in
particular, these measures do not depend on the 
choice of $\vre$ or of the norm). 

\item The measure $\nu^{(\R^d)}$ is the normalized restriction 
  of the Lebesgue measure on $\bR^d$, to the ball of radius
  $\vre$
  around the origin with respect to the norm $\| \cdot
  \|$. 
\end{enumerate}
Furthermore, each of the
coordinate sequences
\begin{equation} \label{eq: coordinate sequences eps}
  \left(\br{\pi^{\bw_k}_{\bR^d}(\bZ^n)}
  \right)
  \subset
 \XX_d, \ \ \ \left(\on{disp}(\theta, 
   \bw_k)\right)
 \subset \R^d, \ \ \  (
 \bw_k)
 \subset
 \widehat{\bZ}^n
 \end{equation}
 equidistributes in its respective space, with respect to the
 pushforward of $\nu^{(\crly{X}_d)}, \nu^{(\bR^d)}, \nu^{(f)}$ respectively.

\end{theorem}

\begin{remark}
A comparison of the two statements reveals that for the measure $\nu$ arising in
Theorem \ref{thm: main epsilon} we have a somewhat simpler description
than for the measure $\mu$ in Theorem \ref{thm: main best}. 
In fact, as will be seen in \S \ref{sec:
  properties measures}, $\mu$ is absolutely continuous with respect to $\nu$, and the density, which will be described explicitly,
  is not a product. 
\end{remark}

\begin{remark}\label{remark: second comparison}
In the case of best approximations, our proof gives more information
on the set of full measure in $\R^d$ for which the conclusions
hold. For instance, as was pointed out to us by Yiftach Dayan, using
\cite{SimmonsWeiss} one
obtains the equidistribution of the first two components of \eqref{eq:
  seq best} with respect to the measures described in 
Theorem \ref{thm: main best}, for a.e.\;$\theta$, with respect to the natural
measure on a self-similar fractal such as Cantor's middle thirds set
or the Koch snowflake. See Remark \ref{remark: for Cantor set}. 
  \end{remark}
Various special cases of these results, dealing with the individual
coordinate sequences \eqref{eq: coordinate sequences} and \eqref{eq:
  coordinate sequences eps}, and mostly for 
$d=1$, were proved in
prior work. Even in case $d=1$, the joint equidistribution of these
sequences is new. We will survey these results in \S \ref{sec: prior
  work}. 
Furthermore, in the sequel we will state a refinement (see Theorem
\ref{thm: refinement Horesh bundle}) where
equidistribution will take place in a torus bundle  
over the product 
$\XX_d \times \bR^d \times \widehat{\bZ}^n$.


Theorems \ref{thm: main best} and \ref{thm: main epsilon} give
information about typical vectors $\theta\in \bR^d$ but as is often
the case, they say nothing about concrete vectors. The
following result deals with the asymptotic statistical properties of
approximation vectors of certain algebraic vectors $\theta$. It shows
that there are  limit laws governing the approximations but that they
are actually quite different from the ones appearing in Theorems
\ref{thm: main best} and \ref{thm: 
  main epsilon}. To the best of our knowledge there are no prior
results of this type.

\begin{theorem}\label{thm: number field}
  Let $\vec{\alpha} \in \R^d$ be a vector of the form $\vec{\alpha} = (\alpha_1,
  \ldots, \alpha_d)$, where $\bK \df \spa_\Q(1, \alpha_1, \ldots, 
\alpha_{d})$ is a totally real number field of degree $n \geq 3$.
Let $\| 
\cdot \|$ be a norm on $\R^d$.
Let
\begin{equation}\label{eq: def epsilon zero}
\index{E@$\vre_0$ -- smallest $\vre$ relevant for Case II, $\vre$-approximations}
  \vre_0 \defi \inf\set{\vre>0 : \textrm{there are infinitely many
      $\vre$-approximations to $\vec{\al}$}}.
  \end{equation}
 Let $\vre>\vre_0$. 
Then there are measures 
$\mu^{(\vec{\alpha})} = \mu^{(\vec{\alpha})}_{\mathrm{best}, \|\cdot\|}$ and
$\nu^{(\vec{\alpha})} = \nu^{(\vec{\alpha})}_{\vre-\mathrm{approx}, \|\cdot\|}$
on $\XX_d \times \bR^d \times \widehat{\bZ}^n$, 
such that the following hold. Let $(\bv_k)_{k \in \N}$ and
$(\bw_k)_{k \in \N}$ denote
respectively the sequence of best approximations and
$\vre$-approximations of $\vec{\alpha}$, with respect to
$\|\cdot\|$. Then:
\begin{itemize}
\item
  The sequence \eqref{eq: seq best} equidistributes with respect to
$\mu^{(\vec{\alpha})}$, provided 
  \begin{equation}\label{eq: condition on norm}
\text{    the norm on } \R^d \text{ is either the Euclidean norm or
  the sup-norm}. 
    \end{equation}
 \item
  For any norm, the sequence \eqref{eq: seq eps} equidistributes with
  respect to $\nu^{(\vec{\alpha})}$. 
\end{itemize}
Furthermore, the supports of the projections of $\mu^{(\vec{\alpha})}$ and
$\nu^{(\vec{\alpha})}$ to $\XX_d, \R^d, \widehat{
  \Z}^n$ are null sets with respect to $m_{\XX_d}, m_{\bR^d},
m_{\widehat{\bZ}^n_{\on{prim}}}$ respectively, and  
in particular, they are singular with respect to the measures
appearing in Theorems \ref{thm: main best} and \ref{thm: main epsilon}. 
\end{theorem}

\begin{remarks}\label{remark: about number fields}
  \begin{enumerate}
  \item
   A version of Theorem \ref{thm: number field} is also true in dimension $d=1$ but
requires slightly different methods (see Remark~\ref{remark: failure for d=1} for details 
on what can go wrong in our proof). 
 Also, with our method one can also treat vectors $\vec{\alpha}$ for which the
 field $\bK$ is real but not
 totally real. However this necessitates some additional arguments and
 some additional conditions on the norm. We hope to return to these
 topics in future work.
    \item
Let $\vec{\alpha} = (\alpha_1, \ldots, \alpha_d)$ be as in Theorem \ref{thm:
  number field}. It is well-known that $\vec{\alpha}$
is badly approximable (see e.g.\,\cite{Schmidt_survey}). That is,
\begin{equation}\label{eq: that is}
    \inf \{\norm{\disp(\vec{\alpha}, \bw)}
    : \bw  \in \Z^d \times \N\}>0 \end{equation}
  (where $\disp(\vec{\alpha}, \bw)$ is as in \eqref{eq: def displacement vector}). 
  This shows that if limit measures $\mu^{(\al)}, \nu^{(\al)}$ as in
  Theorem \ref{thm: number field} exist, then  
  their projection to $\bR^d$ is bounded away from $0$. This is
  already quite different from 
  the typical behavior described in Theorems \ref{thm: main best}, \ref{thm: main epsilon}. 
\item \label{item: for SZ}
  The measures $ \mu^{(\vec{\alpha})}_{\mathrm{best}, \|\cdot\|}$ and
$\nu^{(\vec{\alpha})}_{\vre-\mathrm{approx}, \|\cdot\|}$ admit an
explicit description, see \S \ref{sec: properties measures}. 
Although they are very different from the measures $
\mu_{\mathrm{best}, \|\cdot\|}$ and 
$\nu_{\vre-\mathrm{approx}, \|\cdot\|}$ appearing in Theorem \ref{thm:
main best}
and \ref{thm: main epsilon}, in the recent  paper \cite{Uri_cheng}, the first author
and Zheng exhibit explicit choices of sequences of vectors $\vec{\alpha}_j$,
corresponding to totally real number fields, for which $
\mu^{(\vec{\alpha}_j)}_{\mathrm{best}, \|\cdot\|}\longrightarrow_{j \to \infty}
\mu_{\mathrm{best}, \|\cdot\|}.$ 

\item
  Several authors (see \cite[\S2.3]{Chevallier} and references
  therein) have shown that in some cases of vectors whose coordinates
  generate cubic fields which
  are not totally real, best approximation
  denominators vectors are periodic, or satisfy higher order linear
  recurrences. 
  We show (see Proposition \ref{prop: c and P}) that for any norm, and any
    $\vec{\alpha}$ as in Theorem \ref{thm: number field}, for both best approximations
    and $\vre$-approximations, the sequence of denominators
    $(\log(q_k ))_{k \in \N}$ is strongly asymptotic to a
    one-dimensional cut-and-project set, and the sequence of
    approximation vectors is strongly asymptotic to a generalized 
    cut-and-project set (see \S \ref{subsec: c and p structure} for the definitions of these
    terms).   
\end{enumerate}
\end{remarks}

\ignore{
 2. The measures  $ \mu^{(\vec{\alpha})}_{\mathrm{best}, \|\cdot\|}$ and
$\nu^{(\vec{\alpha})}_{\vre-\mathrm{approx}, \|\cdot\|}$ appearing in
Theorem \ref{thm: number field} depend only on the $\Z$-module in
$\bK$, generated by $1, \alpha_1, \ldots, \alpha_d$; that is, if there
is a matrix $(a_{ij}) \in \GL_d(\Z)$ such that $\beta_i = \sum_j
a_{ij} \alpha_j$,  then for $\vec{\beta} = (\beta_1, \ldots, \beta_d)$
we have
$ \mu^{(\vec{\alpha})}_{\mathrm{best}, \|\cdot\|} =
\mu^{(\vec{\beta})}_{\mathrm{best}, \|\cdot\|}$  and
$\nu^{(\vec{\alpha})}_{\vre-\mathrm{approx}, \|\cdot\|} =
\nu^{(\vec{\beta})}_{\vre-\mathrm{approx}, \|\cdot\|}$. \red{make sure
  this is okay and if so, adapt the statements so that this is proved.}
}

\subsection{Outline of method and structure of the paper}
Our method
closely follows that used by \cite{Cheung_Chevallier}, which in turn was inspired by 
\cite{Arnoux_Nogueira}. We consider a space $X$, a flow $a_t \curvearrowright X$,
 and a subset $\mathcal{S} \subset X$
intersecting all orbits along a discrete infinite countable set of
times. The set $\mathcal{S}$ is called a {\em Poincar\'e section} or
{\em cross-section}, and the flow induces a {\em return time function}
$$\tau: \mathcal{S} \to \R_+, \ \ \tau(x) = \min\{t>0 : a_t x \in
\mathcal{S}\},$$
and a {\em first return
  map}
$$T_{\mathcal{S}} : \mathcal{S} \to \mathcal{S}, \ \
T_{\mathcal{S}}(x) = a_{\tau(x)}x.$$
It has
  been known since classical work of Ambrose and Kakutani
  (see \cite{Ambrose, Ambrose_Kakutani} or \S \ref{section:cs}) that
  there is a bijection $\mu \ 
  \longleftrightarrow \ \mu_{\mathcal{S}}$ between $\{a_t\}$-invariant
  ergodic measures on $X$, and $T_{\mathcal{S}}$-invariant ergodic measures on 
  $\mathcal{S}$. These notions were intensively studied both in the
  setup of standard Borel spaces with Borel actions, and in the setup of smooth
  flows on manifolds, with sections which are codimension one
  submanifolds. In our setup, the space $X$ and the section 
  $\mathcal{S}$ are chosen so that they satisfy the following
  properties.
  \begin{itemize}
  \item
    The space $X$ is an adelic homogeneous space, that is, a quotient
    $X = \mathbf{G}(\mathbb{A})/\mathbf{G}(\Q) $ where $\mathbf{G}$ is a
    $\Q$-algebraic group, and 
    $\{a_t\}$ is a one-parameter real subgroup acting on $X$ by left
    translations. In particular, $X$ is a locally compact and second
    countable topological space, but is not a manifold.
  \item
The section $\mathcal{S}$ is chosen so that for each $\theta$ in
$\R^d$ there is $\wt{\Lambda}_\theta \in X$, such that visits to $\cS$ 
of the trajectory $\{a_t \wt{\Lambda}_\theta : t \geq 0 \}$ in
$X$, correspond in an explicit way to a sequence of approximations. In
particular we will choose distinct (but closely related) sections for dealing with best
approximations and $\vre$-approximations.
\item
  The observables we are interested in, like
  $[\pi^{\bv}_{\R^d}(\Z^n)], \, \disp(\theta, \bv)$, and $\bv$ (seen
  as an element of 
  $\widehat{\Z}^n_{\mathrm{prim}}$), are explicitly given by functions on $\mathcal{S}$.
\item
  For certain dynamically natural measures $\mu$ on $X$,
  the corresponding measures $\mu_{\mathcal{S}}$ on $\mathcal{S}$ can
  be described explicitly. In particular this is true for $\mu = m_X,$
  the unique $\mathbf{G}(\mathbb{A})$-invariant probability measure on $X$.  
\end{itemize}

The specific cross-section  we work with in this paper involves lattices which
contain a primitive vector whose $n$th coordinate is 1. We remark that
in \cite{Cheung_Chevallier} a different cross-section is used.

Recall that for a flow $\{a_t\}$ on a locally compact second countable
space $X$, and an invariant measure $\mu$, the trajectory of a point
$x \in X$ is {\em generic for $\mu$} \index{generic trajectories} if the
orbital averaging measures 
$\frac{1}{T} \int_0^T \delta_{a_tx} dt$ converge weak-* to $\mu$ as $T
\to \infty$. A similar definition can be given with $\mu_\mathcal{S},
T_{\mathcal{S}}$ in place of $\mu, \{a_t\}$. 
Although the relationship between $a_t$-invariant measures on $X$ and $T_\cS$-invariant measures on
$\mathcal{S}$ is well-understood, for our application we need a finer
understanding of the relationship between $\{a_t\}$-generic
trajectories on $(X, \mu)$ and $T_{\mathcal{S}}$-generic trajectories on
$(\mathcal{S}, \mu_{\mathcal{S}})$. We are not aware of any treatment
which is suitable for our purposes and we develop it in detail in this
paper.

In \S \ref{sec: lift functionals} we introduce a 
torus bundle $\crly{E}_n \to \XX_d$, which we call the {\em
  space of lift functionals}. We then  state strengthenings of Theorems
\ref{thm: main best},
\ref{thm: main epsilon}, \ref{thm: number field} in which the measures
$\mu^{(\XX_d)}, \nu^{(\XX_d)}$ are 
replaced with measures $\mu^{(\crly{E}_n)}, \nu^{(\crly{E}_n)}$
on $\crly{E}_n$. 
In \S \ref{sec: prior work} we state more detailed results and compare our
results to those of previous authors. \S \ref{sec: prelims} contains  some measure-theoretic
preliminaries, and  in \S \ref{sec: JM sections} we introduce
{\em reasonable cross sections} for a flow on a lcsc measure space
$(X, \mu)$ (here and throughout, {\em lcsc} is 
an abbreviation for locally compact second countable Hausdorff space, and all measures
on lcsc spaces are 
Borel regular Radon measures, but not necessarily probability
measures). Roughly speaking, the definition of a reasonable cross
section captures the minimal topological and measure-theoretic
structure needed in order to 
establish a link between generic points for $\{a_t\}$ and for
$T_{\mathcal{S}}$.

A further helpful notion introduced in Definition \ref{def: tempered} is that of a {\em
  tempered subset} of a reasonable cross-section. Such a subset is a
cross-section in its own right, and for tempered subsets, the
relation between generic points is clearer: a point in $X$ is generic
for $\mu$ if and only if its $\{a_t\}$-orbit intersects $\mathcal{S}$ in points which
are generic for $T_{\mathcal{S}}$. We remark that the cross-section we will
use for analyzing best approximations is a tempered subset, but the
cross-section we will use for analyzing $\vre$-approximations is
not (see Proposition~\ref{prof: not tempered}). Correspondingly, for best approximations our results will be
slightly stronger (see Remark \ref{remark: second
  comparison} and Remark \ref{remark: about number fields}\eqref{item:
for SZ}). We remark further that the section 
of \cite{Cheung_Chevallier} was also shown to be tempered.

In \S \ref{sec:lifting} we will discuss some useful abstract properties of
reasonable cross-sections; in particular, how to lift a cross-section
from a factor, and continuity properties of cross-sections
measures. Taken together, sections \S\S \ref{sec:
  prelims}-\ref{sec:lifting}  constitute our contribution to the
abstract theory of cross-sections in lcsc spaces. 

In the subsequent sections we will apply this abstract theory to the
spaces which are relevant for us. 
  In \S \ref{sec: spaces and measures} and \S \ref{sec: applications}
  we will introduce the specific real and adelic spaces we will 
  work with, and the flows and cross-sections we will use. Namely, for Theorems
  \ref{thm: main best} and \ref{thm: main epsilon} we will work with
  $\mu = m_{\XXnA}$, the natural measure on the adelic space
  $\SL_n(\bA)/\SL_n(\Q)$, and for Theorem \ref{thm: number field} the
  measure $\mu$ will be a homogeneous measure on an adelic torus-orbit. A
  considerable part of the paper, comprising \S \ref{sec:
    applications}-\S \ref{sec: special subsets}, will
  be devoted to checking that the cross-sections we work with are
  reasonable, and special subsets corresponding to best approximations
  and $\vre$-approximations, are Jordan measurable. In \S \ref{sec: interpreting} we will 
  explain how certain observables associated with approximation
  vectors, can be read off from intersection times with the
  section. In \S \ref{sec: properties measures} we will analyze the
  cross-section measures in detail. In \S \ref{sec: weak stable invariance} we will argue
  that for typical $\theta$, the trajectories $\{a_t \wt{\Lam}_\theta\}$
  relevant for the approximation question, give rise to generic
  intersection with the cross-section. In \S \ref{sec: concluding
    proofs} we will combine all of these ingredients 
 to conclude the proofs of the main results.

  \subsection{Acknowledgements}
  We are grateful to Anish Ghosh and Gaurav Aggarwal for useful
  comments and spotting some inaccuracies  
  in an earlier draft. 
  The authors are also grateful to Alon Agin,
  Chen Frenkel,
  Rishi Kumar, 
  Vika Rudykh, Rotem Yaari and Yuvel Yifrach for helpful comments that
  helped improve the exposition. The authors are very grateful to
  several anonymous referees for detailed and insightful remarks.
  The authors gratefully acknowledge support of the following grants:
  ISF 2019/19, ISF 871/17, ISF-NSFC 3739/21, BSF 2016256. This work
  has received funding from the 
    European Research Council (ERC) under the European Union's Horizon
    2020 Research and Innovation Program, Grant agreement no. 754475. 

  \section{A strengthening: the bundle $\crly{E}_n$ 
    of lift
  functionals over projected lattices} \label{sec: lift
  functionals}
Let $\rho_{\crly{X}_d}$ 
be the map appearing in the first coordinate of \eqref{eq: seq
  best}, that is
$$
\rho_{\crly{X}_d} (\bv ) \df \left[\pi^{\bv}_{\bR^d}(\Z^n) \right] .
$$
In this section we will introduce a probability space
$(\crly{E}_n, m_{\crly{E}_n})$  and maps $\pi_{\XX_d},
\rho_{\crly{E}_n}$ that fit in the following commuting diagram: 
\begin{equation}\label{eq:diagram1159}
\xymatrix{
 &\crly{E}_n \ar[d]^{\pi_{\XX_d}}\\
\bZ^n_{\on{prim}} \ar[r]_{\rho_{\crly{X}_d}}\ar[ru]^{\rho_{\crly{E}_n}} &\crly{X}_d
}
\end{equation}


We first define the spaces and maps group-theoretically, and then use
them to state 
strengthenings of our results. We will then discuss the geometric
information encoded in the space $\crly{E}_n$.

Let \index{E@$\crly{E}_n \cong H\Z^n$ -- bundle of lift functionals}
\begin{equation}\label{eq: all functionals}
  \crly{E}_n \df \set{\Lam \in \crly{X}_n : \mb{e}_n\in
    \Lam_{\on{prim}}}
  \end{equation}
be the set of lattices containing $\mathbf{e}_n$ as a primitive
vector, and let \index{H@$H$ -- the non-expanding group}
\begin{equation}\label{eq: non-expanding group}
  H \df \set{\smallmat{A& 0\\ \mb{x} &1}\in\SL_n(\bR) :A\in\SL_d(\bR), \ 
    \mb{x}^{\on{t}}\in \bR^d}.
  \end{equation}
Then the lattice $\Z^n$ is contained in $\crly{E}_n$, the group $H$
acts transitively on $\crly{E}_n$, and $H(\Z)$ is the stabilizer of
$\Z^n$ for this action. Thus we may identify
$\crly{E}_n\simeq H/H(\bZ)$. Since $H(\bZ)$ is a lattice in 
$H$, there is a unique $H$-invariant probability measure on
$\crly{E}_n$ which we denote  \index{M@$m_{\crly{E}_n}$ -- the natural
measure on $\crly{E}_n$}
$m_{\crly{E}_n}$. We define \index{P@$\pi_{\XX_d}$ -- the projection
  $\crly{E}_n\to \crly{X}_d$} 
\begin{equation}\label{eq: def pi XXd}
  \pi_{\XX_d}:\crly{E}_n\to \crly{X}_d, \ \ \ \  \ \pi_{\XX_d}(\Lam) \df
  \pi_{\R^d}(\Lam);
  \end{equation}
that is, the projection of $\Lam$ along $\mathbf{e}_n$. The condition
$\mathbf{e}_n \in \Lam_{\prim}$ ensures that the image of this map is
indeed in $\XX_d$. It is clear that $\pi_{\XX_d}$ intertwines the
$\SL_d(\bR)$-actions on $\crly{E}_n$ and $\crly{X}_d$ (where we view
$\SL_d(\bR)$ as embedded in $H$ by taking $\bx=0$ in \eqref{eq:
  non-expanding group}). From the uniqueness 
of invariant probability measures for transitive actions
we obtain
\begin{equation}\label{eq: measures pushed}
  \left( \pi_{\XX_d}\right)_* m_{\crly{E}_n} = m_{\XX_d}.
\end{equation}

We now define the map
$\rho_{\crly{E}_n}:\bZ^n_{\prim}\to \crly{E}_n$. For 
$v\in \bR^d$ and $t \in \R$, let \index{U@$u(v)$ -- an element of the expanding group $U$}
\begin{equation}\label{eq: def u(v)}
  u(v) \df \smallmat{I_d &v\\ 0 &1}\in \SL_n(\bR)
\end{equation}
and \index{A@$a_t = \diag{e^t,
  \ldots, e^t, e^{-dt}}$ -- one-parameter diagonal group} 
\begin{equation}\label{eq: def at}
  a_t \df \diag{e^t,\dots, e^t, e^{-dt}}\in \SL_n(\bR).
\end{equation}
Given $\mb{v} = (v_1, \ldots, v_n)^{\on{t}} \in \bZ^n_{\on{prim}}$, we set
$$v \df - \frac{\pi_{\R^d}(\bv)}{v_n} = -\frac{1}{v_n}(v_1,\dots,
v_d)^{\on{t}} \ \ \text{ and } t \df 
\frac{1}{d}\log \av{v_n},$$
and define \index{R@$\rho_{\crly{E}_n}$ -- the map from $\Z^n$ to the
  bundle of lift functionals}
\begin{equation}
  \label{eq: def rho En}
  \rho_{\crly{E}_n}(\mb{v}) \df a_t u(v) \bZ^n.
  \end{equation}
That is, we use $u(v)$ and $a_t$ to deform $\bZ^n$ to a lattice in
$\crly{E}_n$ by moving $\mb{v}$ to $\mb{e}_n$, using maps which do not
change the homothety class of the 
projection to the horizontal space. Recall that $[\cdot ]$ is our
notation for the covolume one representative of the homothety equivalence class, and note that $a_t$ acts
by dilations on $\R^d$ and $u(v)$ acts trivially on $\R^d$. Thus, this
choice ensures that 
$$
\left[\pi_{\bR^d}^{\mb{v}}(\bZ^n) \right] = \left[ a_t u(v)
  \pi_{\bR^d}^{\mb{v}}(\bZ^n)\right] =
\left[\pi_{\bR^d}(a_tu(v)\bZ^n) \right]=
\left[\pi_{\XX_d}(\rho_{\crly{E}_n}(\mb{v}))\right]. $$
This is  the commutativity of the diagram \eqref{eq:diagram1159}.

With this notation in hand we now state a strengthening of
Theorems \ref{thm: main best}, \ref{thm: main epsilon} and \ref{thm:
  number field}:

\begin{theorem}\label{thm: refinement Horesh bundle}
  Let $\| \cdot \|$ be a norm on $\R^d$, let $\vre>0$, and let
  $\vec{\alpha}$ be as in Theorem \ref{thm: number field}.
  Then there are measures
  $$\mu^{(\mathbf{e}_n)}, 
  \ \mu^{(\mathbf{e}_n, \vec{\alpha})},
  \ \nu^{(\mathbf{e}_n)}
   \text{ and }
  \nu^{(\mathbf{e}_n, \vec{\alpha})}
  $$
  on $\crly{E}_n \times \bR^d \times 
\widehat{\bZ}^n $ such that, denoting by $(\bv_k)_{k \in \N}, \,
  (\bw_k)_{k \in \N}$ the sequence of best approximations and
  $\vre$-approximations of $\theta \in \R^d$, the sequences 
\begin{equation}\label{eq: seq bundle best}
  \left(\rho_{\crly{E}_n}(\bv_k),  \on{disp}(\theta, \mb{v}_k), 
 \mb{v}_k\right)_{k \in \N} \in \crly{E}_n \times \bR^d \times 
\widehat{\bZ}^n 
\end{equation}
and 
\begin{equation}\label{eq: seq bundle eps}
  \left( \rho_{\crly{E}_n}(\bw_k), \on{disp}(\theta, \bw_k), 
 \bw_k\right)_{k \in \N} \in \crly{E}_n \times \bR^d \times 
\widehat{\bZ}^n 
\end{equation}
 equidistribute with respect to $\mu^{(\mathbf{e}_n)}$ 
 and $\nu^{(\mathbf{e}_n)}$
   respectively for Lebesgue a.e. $\theta$, and with respect to
   $\mu^{(\mathbf{e}_n, \vec{\alpha})}$ (provided \eqref{eq: condition
     on norm} holds) and 
 $\nu^{(\mathbf{e}_n, \vec{\alpha})}$ (provided $\vre > \vre_0$, for
 $\vre_0$ as in \eqref{eq: def epsilon zero}) for $\theta =
 \vec{\alpha}$. Furthermore, the properties of these measures, 
   listed in Theorems \ref{thm: main best}, \ref{thm: main
     epsilon} and \ref{thm: number field}, remain valid, provided we
   replace everywhere $\XX_d$ with 
   $\crly{E}_n$. 


  \end{theorem}

  We now discuss the additional information encoded by
  $\rho_{\crly{E}_n}(\mb{v})$, besides the information encoded by
  $\rho_{\crly{X}_d}(\mb{v})$. 
  Let $\bv \in \R^n \sm \R^d$ and suppose $\Lambda \in \XX_n$ contains
$\bv$ as a primitive vector (for example $\Lam \in
\crly{E}_n$ and $\bv = \mathbf{e}_n$). 
We will introduce  an additional geometric invariant associated with
$\bv$ and $\Lam$, which we will call the {\em
  lift functional}. \index{lift functional} It records 
  the information required to reconstruct $\Lambda$ from its
  projection onto the hyperplane $\bR^d$ along $\bv$.
Let $\Lambda'' \df \pi_{\R^d}^{\bv}(\Lam)$ and let $\Lam'$ be the
rescaling of $\Lam''$ which has covolume one. That is, if $v_n$ is the
$n$th coefficient of $\bv$, then
\begin{equation}\label{eq: def Lam one}
  \Lam'  = |v_n|^{1/d}\pi_{\bR^d}^{\bv} (\Lam)\in \XX_d.
\end{equation}
Note that $\Lam' = \pi_{\XX_d}(\Lam)$ if $\Lam \in \crly{E}_n$ and $\bv =
\mathbf{e}_n$. 

We claim that there is a linear
  functional $f: 
  \R^d \to \R$ such that
\begin{equation}\label{eq: lift functional} 
 \forall \bw \in \Lambda \ \exists k \in \Z \text{ such that } 
  \bw = \pi^\bv_{\R^d}(\bw) + \left(f(\pi^\bv_{\R^d}(\bw) ) + k \right) \bv.
\end{equation} 
To see this, complete $\bv$
  to a $\Z$-basis $\bv, \bv_1, \ldots, \bv_d$ of $\Lambda$, so that
  $\bw_i \df \pi^{\bv}_{\R^d}(\bv_i) \ (i=1, \ldots, d)$ is a $\Z$-basis
  of $\Lambda''$. For each $i$ there is $c_i \in \R$ such that $\bv_i =
  \bw_i + c_i \bv$, and we can define $f \in (\bR^d)^*$ via the requirement 
$$f(\bw_i) =c_i \ \ \ (i=1, \ldots, d).$$
This construction depends on the choice of the vectors $\bv_i$, but
any two functionals $f_1, f_2$ satisfying \eqref{eq: lift
  functional} will satisfy that $f_1(\bu) - f_2(\bu) \in \Z$
for any $\bu \in \Lambda''$, 
that is they differ
by an element of the dual lattice $(\Lambda'')^*$. It will be more
convenient to work with lattices of covolume 1, so we will replace
$\Lam''$ with $\Lam'$ and $f$ with
\begin{equation}\label{eq: def phi one}
  f'(x) \df f\left(
  |v_n|^{-1/d}x \right),
\end{equation}
so that $f'$ is well-defined as  a class in the
torus $\bT_{\Lam'} $, where for $\Delta \in \XX_d$ we denote
\index{T@$\bT_{\Delta} =(\R^d)^*/\Delta^*$ torus of lift functionals
  over a $\Delta \in \XX_d$}
\begin{equation}\label{eq: torus of functionals} 
  \bT_{\Delta} \df \left(\R^d \right)^*/\Delta^*.
  \end{equation}
 This class $[f']$ is the lift
functional. 

Consider the set
$$\crly{E} \df \left\{(\Lam', f) : \Lam'\in \XX_d , \ f\in  \bT_{\Lam'}
\right\},$$
which is a torus bundle over $\XX_d$. Then a pair $(\Lam',
f)$ constructed as above by projecting $\Lam$ along $\bv \in
\Lam_{\prim} \sm \R^d$ is an
element of $\crly{E}$. Moreover, using \eqref{eq: lift functional}, we
can recover $\Lam$ from $\Lam', 
f$ and $\bv$, as follows: \index{L@$\Lam(\Lam', f, \bv)$ -- the
  lattice reconstructed from the projected lattice, lift functional,
  and projection vector}
\begin{equation}\label{eq: recovering from lift}
\Lam = \Lam(\Lam', f, \bv) = \set{\av{v_n}^{-\frac{1}{d}}\bx + (f(\bx)+k)\bv: \bx \in \Lam', \, k \in \Z}.
\end{equation}
 Thus for each fixed $\bv \in \R^n \sm \R^d$ we have an
identification of $\crly{E}$ with the set of lattices in $\XX_n$ which
contain $\bv$ as a primitive vector. In particular, choosing $\bv =
\mathbf{e}_n$ we obtain $\crly{E}_n$ as the image of the map
$\Lam(\cdot, \cdot, \mathbf{e}_n)$ as in \eqref{eq: recovering from
  lift}.

In order to make the connection between $\crly{E}_n = H\Z^n$
and lift functionals more concrete, we write any element of $ h \in H$  in the form 
$h = \smallmat{ A & \mathbf{0} \\ \bx & 1 },$
as in \eqref{eq: non-expanding group}.
Writing $h = (A, \bx)$ we see that $H$ is a semi-direct product
$\SL_d(\R) \ltimes \R^d$ 
where $\SL_d(\R)$ acts on $\R^d$ by right-multiplication on row
vectors. Furthermore, for $h_i = (g_i,
\bx_i)$ we have $h_1\Z^n = h_2 \Z^n$ if and only if
 $g_1  \in g_2\SL_d(\Z)$  and $\bx_1g_1^{-1}  \in
\bx_2g_2^{-1} + \Z^dg_1^{-1}$. In particular, the lattice
$\Lambda' = g_i \Z^d$ and 
the equivalence class of the functional $f(\bw)
= \bx_i g_i^{-1} \bw,$ seen as an element of $\mathbb{T}_{\Lambda'}$,
are well-defined 
independently of $i=1,2$. They represent respectively the projected lattice
$\pi^{\bv}_{\R^d}(h_i\Z^n)$, and the lift
functional.

\section{Consequences, and prior work}\label{sec: prior work}
%
%
Fix a norm $\| \cdot \|$ on $\R^d$ and fix $\vre>0$. 
Let $\mu, \, \nu$ be as in Theorem \ref{thm: main best} and \ref{thm:
  main epsilon} respectively, and let \index{N@$\nu^{(\XX_d)}, \nu^{(\R^d)}, \nu^{(f)}$ -- projected limit
  measures, $\vre$-approximations} \index{M@$\mu^{(\XX_d)}, \,\mu^{(\R^d)},
  \mu^{(f)}$ -- projected limit measures, best approximations}
$\mu^{(\XX_d)}, \,\mu^{(\R^d)},
\mu^{(f)}, \nu^{(\XX_d)}, \nu^{(\R^d)}, \nu^{(f)}$ 
be the projections of these
measures on the respective factors.  We 
discuss prior work about equidistribution on these spaces.
Throughout this section we fix $\theta
\in \R^d$ for which the conclusions of Theorems \ref{thm: main best} and
\ref{thm: main epsilon} are valid, and let $(\bv_k), \, (\bw_k)$
denote the sequence of best approximations and $\vre$-approximations
of $\theta$.

\subsection{Growth of denominators, the Khinchin-L\'evy constant, and a
  `Khinchin-L\'evy distribution'}
The \index{Khinchin-L\'evy constant} asymptotic rate of growth of the
denominators $q_k$ (i.e., the 
last coordinate of the vector $\bv_k$ in case of best approximations,
and of the vector $\bw_k$ in case of $\vre$-approximations) has been a
topic of detailed study. 
Our analysis allows us to deduce the following result, which will be derived at the end of \S\ref{sec: concluding proofs}.
\begin{corollary}\label{cor: KL distribution}
Let $\| \cdot \|$ be a norm on $\R^d$ and let $\vre>0$. Then there 
are probability measures $\lambda^{(KL)}_{\| \cdot \|, \mathrm{best}} , \,
\lambda^{(KL)}_{\| \cdot \|, \vre}$ on $\R_+$ with finite expectations 
$\ga^{(KL)}_{\norm{\cdot}, \on{best}}, \ga^{(KL)}_{\norm{\cdot},\vre}$
respectively, such that for a.e.\;$\theta\in \R$,
the following holds. Let $(\bv_k)_{k \in \N}, \, (\bw_k)_{k \in \N}$
denote respectively the best approximation and $\vre$-approximation
sequence, and let $q_k^{(\mathrm{best})}, q_k^{(\vre)}$ denote the
corresponding denominators. 
\begin{enumerate}
\item\label{item: Cor KL 1} 
The following holds
\begin{align*}
\lim_k\pa{ q_k^{(\on{best})}}^{1/k} &=  \exp\pa{d\cdot \ga^{(KL)}_{\norm{\cdot},\on{best}}}\\
\lim_k\pa{ q_k^{(\vre)}}^{1/k}&= \exp\pa{d\cdot \ga^{(KL)}_{\norm{\cdot},\vre}}.
\end{align*}
\item\label{item: Cor KL 2} 
For any $k \geq 1$, let 
\begin{equation}\label{eq: seq KL}
s_k^{(\mathrm{best})} \df\frac{1}{d} \log
\left(\frac{q^{(\mathrm{best})}_{k+1}}{q^{(\mathrm{best})}_{k}}
  \right) \ \ \text{ and } \ \ s_k^{(\vre)} \df \frac{1}{d}\log
\left(\frac{q^{(\vre)}_{k+1}}{q^{(\vre)}_{k}}
  \right). 
  \end{equation}
  Then, the uniform probability measures on the sets
  $$\left\{s_k^{(\mathrm{best})} : 1
  \leq k \leq N \right\}, \ \ \ \ \  \left\{s_k^{(\vre)} : 1
  \leq k \leq N \right\}$$
converge respectively to $\lambda^{(KL)}_{\|
  \cdot \|, \mathrm{best}} , \, 
\lambda^{(KL)}_{\| \cdot \|, \vre}$, in the weak-* topology as $N \to
\infty$.
\end{enumerate}
Finally, 
$\ga^{(KL)}_{\| \cdot \|, \vre}=\frac{\zeta(n)}{d\cdot V_{d, \|\cdot\|} \, \vre^d}$, where $V_{d,
  \|\cdot\|}$ is the Lebesgue measure of the unit ball of the norm, in $\R^d$. 
  \end{corollary}

Some remarks are in order to put the results in the above Corollary in context.
For $d=1$, in the case of best
approximations and where there is only one norm, Khinchin \cite{Khinchin} showed in 1936 that the limit $\lim_{k
 \to \infty}  q_k^{1/k}$ exists, and in the same year, L\'evy
\cite{Levy} proved it is equal to $\pi^2/12\log 2$ (known today as the Khinchine-L\'evy constant).
In 1986 Jager proved equidistribution as in Corollary~\ref{cor: KL
  distribution}\eqref{item: Cor KL 2},  
and gave an explicit formula for
the limit measure (he actually also established a joint
equidistribution result for the 
sequence of 
pairs $\left(t_k^{(\mathrm{best})}, \disp(\theta, \bv_k)\right)_{k \in
  N}$ , see
\cite[Lemma 5.3.11]{Dajani_Kraaikamp}).  For best approximations and arbitrary $d>1$, 
for the Euclidean norm, the content of the Corollary 
  was recently proved by Cheung and
Chevallier~\cite{Cheung_Chevallier}. Our work extends this to arbitrary
norms, answering \cite[Question 1]{Cheung_Chevallier};
the question of 
Cheung and Chevallier was also answered in work of Aggarwal and Ghosh
\cite{Aggarwal_Ghosh_KL}, 
in a more general setting. In a subsequent work
with Xieu \cite{Xieu}, the constant $\ga^{(KL)}_{\norm{\cdot},
  \on{best}}$ was evaluated in the case $d=2$. 

Lagarias
\cite{Lagarias2} constructed examples showing irregular growth of the
sequence $q_k^{1/k}$; in particular showing that the a.e.\;behavior can
fail dramatically for some special choices of $\theta$. A different way to
measure the growth of approximations is due to Borel and Bernstein in
dimension $d=1$, 
and was extended by Chevallier to the case $d>1$. See
\cite{Chevallier} for more details on this, and for a survey of more results on best
approximations. 

 As noted in \cite[Question 1]{Cheung_Chevallier}, it would be
 interesting to compute  $\ga^{(KL)}_{\|\cdot\|,
 \mathrm{best}}$ for different dimensions and different norms. It would
 also be interesting to obtain additional statistical informations
 about these distributions, such as higher moments and tail bounds.

\subsection{Equidistribution of (shapes of) projections in
  $\XX_d$}\label{subsec: Schmidt shapes}
Our theorems assert that the sequence of projected lattices 
 $\left(\br{\pi^{\bw_k}_{\bR^d}(\bZ^n)}
  \right)_{k \in \N}$ equidistributes with respect to the natural measure $\nu^{(\XX_d)}
  = m_{\XX_d}$ on $\XX_d$, and the sequence $ \left(\br{\pi^{\bv_k}_{\bR^d}(\bZ^n)}
  \right)_{k \in \N}$ also equidistributes  with respect to a measure $\mu^{(\XX_d)}$ on
  $\XX_d$, and $\mu^{(\XX_d)} \ll m_{\XX_d}$. We stress that $\mu^{(\XX_d)} \neq
  m_{\XX_d}$. We will give an explicit description of $\mu^{(\XX_d)}$ in \S
 \ref{sec: properties measures}. In particular we will show that
 $\mu^{(\XX_d)}$ gives less weight to 
  lattices with short vectors than $m_{\XX_d}$. 

  The first equidistribution result for projected lattices appears in 
  beautiful papers  of Roelcke  \cite{Roelcke} and 
  Schmidt \cite{Schmidt_shapes}. Let
  $$\SXX_d \df \SO_d(\R)
  \backslash \SL_d(\R) / \SL_d(\Z)$$
  denote the {\em space of shapes of
    lattices}, \index{S@ $\SXX_d = \SO_d(\R)
  \backslash \SL_d(\R) / \SL_d(\Z)$-- space of shapes of lattices in
    $\R^d$} that is, the space of similarity classes of lattices in 
  $\R^d$ (two lattices are {\em similar} if they differ from each
  other by a conformal linear map, that is a composition of a dilation
  and an orthogonal transformation). This space is equipped with a natural measure
  $m_{\SXX_d}$ \index{M@$m_{\SXX_d}$ -- natural measure on $\SXX_d$}
  constructed from the Haar measure on 
  $\SL_d(\R)$. There is a natural projection map $\XX_d \to
  \SXX_d$ with compact fiber, and we can define
  $m_{\SXX_d}$ concretely as the image of $m_{\XX_d}$
  under this projection. 

  For any $\bv \in \Z^n \sm \{0\}$, the projected lattice
  $\pi^{\bv}_{\bv^\perp}(\Z^n)$ is a lattice in the $d$-dimensional
  subspace $\bv^{\perp} \subset \R^n,$ and represents an equivalence
  class in $\SXX_d$ (which we continue to denote by 
  $\pi^{\bv}_{\bv^{\perp}}(\Z^n)$). Schmidt showed that the uniform 
  measures on the finite sets $\left\{ \pi^{\bv}_{\bv^{\perp}}(\Z^n) : \bv
  \in \Z^n_{\mathrm{prim}}, \, \| \bv \| \leq 
  T\right\}$  are equidistributed with respect to $m_{\SXX_d}$, in
  the limit as $T \to \infty$\footnote{In fact Schmidt considered the
    lattices obtained by interecting $\Z^n$ with $\bv^\perp$, but the
    results are equivalent, as can be easily seen by passing to dual lattices.}. 
  There has been a recent surge of interest in equidistribution
  results for shapes of projected lattices, see \cite{Marklof-Inventiones, EMSS,
    AES2} and references therein. In these finer results, one typically considers
  projected lattices along a sparser sequence of vectors, than the one considered
  by Schmidt. 

Our results deal with the projections 
  along vectors close to a line and to the best of our knowledge these have never been considered.
   Note that we do not work in
  $\SXX_d$ but in the larger  space $\XX_d$. Indeed, in
  our case the directions of spaces 
  in which the projected lattices lie converges, and there is no 
  reason to reduce lattices modulo similarity.

  \subsection{Equidistribution in the bundle of lift functionals}
In case $d=1$ the space $\XX_d$ reduces to a point and the bundle
$\crly{E}_n$ is isomorphic to the fiber over this point, namely the one-dimensional torus
$\mathbb{T}_1=\R/\Z$. It thus makes sense to consider an analogous question to
Schmidt's result mentioned in \S \ref{subsec: Schmidt shapes}, namely
ask whether in the limit as $T \to \infty$, 
the uniform measures on the finite collection of  lift functionals
associated with vectors of $\Z^2_{\mathrm{prim}}$ of length at most
$T$, converges to the uniform measure on $\mathbb{T}_1$. Indeed, such
results were obtained by Dinaburg-Sinai \cite{Dinaburg_Sinai} and
Risager-Rudnick \cite{RR}. The 
analogous question for lift functionals arising from a sequence of
approximations to a line in $\R^2$ has not been previously considered,
as far as we know. In this very special case, our result for the
equidistribution of the sequence $\left(\rho_{\crly{E}_n}(\bv_k)
\right)_{k \in 
\N}$ takes the following form:
\begin{corollary}\label{cor: RR for lines}
  Let $\vre>0$. 
For a.e. $\theta \in \R$, let $(p_k, q_k)$ be the sequence of best
approximations, let $(p'_k, q'_k)$ be the sequence of
$\vre$-approximations, and let $(u_k, v_k )$ and $(u'_k , v'_k)$
denote
respectively
the unique solutions in $\bZ^2_{\on{prim}}$ of

%
\begin{align*}
&\det \left( \begin{matrix} p_k & u_k \\ q_k & v_k \end{matrix} \right) = 
\det \left( \begin{matrix} p'_k & u'_k \\ q'_k & v'_k \end{matrix} \right)
 =1,\\
& -\frac{1}{2}\le \frac{v_k}{q_k}, \frac{v_k'}{q_k'} < \frac{1}{2}. 
\end{align*}
 Then, the sequence 
$  (\frac{v_k'}{q_k'})_{k\in \N}
 $ is
 equidistributed in the interval $\left[ -\frac12, \frac12\right]$, with respect to
 Lebesgue measure; and the sequence
 $(\frac{v_k}{q_k})_{k\in \N}$ is equidistributed in $\left[ -\frac12,
   \frac12\right]$, with respect to an absolutely continuous measure
 whose density is given by
 $$
 F: \left[ - \frac{1}{2}, \frac12\right] \to \R, \ \ F(t) =
 \frac{1}{2\log 2} \, \left(\frac{1}{2-|t|}+\frac{1}{1+|t|} \right). 
  $$
\end{corollary}
See \cite{Dinaburg_Sinai, RR, HN} for related work and for an interpretation
of $(u_k, v_k)$ as the \index{shortest solutions to the gcd equation }
shortest solutions to the gcd equation.

\subsection{Equidistribution of displacement vectors in $\R^d$}
In case $d=1$, for best approximations, a very closely
related result is the so-called {\em
  Doeblin-Lenstra conjecture}, proved \index{Doeblin-Lenstra
  conjecture } by Bosman, Jager and Wiedijk (see 
\cite{Doeblin, BJW} for the original papers, and \cite[Chap. 5.3.2]{Dajani_Kraaikamp} and 
\cite[Chap. 4]{Iosifescu_Kraaikamp} for more detailed treatments
and related results). They showed that for a.e.\;$\theta \in \R$, with
best approximations 
$(q_k, p_k)$, the uniform measures on the sets 
$$\{|q_k(q_kx-p_k)| : n =1, \ldots, N\}$$
converge weak-* to the measure $\nu$ on
$[0,1]$ whose density function is given by
\begin{equation}
  \label{eq: cumulative Doblin}F: [0,1] \to \R, \ \ \ 
F(t ) = \left\{ \begin{matrix} \frac{1}{\log 2} & 0 \leq t \leq
    \frac{1}{2} \\ \frac{1/t-1}{\log 2}  & \frac{1}{2} <
    t \leq 1. \end{matrix} \right. 
\end{equation}
In contrast, our result, is about the uniform measures
on $$\{q_k(q_kx-p_k): n = 1, \ldots, N\}$$
(i.e., we have removed the absolute values), and we show weak-*
convergence to the measure on $[-1,1]$ whose density 
function is
$$\bar F : [-1,1] \to \R, \ \ \bar{F}(t) = \frac{1}{2} F(|t|).$$
That is, the measure $\mu'$ we consider is invariant under $t \mapsto
-t$ and projects to the measure $\mu$
 given by \eqref{eq: cumulative Doblin} under the map $t\mapsto |t|$.

As far as we are aware, in dimension $d>1,$ the distribution of the
displacement sequence 
$\disp(\theta, \bv_k)$ for best approximations has not been
studied. Similarly, 
the distribution of displacements $\disp(\theta, \bw_k)$ for
$\vre$-approximations has not been studied, even in the one
dimensional case. However, there has been some study of the
{\em direction} of the displacement of $\vre$-approximations. Namely,
let $\mathbb{S}^{d-1}  = \{\vec{x} \in \R^d: \| \vec{x}\|=1\}$ be 
  the unit sphere with respect to the given norm, let 
\index{P@$\mathrm{Proj}$ -- the projection $\R^d \to
      \bS^{d-1}$}
    \begin{equation}\label{eq: def Proj}
      \mathrm{Proj}: \R^d \to \mathbb{S}^{d-1}, 
      \end{equation}
      and consider the sequence
      \begin{equation}\label{eq: spiraling}
  \left(\mathrm{Proj}(\disp(\theta, \bw_k)) 
  \right)_{k
  \in \N} \subset \mathbb{S}^{d-1}.
\end{equation}
In \cite{Athreya_et_al_spiraling}, Athreya, Ghosh and Tseng showed that
  if $\| \cdot \|$ is the Euclidean norm, and $\vre$-approximations
  are defined without requiring $\gcd(p_1, \ldots, p_d, q)=1$, 
  then the sequence \eqref{eq:
    spiraling} is equidistributed with respect to the unique
  $\SO_d(\R)$-invariant measure on $\mathbb{S}^{d-1}$. From Theorem
  \ref{thm: main epsilon} we immediately obtain the following statements, valid
  for any norm:
  \begin{corollary} 
       For any norm $\| \cdot \|$ on $\R^d$, let $m|_B$ be the
    normalized restriction of Lebesgue measure to the unit norm-ball $\{\bu
    : \| \bu \| \leq 1\}$, and let $\nu^{(\mathbb{S}^{d-1})} = (\mathrm{Proj})_*
    (m|_B)$. Also let $\mu^{(\mathbb{S}^{d-1})} = (\mathrm{Proj})_* \mu^{(\R^d)}$, where
    $\mu^{(\R^d)}$ is as in Theorem \ref{thm: main best}. Then for
    any $\vre>0$, the
    sequence \eqref{eq: spiraling} of directions of
    $\vre$-approximations is equidistributed with respect to
    $\nu^{(\mathbb{S}^{d-1})}$, and the sequence
     \begin{equation}\label{eq: sequence of directions}
       \left(
         \mathrm{Proj} (\disp(\theta, \bv_k))
       \right)_{k
  \in \N} 
\end{equation}
of directions of best approximations, 
is equidistributed with respect to $\mu^{(\mathbb{S}^{d-1})}$.
    \end{corollary}
%
%
Note that by Theorem \ref{thm: main best}\eqref{item: 5}, when $\|
\cdot \|$ is the Euclidean norm, the measures
$\mu^{(\mathbb{S}^{d-1})}$ and $\nu^{(\mathbb{S}^{d-1})}$ coincide,
but this is not true for general norms. 

In another direction, Moschevitin \cite{Moshchevitin} has
characterized the possible sets that can be obtained as the closures
of the sequences \eqref{eq: sequence of directions}, for various
$\theta$. 

\subsection{Equidistribution in $\widehat{\Z}^n$, and congruence
  constraints} 
According to \eqref{eq: coordinate sequences} and \eqref{eq:
  coordinate sequences eps},  for Lebesgue a.e. $\theta$, the sequences
$\left(\bv_k \right)_{k \in \N}, \left(\bw_k\right)_{k \in \N}$ of
best approximations and $\vre$-approximations, considered as elements
of $\widehat{\Z}^n_{\mathrm{prim}}$, are both equidistributed with
  respect to $m_{\widehat{\Z}^n_{\mathrm{prim}}}. $ 
  In case $d=1$, such results have a long history.
  In 1962, Sz\"{u}sz \cite{Szusz}
  proved the following. For
  any $a$ and $m$, there is $c'$ such that for any $\vre>0$, for a.e.\;$\theta \in \R$,
  the sequence 
$(p_k , q_k)$ of $\vre$-approximations satisfies
  $$
\frac{1}{N} |\{1 \leq k  \leq N : q_k \equiv a \mod m \}|\to_{N\to\infty} 
c'.
$$
He also proved a similar statement for sequences of rationals arising
from more general approximation
functions. In 1982 Moeckel \cite{Moeckel} proved a dynamical result
which implies the following. For a
positive integer $m$ and integers $a,  b$, there is $c$ such that for
a.e.\;$\theta \in \R$, the sequence 
$(p_k, q_k)$ of best approximations satisfies
$$\frac{1}{N} |\{1 \leq k \leq N : q_k \equiv a\ \& \  p_k \equiv
b\mod m\}|\to_{N\to\infty} 
c.
$$
In these results, the constants $c, c'$ can be easily computed. 
See also \cite{JagerLiardet, Dajani_Kraaikamp, 
  FisherSchmidt} and references therein.

For a positive integer $m$ and a vector $\vec{a} \in (\Z/m\Z)^n$, we say
that $\vec{a}$ is {\em primitive (mod $m$)} if \index{primitive (mod
  $m$)} there are no $b, d_1,
\ldots, d_n \in \Z/m\Z$ with $b \notin(\bZ/m\bZ)^\times$ and  $a_i = bd_i$ for
all $i$. Let $N_{m,n}$ be the cardinality of the set of primitive vectors in
$(\Z/m\Z)^n$; it is not hard to verify that
$$N_{m,n} = \prod_{i=1}^k \left(p_i^{nr_i} - p_i^{n(r_i-1)} \right), \
$$
where $ m = p_1^{r_1} \cdots p_k^{r_k}$ is the prime factorization of
$m$. We generalize the above results to any dimension $d\ge 1$ in the following Corollary, which follows immediately
from Theorems~\ref{thm: main best}, \ref{thm: main epsilon}.  
\begin{corollary}\label{cor: adelic equidistribution}
  Let $\|\cdot \|$ be a norm on $\R^d$.
Let $ m$ be a 
positive integer and let $\vec{a} \in (\Z/m\Z)^n$. 
Then
for any $\vre>0$, 
a.e.\;$\theta \in \R^d$, the sequences  
$\left(\bv_k \right)_{k \in \N}, \left(\bw_k\right)_{k \in \N}$
of best approximations and $\vre$-approximations to $\theta$ both satisfy
$$
\frac{1}{N} |\{1 \leq k \leq N : \bu_k\equiv \vec{a}\mod m \}|\to_{N\to\infty}  c $$
(in both cases $\bu_k = \bv_k$ and $\bu_k= \bw_k$), where
$$c  \df \left\{\begin{matrix} \frac{1}{N_{n,m}}  & \vec{a} \text{ is
      primitive (mod $m$)} \\ 0 
  & \text{otherwise.} \end{matrix}  \right. 
$$
  \end{corollary}
As far as we know, prior to
this work, 
there were no such results for dimension $d>1$; this despite a
remark of Cassels (Math Review of \cite{Szusz}) that `it would be
interesting to extend these results to simultaneous approximation'. A
related work of Berth\'e, Nakada and Natsui \cite{Berthe_Nakada_Natsui}
establishes similar 
properties for the convergents arising
from the $d$-dimensional Jacobi-Perron algorithm. 

 \section{Preliminaries}\label{sec: prelims}
 In this section we recall some standard results. Some of them will be
 valid in a general Borel measurable setup, and for some we will
 require some topological assumptions. To distinguish these cases we
 will use the following terminology. We will have a space $X$ on which
 a one-parameter group $\{a_t\}$ acts, a
 $\sigma$-algebra\index{B@$\cB_X$ -- Borel $\sigma$-algebra on $X$} $\cB_X$
 of subsets of $X$, and a measure $\mu$ which is $a_t$-invariant. 
 \begin{itemize}
\item {\em The Borel setup} is the \index{Borel setup} one in which $(X, \cB_X)$ is a
  standard Borel space, i.e., there is a Polish 
structure (the topology induced by a complete and separable metric) on
$X$ for which $\cB_X$ is the Borel 
$\sigma$-algebra. In this setup the action map $\R \times X \to X,  \
\  (t, x) \mapsto a_tx$
is assumed to be Borel, and the measure $\mu$ is assumed to be 
$\sigma$-finite.
\item
  {\em The lcsc setup} is \index{lcsc setup} the one in which in addition, $X$ is a
  Hausdorff locally compact
  second countable (lcsc) space, $\cB_X$ is the Borel $\sigma$-algebra for
  the underlying topology, and the action map is 
  continuous. Moreover the measure $\mu$ is assumed to be Radon
  (regular and finite on compact sets).

   \end{itemize}

   The above conditions will be assumed throughout the
   paper without further mention. Additionally, we will usually
   assume that $\mu$ is finite, but will make explicit mention of this
   when we do. 
   
\subsection{Measurable cross-sections}\label{section:cs}
In this section we recall
classical definitions and results 
in the Borel setup. We start with the definition of a Borel cross-section which do not involve a measure
and continue with the definition of a $\mu$-cross-section which is relevant when a measure is present.
\begin{definition}\label{def: Borel cs}
Let $(X, \cB_X, \set{a_t})$ be as in the Borel setup.
A {\em
  Borel cross-section} is \index{Cro@ Borel cross-section} a Borel
subset $\cS$ with the following 
properties:
\begin{itemize}
\item
For any $x\in X$,
  the sets \index{Y@$\mathcal{Y}_x $ -- set of visit times to a
    cross-section} of \index{visit times} {\em visit times}
  \begin{equation}\label{eq: visit times}
    \mathcal{Y}_x \df
    \{t \in \R: a_t x \in \cS\}
    \end{equation}
  are all 
  discrete and totally unbounded; that is, for any $T>0$,
  
  {\small
  \begin{equation}\label{eq: requirements visit times} \mathcal{Y}_x \cap (T, \infty) \neq \varnothing,  \ 
  \mathcal{Y}_x \cap (-\infty, -T) \neq \varnothing, \text{ and } \ \#
  \left( \mathcal{Y}_x \cap (-T, T) \right) < \infty  .\end{equation}}
\item
  The \index{return time function} {\em return time function}
  \begin{equation}\label{eq: def return time}
\tau_{\cS} : \cS \to \R_+, \ \ \ \ \tau_{\cS}(x) \df \min \left(\mathcal{Y}_x
\cap \R_+\right)
\end{equation}
is Borel. \index{T@$\tau_{\cS}$ -- return time function}
\end{itemize}
\end{definition}
If $\cS$ is a Borel cross-section we will denote by $T_{\cS} : \cS \to \cS$ the {\em
  first return map} \index{first return map} \index{T@ $T_{\cS}$ --
  first return map to a cross-section} defined by
$$
T_{\cS} (x) \df a_{\tau_{\cS}(x)}x.
$$
\begin{definition}\label{def:cs}
Let $(X,\cB_X,\set{a_t})$ be as in the Borel setup and let $\mu$ be an
$\set{a_t}$-invariant measure. 
We will say that $\cS\in \cB_X$ {\em is a $\mu$-cross-section}
\index{M@ $\mu$-cross-section}
  if there is an $\set{a_t}$-invariant set $X_0
\in \cB_X$ such that
   \begin{itemize}
\item $\mu(X \sm X_0)=0$.
\item $\cS \cap X_0$
is a Borel cross-section for $(X_0,
\cB_{X_0})$  (where $\cB_{X_0}$ is
the restricted $\sigma$-algebra defined by
$\cB_{X_0}  \df \{A \cap X_0: A \in \cB_X\}$).
\end{itemize}
\end{definition}

If $\cS$ is a $\mu$-cross-section and $X_0$ is as in
Definition~\ref{def:cs}, then we denote $\tau_\cS: \cS\to \bR_+\cup
\set{\infty}$ 
the function 
$$\tau_\cS(x) \df 
\Big\{
\begin{array}{ll}
\tau_{\cS\cap X_0}(x) &\text{ if }x\in \cS\cap X_0,\\
\ \infty&\text{  if }x\in \cS\smallsetminus X_0.
\end{array}$$ 
We denote $T_{\cS\cap X_0}$ by $T_\cS$. It is a well defined Borel map $\cS\cap X_0\to \cS\cap X_0$, and we say that a measure on $\cS$ is $T_\cS$-invariant if $\cS\smallsetminus X_0$ is a null-set and its restriction to $\cS\cap X_0$ is $T_\cS$-invariant.

\begin{remark}
  The restricted Borel space $(X_0, \cB_{X_0})$
    appearing in Definition \ref{def:cs} is a standard Borel space,
    see \cite[Chap. 13]{Kechris}.
    
\end{remark}

For $\vre >0$ we \index{S@ $\cS_{< \vre}  $ -- sub-level set of
  return time function} let 
\begin{equation}\label{eq: sublevel sets}
  \cS_{\ge \vre} \df \set{x\in \cS: \tau_\cS(x)\ge \vre} \ \text{ and } \
  \cS_{<\vre} \df \cS \sm \cS_{\ge \vre} .
  \end{equation}
The sets $\cS_{\ge\vre}$ are an increasing collection of Borel sets
whose union is $\cS$. Given $E\in  \cB_X$ and $I\subset \bR$ we let 
\begin{equation}\label{eq: thickening} E^I \defi \set{a_tx:x\in E,
    t\in I},
  \end{equation}
and \index{T@ $E^I = \bigcup_{t \in I} a_t(E)$ -- $I$-thickening of
  $E$} let $m=m_{\bR}$ denote the Lebesgue measure on $\bR$.

The following structure theorem is due to
Ambrose and Kakutani \cite{Ambrose, Ambrose_Kakutani}, see also
\cite{Nadkarni, Wagh}: 
\begin{theorem}\label{cor:csmeasure} Let $(X, \cB_X, \{a_t\})$ be as in
  the Borel setup and let $\cS\in \cB_X$. 
Then, for any finite
$\{a_t\}$-invariant 
measures $\mu$ on $(X, \cB_X)$ for which $\cS$ is a
$\mu$-cross-section, there exists a 
$T_\cS$-invariant measures $\mu_\cS$ on
$(\cS, \cB_{\cS})$ such that  the following holds for any
Borel set $E \subset \cS$:  
\begin{enumerate}[(i)]
\item\label{1204-2} If $E\subset \cS_{\ge \vre}$, and $I$ is an
  interval of length $<\vre$ then \index{M@$\mu_{\cS}$ -- induced
    measure on $\cS$}   
$$\mu_\cS(E) = \frac{\mu(E^I)}{m(I)}.$$
In particular, for any $\vre>0$ we have $\mu_\cS(\cS_{\geq \vre})<\infty$. 
\item\label{1204-1} For any interval $I$, $\mu_\cS(E)
  \ge \frac{\mu(E^I)}{m(I) }.$ 
\item\label{1204-3} In general 
\begin{equation}\label{eq:1019}
\mu_\cS(E) = \lim_{\vre \to 0} \frac{1}{\vre} \, \mu\left(E^{(0,\vre)}\right).
\end{equation}
\item\label{1204-4} 
  We have $\mu_\cS(E)=0$ if and only if $\mu(E^\bR) = 0$.
  \item\label{item: ergodic decomposition} If $\mu = \int \eta
   \, d\Theta(\eta)$ for a measure $\Theta$ on $\mathcal{P}(X)$, where
    $\Theta$-a.e. $\eta$ is $\{a_t\}$-invariant, then $\cS$ is an
    $\eta$-cross-section for $\Theta$-a.e. $\eta$,
    and $\mu_{\cS}  =
    \int \eta_{\cS} \,
    d\Theta(\eta)  .$ Moreover, $\mu$ is
    $\{a_t\}$-ergodic if and only if 
$\mu_\cS$ is $T_\cS$-ergodic. 
\item\label{item: Kac formula}
We have  
  $$
\mu(X) = \int_{\cS} \tau_\cS \, d\mu_\cS.
$$
\item\label{item: one more}
  If $K \curvearrowright (X, \mu)$ is a group-action commuting with
  the $\{a_t\}$-action and preserving $\cS$, then the $K$-action also
  preserves $\mu_{\cS}$. 
\end{enumerate}
\end{theorem}
We will refer to item \eqref{item: Kac formula} as \index{Kac formula} the {\em Kac
  formula.} Note that \eqref{item: one more} is not mentioned in the
above references but follows immediately from \eqref{eq:1019}. 

\begin{definition}
Let $(X,\cB_X,\set{a_t})$ be as in the Borel setup and let $\mu$ be a
finite $\set{a_t}$-invariant measure. Let 
$\cS$ be a $\mu$-cross-section. Then the measure $\mu_\cS$ 
from Theorem~\ref{cor:csmeasure} is called the {\em cross-section 
measure of $\mu$.} 
\end{definition}

Note that in Theorem \ref{cor:csmeasure}, $\mu$ is finite but $\mu_{\cS}$ need not
be. However, from \eqref{item: Kac formula} one easily sees that 
if $\tau_{\cS}$ is bounded below, or more generally, if there is $c>0$
such that $\tau_{\cS}(x) \geq c$ for $\mu_\cS$-a.e. $x$, then
$\mu_{\cS}$ is finite. 
\begin{remark}
In Theorem \ref{cor:csmeasure} we defined a map $\mu\mapsto
\mu_\cS$. Under suitable conditions this map is injective and  
its image has an explicit description. See \cite{Nadkarni}.
\end{remark}
\begin{definition}\label{def: tempered}
Let $(X,\cB_X,\set{a_t})$ be as in the Borel setup and let $\mu$ be a
finite $\set{a_t}$-invariant measure. Let 
$\cS$ be a $\mu$-cross-section. 
Let $M \in \N$ and let $E \subset \cS$ be Borel. We say
\index{tempered subset of a cross-section} that {\em $E$
  is $M$-tempered} if for every $x \in \mathcal{S}$,
$$
\# \{ t \in [0,1]: a_tx \in E\} \leq M.
$$
We say that $E$ is {\em tempered} if it is $M$-tempered for some $M$,
and that $\cS$ is a {\em tempered $\mu$-cross-section} if this condition holds for
$E = \cS$. 
\end{definition}

Clearly $\cS$ is tempered if $\cS = \cS_{\geq \vre} $ for some
$\vre>0$. It is not hard to show (see \cite[Proof of
Prop. 19]{Cheung_Chevallier}) that $\mu_{\cS}(E)<\infty$ for any tempered
subset.

 \subsection{Tight convergence of measures on lcsc spaces}
 Let $X, \cB_X$ be as in the lcsc setup, let $C_b(X)$
 denote the collection of 
bounded continuous real-valued functions on $X$, and let $\cM(X)$ denote
the collection of finite regular
Borel measures on $X$. Whenever we discuss convergence of measures in
this paper, it  
will be assumed that the measures belong to $\cM(X)$. In
particular, although infinite measures may appear in the discussion,
convergence of 
measures will only be discussed for finite measures. For $\mu \in
\cM(X)$ and $f \in C_b(X)$, we denote the
integral $\int_X f \, d\mu$ by $\mu(f)$. 
We will use the so-called \textit{tight topology} \index{tight
  topology on measures} (sometimes referred
to as {\em strict topology}) on $\cM(X)$ for
which convergence $\mu_k\to\mu$ is defined 
by either of the following equivalent requirements (see
\cite[\S 5, Prop. 9]{Bourbaki_translation} for the equivalence):
\begin{enumerate}[(i)]
\item For all  $f \in C_b(X)$, $\mu_k(f)\to
  \mu(f)$. 
\item\label{1315top} For any compactly supported continuous function
  $f:X\to\bR$, $\mu_k(f)\to \mu(f)$ and $\mu_k(X)\to \mu(X)$. 
\end{enumerate}
Note that this is not equivalent to weak-* convergence\index{weak-*
  topology on measures}, in which $f$
is taken to be compactly supported. Due to the
characterization~\eqref{1315top}, the   
weak-* topology is coarser than the tight topology when $X$ is not
compact. Nevertheless, when the total masses $\mu(X), \mu_k(X)$ are the
same (e.g., when  
they are probability measures), these notions of convergence coincide.

For a topological space $X$ and $E \subset X$,  $\on{int}(E), \,
\on{cl}(E)$ and $\partial (E)  = \on{cl}(E) \cap \on{cl}(X \sm E)$ denote
respectively its topological  interior, closure, and
boundary.
\index{I@$\on{int}_X(E) $ -- relative interior of $E$}
\index{D@$\partial_X(E) $ -- relative boundary of $E$ in $X$}
\index{C@$\on{cl}_X(E) $ -- relative closure of $E$}
Since we will work
with several topological spaces, if we want 
to stress the dependence on $X$ we will write $\on{int}_X(E), \,
\on{cl}_X(E), \, \partial_X( E)$.

\begin{definition}
Let $X$ be an lcsc space and $\mu\in \cM(X)$. We say \index{Jordan
  measurable with respect to a measure} that $E\in \cB_X$
is {\em Jordan measurable with respect to $\mu$} (abbreviated {\em
  $\mu$-JM}) if \index{J@ $\mu$-JM -- Jordan measurable with respect
  to $\mu$} $\mu(\partial_X (E)) = 0$.   
\end{definition}

The collection of $\mu$-JM sets forms a sub-algebra of $\cB_X$, and
this  algebra is rich enough to capture the tight convergence to $\mu$. More 
precisely we have the following.
\begin{lemma}[See \cite{Billingsley}, Thms. 2.1 \& 2.7, or
  \cite{Bourbaki_integration_translation},
  Chap. 4] \label{lem:convergence-equivalence}  
If $\mu_k,\mu\in \cM(X)$ then $\mu_k\to \mu$ tightly if and only if
for any $\mu$-JM set $E$ one has  
$\mu_k(E)\to\mu(E)$. 

Moreover, if $Y$ is also an lcsc space and $\psi: X \to Y$ is a measurable function, then the push-forward map
$\psi_*:\cM(X)\to \cM(Y)$ is continuous at a measure $\mu\in \cM(X)$ (with respect to the tight topologies) provided that $\psi$ is continuous $\mu$-almost everywhere. 
\end{lemma}

\section{Equidistribution of visits to a cross-section}\label{sec: JM sections}
Throughout this section we let $(X, \cB_X, \{a_t\}, \mu)$ be as in the lcsc setup and
$\cS$ be a $\mu$-cross-section as in Definition~\ref{def:cs}.
 We further assume  that
\begin{enumerate}[(A)]
\item\label{item: probability} $\mu$ is a probability measure.
\item\label{item: B}
  $\mu_\cS$ is finite.
\item\label{item: C}
  $\cS$ is lcsc (with respect to its subset topology induced by the
  topology on $X$). 
  \end{enumerate}
Let $\chi_E$ be the
indicator of $E \subset X$. 
 For $x\in X$ and $E\subset \cS$ we let \index{N@$N(x,T,E)$ -- number
   of visits of orbit of $x$ to $E$ up to time $T$}
\begin{align*}
N(x,T,E) &=\#\set{t\in [0,T] : a_tx\in E}.
\end{align*}

\begin{definition}\label{def: genericity}\;
\begin{enumerate}
\item We say that a point \index{G@$(a_t,\mu)$-generic}
  \index{G@$(a_t,\mu_{\cS})$-generic}
$x\in X$ is {\em
    $(a_t,\mu)$-generic} if 
  $\frac{1}{T}\int_0^T\del_{a_tx}dt\to_{T \to \infty} \mu$.   

\item We say  
that $x\in X$ is {\em $(a_t, \mu_\cS)$-generic} 
if the sequence of visits of the orbit $\{a_tx : t>0\}$ to $\cS$ 
equidistributes with respect to
$\frac{1}{\mu_\cS(\cS)}\mu_\cS$. 
\item For a Borel subset $\cS'\subset \cS$ which is $\mu_\cS$-JM and of positive $
\mu_\cS$-measure, we say that $x\in X$ is
{\em $(a_t, \mu_\cS|_{\cS'})$-generic} if the sequence of visits of
the orbit $\set{a_tx:t>0}$ to $\cS'$ equidistributes with respect to 
$\frac{1}{\mu_\cS(\cS')}\mu_\cS|_{\cS'}$.
\end{enumerate}
\end{definition}
Since, by \eqref{item: probability}, both cases above concern
probability measures, we can understand the
equidistribution equivalently as either weak-* convergence or as tight convergence. 
Also, by 
 Lemma~\ref{lem:convergence-equivalence}, $x$ is $(a_t,
\mu)$-generic if and only if 
\begin{equation}\label{eq:0034}
\text{for any } \mu\textrm{-JM set } 
E\subset X, \ \ \ \frac{1}{T}\int_0^T \chi_E(a_tx) dt \to_{T \to \infty} \mu(E).
\end{equation}
Similarly, by Lemma~\ref{lem:convergence-equivalence}, and using
\eqref{item: B} and \eqref{item: C}, $x $ is
$(a_t, \mu_\cS)$-generic if and only if  
\begin{equation}\label{eq:0022}
\text{for any } \mu_\cS\textrm{-JM set }E\subset \cS, \ \ \ \ 
\frac{N(x,T,E)}{N(x,T,\cS)}\to_{T \to \infty} 
\frac{\mu_\cS(E)}{\mu_\cS(\cS)}.
\end{equation}

\begin{remark}
Note that we do not define genericity with respect to the first return
map $T_\cS: \cS \to \cS$, but the reader will note the relationship between $(a_t,
\mu_\cS)$-genericity and $(T_\cS, \mu_\cS)$-genericity: $x \in X$ is $(a_t,
\mu_\cS)$-generic if and only if $x' = a_{\tau_\cS(x)}x$ is
$\left(T_\cS, \frac{1}{\mu_{\cS}(\cS)} \,
\mu_\cS\right )$-generic (in the natural sense).
\end{remark}

In this section we will study the relationship between
$(a_t,\mu)$-genericity and $(a_t,\mu_\cS)$-genericity.  
As a motivating example, consider the following simplified
situation. Assume that (i) for any $\mu_\cS$-JM 
set $E\subset \cS$, the thickened set $E^{(0,1)}$ (defined via
\eqref{eq: thickening}) is $\mu$-JM; and  (ii) the return time
function $\tau_\cS$ is bounded below by 1. Then it
follows that if $x\in X$ is $(a_t,\mu)$ generic, then it is
$(a_t,\mu_\cS)$-generic. Indeed, we have 
$$\frac{1}{T}N(x,T,E) = \frac{1}{T} \left( \int_0^T
  \chi_{E^{(0,1)}}(a_tx)dt + O(1) \right)$$
and the claim can be deduced using~ \eqref{eq:0034}, \eqref{eq:0022},
and Theorem~\ref{cor:csmeasure}\eqref{1204-2}. 

Since  the notions of cross-sections in Definitions~\ref{def: Borel cs}, \ref{def:cs} refer only to the Borel
structure,
while equidistribution is a topological notion, it is
natural to expect a topological assumption like (i), and indeed, we
will require some similar additional assumptions. Regarding (ii), when
trying to remove it, one 
encounters a complication involving the relation between the
continuous time parameter $T$  
and the number of visits $N(x,T,\cS)$. Suppose $x$ is
$(a_t,\mu)$-generic but the number of visits up to time $T$ 
is large compared to $T$ (e.g., on the order of
$T^2$). This implies frequent visits to $\cS_{<\vre}$. In 
Theorem~\ref{thm:Sgenericity} below, which is the main result of this
section, we use this observation to define a concrete $\mu$-null set, and show
that for trajectories
avoiding this set, $(a_t,\mu)$-genericity implies  
$(a_t,\mu_\cS)$-genericity.

\begin{definition}\label{def:mureasonable}
 Let $\cS\subset X$ be a $\mu$-cross-section satisfying \eqref{item:
   probability}, \eqref{item: B} and
 \eqref{item: C}. We say that $\cS$ is a \index{reasonable cross-section}
 \index{M@$\mu$-reasonable 
   cross-section} {\em $\mu$-reasonable} 
 if in addition, the following hold: 
  \begin{enumerate}
\item\label{item: reasonable 1} For all sufficiently small $\vre$, the sets $\cS_{\ge\vre}$ are
  $\mu_\cS$-JM.  
\item\label{item: reasonable 2}
There is a relatively open subset $\cU\subset\cS$ such that 
the following two conditions hold:
\begin{enumerate}
\item The
map $(0,1)\times \cU\to X$, $(t,x)\mapsto a_tx$ is open;
\item 
$\mu\pa{(\on{cl}_X(\cS)\smallsetminus \cU)^{(0,1)}} =0$.
\end{enumerate}
  \end{enumerate}
\end{definition}
\begin{remark}
We note that the interval $(0,1)$ in Definition~\ref{def:mureasonable} can be replaced by any fixed small interval. We also note
that it is possible to obtain our results while replacing Definition
\ref{def:mureasonable}\eqref{item: reasonable 1} 
with  the  weaker requirement that there exists an increasing
collection 
of $\mu_\cS$-JM sets 
$\cF_k$ such that $\cS=\bigcup_k\cF_k$ 
and such that $\cF_k\subset \cS_{\ge \vre_k}$ for some
$\vre_k>0$. This more flexible framework requires slightly more
involved arguments, but we will not need it and leave the details
to the dedicated reader. 
\end{remark}

The following elementary  lemma shows that Definition
\ref{def:mureasonable} implies 
the property (i) used in the preceding discussion.
\begin{lemma}\label{lem:JM}
If $\cS$ is $\mu$-reasonable then for any $\mu_\cS$-JM set $E\subset
\cS$ and any interval $I\subset [0,1]$,  $E^I$ is $\mu$-JM. 
\end{lemma}

\begin{proof}
Assume for concreteness that $I$ is a closed interval, the other cases
being similar, and write $I = [\tau_1, \tau_2]$. Let $\cU\subset \cS$ be the
relatively open set appearing in  
Definition~\ref{def:mureasonable}\eqref{item: reasonable 2}.
We need to show that $\mu \left( \partial_X \left (E^I\right)
  \right)=0$. This will follow once we show 
\begin{equation}\label{eq:042}
\partial_X \left(E^I \right)\subset (\on{cl}_X(\cS)\smallsetminus \cU)^{I}\cup
a_{\tau_1}\cS\cup a_{\tau_2} \cS\cup (\partial_\cS (E))^I
\end{equation}
Indeed, all sets on the RHS of~\eqref{eq:042} are $\mu$-null: the
first by the choice of $\cU$, the fourth because of
Theorem~\ref{cor:csmeasure}\eqref{1204-4}
and the assumption that   
$\partial_\cS (E)$ is $\mu_\cS$-null, and the second and third
sets are $\mu$-null by standing assumption \eqref{item: B} and Theorem
\ref{cor:csmeasure}\eqref{1204-1}. 

We prove~\eqref{eq:042}. Let $x\in \partial_X\left (E^I \right)$, so
that there is a sequence $t_k\in[\tau_1,\tau_2]$ and
$y_k\in E$ such that $a_{t_k}y_k\to x$. By passing to subsequences we may 
assume that $t_k\to t_0 \in I$ and $y_k\to
y_0 \in \on{cl}_X(E)\subset \on{cl}_X(\cS)$. We distinguish several cases. 
If $y_0\notin \cU$, then $x$ clearly belongs to the RHS of
\eqref{eq:042}. 
Thus we assume that $y_0\in \cU$. If $t_0 \in\set{\tau_1,\tau_2}$ then again 
$x$ clearly belongs to the RHS of
\eqref{eq:042}.
Assume
that $\tau_1<t_0<\tau_2$ and $y_0\in\cU$, hence in particular $y_0\in
\on{cl}_{\cS}(E)$. If $y_0\notin \on{int}_\cS(E)$ then by definition
$y_0\in \partial_{\cS} (E)$ 
and again $x$ belongs to the RHS of~\eqref{eq:042}. The only remaining
possibility is that $y_0\in \on{int}_{\cS}(E)$ but this is impossible
since $x \notin \on{int}_X\left(E^I \right)$ but the map $(t,y)\mapsto
a_ty$ is open from $(0,1)\times\cU$ to $X$.
\end{proof}

For a set $E\subset \cS_{\ge\vre}$ the relation between 
$N(x,T,E)$ and $\int_0^T\chi_{E^{(0,\vre)}}(a_tx)dt$ is simple:
\begin{proposition}\label{prop:JM}
Let $\cS$ be a $\mu$-reasonable cross-section. Then
 for any $x\in X$ which is $(a_t,\mu)$-generic, for any $\vre>0$ and
 for any $\mu_\cS$-JM set $E\subset \cS_{\ge\vre}$,  
 \begin{align}\label{eq:1230}
\lim_{T\to\infty}\frac{1}{T}N(x,T, E) &=\mu_\cS(E).
\end{align}
\end{proposition}
\begin{proof}
By Lemma~\ref{lem:JM} we have that
$E^{(0,\vre)}$ is $\mu$-JM, and thus by
Lemma~\ref{lem:convergence-equivalence},
$ \frac{1}{T}\int_0^T \chi_{E^{(0,\vre)}}
(a_tx)dt \to_{T \to \infty} \mu\left(E^{(0,\vre)} \right).$ Since
$E\subset \cS_{\ge \vre}$, we obtain (using Theorem \ref{cor:csmeasure}) that
$$ \mu(E^{(0,\vre)}) = \vre\mu_\cS(E) \ \text{ and } \ \int_0^T
\chi_{E^{(0,\vre)}}  
(a_tx)dt = \vre N(x,T,E) + O(1),$$ and \eqref{eq:1230} follows. 
\end{proof}
The following sets will be useful for analyzing trajectories $\{a_tx: t>0\}$ visiting 
$\cS_{<\vre}$ with abnormally large frequency. 
 Let \index{D@$\Del_{\cS,\del} $ -- points with $\delta$-frequent
   abnormal visits to all $\cS_{<\vre}$} \index{D@$\Del_{\cS} $ -- points with
   abnormal visits to all $\cS_{<\vre}$}
\begin{align}\label{def:badset}
\Del_{\cS,\del} & \df  \set{x\in \cS :  \forall \vre >0, \ \limsup_{T \to
                  \infty} \frac{1}{T} N(x,T, \cS_{<\vre})>\del},\\ 
\nonumber\Del_{\cS}& \df  \bigcup_{\del>0} \Del_{\cS,\del}.
\end{align}

We have the following variant of Proposition~\ref{prop:JM}, in which
we do not require that $E\subset \cS_{\ge\vre}$. 
\begin{proposition}\label{prop:JMonallS}
 Let $\cS$ be a $\mu$-reasonable cross-section.  Assume that  $x\in X\smallsetminus
\Del_\cS^\bR$ is $(a_t,\mu)$-generic (where 
$\Delta_\cS^{\bR}$ is the thickening of $\Delta_\cS$ as in \eqref{eq:
  thickening}). Then \eqref{eq:1230} holds
for any $\mu_\cS$-JM set $E\subset \cS$. 
\end{proposition}
\begin{proof}
 By Definition~\ref{def:mureasonable}, for any small enough $\vre>0$ we have that $E\cap\cS_{\ge\vre}$ is $\mu_\cS$-JM, as
 an intersection of two such sets  
and so by Proposition~\ref{prop:JM},
$$\liminf_{T \to \infty} \frac{1}{T}N(x,T,E)\ge \lim_{T \to \infty} \frac{1}{T}N(x,T,E\cap
\cS_{\ge\vre}) = \mu_\cS(E\cap \cS_{\ge\vre}).$$  
Since the sets $\cS_{\ge\vre}$ exhaust $\cS$, and $\vre$ can be chosen
arbitrarily small, we find 
\begin{equation}\label{eq:214}
\liminf_{T\to \infty}\frac{1}{T}N(x,T,E)\ge \mu_\cS(E).
\end{equation}

Fix $\del>0$. Since $x \notin \Delta_{\cS}^\bR$, we have $x\notin
\Del_{\cS,\del}^\bR.$ By \eqref{def:badset} 
there exists $\vre>0$ so that
\begin{equation}\label{eq:206}
\limsup_{T \to \infty} \frac{1}{T}N(x,T,\cS_{<\vre})\le \del,
\end{equation}
and clearly we may take $\vre$ arbitrarily small to ensure that $\cS_{\ge \vre}$ is $\mu_\cS$-JM.
Taking~\eqref{eq:206} into account and applying again
Proposition~\ref{prop:JM} we get  
\begin{align}\label{eq:2222}\nonumber
\limsup_{T \to \infty} \frac{1}{T}N(x,T,E)  &=
                                         \limsup_{T\to \infty} \frac{1}{T}\pa{N(x,T,E\cap
                                         \cS_{\ge \vre})+N(x,T,E\cap
                                         \cS_{<\vre})}\\ 
&\le \lim_{T\to\infty} \frac{1}{T}N(x,T,E\cap \cS_{\ge \vre})
           +\limsup_{T \to \infty}
           \frac{1}{T}N(x,T,\cS_{<\vre}) \\ 
\nonumber & \le \mu_\cS(E\cap\cS_{\ge \vre}) +\del \le \mu_\cS(E)
            +\del. 
\end{align}
Since $\del$ was arbitrary, we get \eqref{eq:1230}.
\end{proof}

The following lemma shows that the extra assumption in
Proposition~\ref{prop:JMonallS}, namely that $x\notin\Del_\cS^\bR$, 
is almost harmless. 
 \begin{lemma}\label{lem:bad-is-null}
Let $\cS$ be a $\mu$-reasonable cross-section. Then
$$\mu_\cS \left(\Del_\cS \right) =0 \  \text{ and } \ \mu
\left(\Del_\cS^{\bR} \right) =0.$$
\end{lemma}
\begin{proof}
By
Theorem~\ref{cor:csmeasure} items \eqref{1204-4} and \eqref{item: ergodic decomposition} it is enough to assume 
$\mu$ is $a_t$-ergodic and to show that
$\mu_\cS\left(\Del_\cS \right) =0$. 
 We fix $\del>0$ and show that  
$\mu_\cS\left(\Del_{\cS,\del}\right) = 0$.

Take $0<\vre_1<\vre_0$ small enough so that  
$\mu_\cS(\cS_{\ge \vre_0})>0$ and  
$\mu_\cS(\cS_{<\vre_1})<\del$. This is possible by \eqref{item: B} and
because $\cS = \bigcup_{\vre>0}  \cS_{\geq \vre}.$
We will show
that 
for $\mu$-a.e. $x$, the cross-section measure
$\mu_\cS$ satisfies 
\begin{equation}\label{eq:1715}
\lim_{T \to \infty} \frac{1}{T}N(x,T,\cS_{<\vre_1}) = \mu_\cS(\cS_{< \vre_1}).
\end{equation}
This will imply $\mu_\cS\left(\Del_{\cS,\del} \right)=0$ as desired.

Note that for $E\subset \cS$, the ratio $\frac{N(x,T,E)}{N(x,T,\cS)}$
is an ergodic average for $\chi_E$ in the dynamical system
$(\cS, \frac{1}{\mu_\cS(\cS)}\mu_\cS, T_\cS),$ which is ergodic by Theorem~\ref{cor:csmeasure}\eqref{item: ergodic decomposition}.  
Applying the pointwise ergodic theorem to the characteristic
functions of the sets $\cS_{\ge \vre_0}, \cS_{<\vre_1}$ we deduce that
there is a set $F\subset \cS$ of full $\mu_\cS$-measure such that for 
any $y\in F$,
\begin{align*} &\frac{N(y,T,\cS_{\ge \vre_0})}{N(y,T,\cS)} \to_{T \to
                 \infty} \frac{\mu_\cS(\cS_{\ge\vre_0})}{\mu_\cS(\cS)}\\ 
               & \frac{N(y,T,\cS_{<\vre_1})}{N(y,T,\cS)}\to_{T \to \infty}
                                                                             \frac{\mu_\cS(\cS_{<\vre_1})}{\mu_\cS(\cS)},
\end{align*} 
and thus
\begin{equation}\label{eq:2359}
\lim_{T \to \infty}\frac{\frac{1}{T} N(y,T,\cS_{<\vre_1})}{\frac{1}{T}N(y,T,\cS_{\ge \vre_0})} =
\lim_{T \to \infty}     
\frac{N(y,T,\cS_{<\vre_1})}{N(y,T,\cS_{\ge\vre_0})} =
\frac{\mu_\cS(\cS_{<\vre_1})}{\mu_\cS(\cS_{\ge \vre_0})} . 
\end{equation}
Moreover, replacing $F$ with a smaller set of full $\mu_\cS$-measure, we can
assume that each $y \in F$ is also $(a_t, \mu)$-generic. 
Applying Proposition~\ref{prop:JM} for $E=\cS_{\ge\vre_0}$ we get that
the denominator on the LHS of \eqref{eq:2359}  converges to the denominator on  
the RHS. This implies \eqref{eq:1715} for $y \in F$.
\end{proof}

For tempered subsets of $\cS$ (see Definition \ref{def: tempered}) we can prove a
version of Proposition~\ref{prop:JM} without reference  
to the problematic set $\Del_\cS$.
\begin{proposition}\label{prop:JMtempered}
Let $\cS$ be a $\mu$-reasonable cross-section and
let $E\subset \cS$ be a tempered subset which is $\mu_\cS$-JM. If
$x\in X$ 
is $(a_t,\mu)$-generic then \eqref{eq:1230} holds.
\end{proposition}
For the proof we will need the following lemma which provides the substitute
for the assumption $x\notin \Del_\cS^\bR$. 
\begin{lemma}\label{lem:JMtempered}
Let $\cS$ be a $\mu$-reasonable cross-section for $(X,a_t,\mu)$ and
let $F\subset \cS$ be a $\mu_\cS$-JM set which is  
$M$-tempered. Then for any $x$ which 
is $(a_t,\mu)$-generic we have 
\begin{equation}\label{eq:2245}
\limsup_{T \to \infty} \frac{1}{T}N(x,T,F) \le M
\mu\left(F^{(0,1)} \right).
\end{equation}
\end{lemma}
\begin{proof}
Let $I = (0,1)$. By Lemma~\ref{lem:JM}, the set $F^I$ is $\mu$-JM and
hence $
\frac{1}{T}\int_0^T\chi_{F^I}(a_tx)dt \to_{T \to
  \infty}\mu\left(F^I\right)$. Thus~\eqref{eq:2245} will 
follow from 
\begin{equation}\label{eq: from} N(x,T,F)\le M \cdot  m \left (\set{t\in[0,T] :
      a_tx\in F^I} \right) + M, \end{equation}
for each $T>0$ (where $m$ is the Lebesgue measure on $\R$). 
Let $k = N(x,T,F)$ and let $t_1<\dots <t_k$ 
be an ordering of $\set{t\in[0,T]: a_tx\in F}$. Then, the $M$-temperedness implies 
that there is a subset $J\subset \set{1,\dots, k}$ of cardinality at
least $\frac{k}{M}$,
such that for any $j_1<j_2$ in $J$ one has  
$t_{j_2}-t_{j_1} \ge 1$. For each $t_j \in J$ except perhaps the
largest, and for any $t \in (t_j, t_j+1),$ we have $t \leq T$ and
$a_tx\in F^I$. This implies that $\set{t\in[0,T] : a_tx\in
  F^I}$ contains at least $\left \lfloor \frac{k}{M} \right \rfloor -1$ disjoint
intervals of length 1, which implies \eqref{eq: from}. 
\end{proof}

\begin{proof}[Proof of Proposition~\ref{prop:JMtempered}]
The proof is very similar to that of Proposition~\ref{prop:JMonallS}. 
The inequality 
$\liminf_{T \to \infty} \frac{1}{T}N(x,T,E) \ge \mu_\cS(E)$
follows as in \eqref{eq:214}. On the other hand, similarly to
\eqref{eq:2222}, using Proposition~\ref{prop:JM},  
for any sufficiently small $\vre>0$
we have 
$$\limsup_{T \to \infty} \frac{1}{T}N(x,T,E) =\mu_\cS(E\cap
\cS_{\ge\vre}) + \limsup_{T \to \infty}
\frac{1}{T}N(x,T, E_{\vre}),$$
where  $E_\vre \df  E\cap \cS_{<\vre}. $ The sets $E_\vre$ are
$M$-tempered and so by 
Lemma~\ref{lem:JMtempered}, 
$\limsup_{T \to \infty} \frac{1}{T}N(x,T,E) \le \mu_\cS(E) + M \, \mu \left(E_\vre^{(0,1)}
  \right)$. But since $\bigcap_{\vre >0} E_\vre^{(0,1)} =
  \varnothing$, we have 
$$\mu\left(E_\vre^{(0,1)}\right)\to_{\vre \to 0} 0,$$
and hence
$\limsup_{T \to \infty} \frac{1}{T}N(x,T,E) \le \mu_\cS(E)$. Putting
these inequalities together 
gives~\eqref{eq:1230}.
\end{proof}

We summarize the results of this section in the following theorem.
\begin{theorem}\label{thm:Sgenericity}
Let $\cS$ be a $\mu$-reasonable cross-section, and let $\cS' \subset
\cS$ be $\mu_\cS$-JM such that
$\mu_{\cS}(\cS')>0$. Suppose $x \in X$ 
is $(a_t,\mu)$-generic. Assume in addition that one of the following also hold:
\begin{enumerate}[(i)]
\item \label{item: tempered i}
 $x \notin \Del_\cS^\bR$;
\item \label{item: tempered ii} $\cS'$ is tempered.
\end{enumerate}
Then $x$ is $(a_t, \mu_{\cS}|_{\cS'})$-generic.
\ignore{

If $x\in X\smallsetminus \Del_\cS^\bR$ is $(a_t,\mu)$-generic then it
is $(a_t,\mu_\cS)$-generic.  
In particular, if $\mu$ is $\{a_t\}$-ergodic, then 
$\mu$-almost any point 
is $(a_t,\mu_\cS)$-generic.

Moreover, if $\cS'\subset \cS$ is a subset which is lcsc,
$\mu_\cS$-JM, and tempered, and is also a cross-section for
$(X,a_t,\mu)$, then
any $x\in X$ which is $(a_t,\mu)$-generic is  
also $(a_t,\mu_{\cS'})$-generic.
}
\end{theorem}

\begin{proof}
Let $x$ be $(a_t,\mu)$-generic and assume \eqref{item: tempered i}. 
Then  Proposition~\ref{prop:JMonallS} and
Lemma~\ref{lem:convergence-equivalence}  imply that $x$ is
$(a_t,\mu_\cS)$-generic.
%
Showing that $x$ is $(a_t,\mu_{\cS}|_{\cS'})$-generic
is equivalent to showing that for any  
$\mu_\cS$-JM set $E\subset \cS'$,
$$\lim_{T \to \infty} \frac{N(x,T,E)}{N(x,T,\cS')} =
\frac{\mu_\cS(E)}{\mu_\cS(\cS')},$$
which readily follows using \eqref{eq:1230} in both numerator and
denominator. 

Now assume \eqref{item: tempered ii}. Since $\cS'$ is tempered, so is
$E \subset \cS'$, and by
Proposition~\ref{prop:JMtempered}, we once again have
\eqref{eq:1230} for both the numerator and the denominator.
\end{proof}

\ignore{
\begin{remark}\red{Uri: this remark should be made into a
    corollary. Barak: This is obvious from Kac formula, do we really
    want it? }
The convergence $\frac{1}{T}N(x,T,E)\to \mu_\cS(E)$ has a nice
interpretation in terms of the integral of the first return: 
We think of $E$ as a cross section on its own and thus it has a return
time function $\tau_E$ and a first return map $T_E$.  
If we define for $k\in \bN$
the time $T_k = \sum_{i=0}^{n-1} \tau_E(T_E^i(x))$ then $N(x,T_k,E) = k$ and so 
$$\frac{1}{T_k}N(x,T_k,E) = \pa{\frac{1}{n}\sum_{i=1}^{n-1} \tau_E(T_E^i(x))}^{-1}$$
and so we obtain that the ergodic averages of the return time to $E$ 
converge to $\mu_\cS(E)^{-1}$.  
\end{remark}
}

\ignore{
  If we are in a setting where there exists a weak-stable group $H$ acting on $X$, since 
the set of $a_t$-generic points is $H$-invariant then we have the following:
\usnote{replace relative generic with $\cS$-generic}
\begin{corollary}
If $x$ is $a_t$-generic then $Hx\smallsetminus \Del_\cS^{\bR}$
consists of $a_t$-generic points relative to $\cS$. 
\end{corollary}

}

\subsection{Continuity  of the cross-section measure construction} 
In this subsection we prove the following continuity property of the
map $\mu \mapsto \mu_{\cS}$ that will be used at the very end of the paper in Theorem \ref{thm:for SZ}: 
\begin{proposition}\label{prop:continuity}
Let $X, \cB, \{a_t\}$ be as in the  lcsc setup, and let $\mu_k,
\mu \in \mathcal{P}(X)$. Suppose $\mu$ and each $\mu_k$ are
$\{a_t\}$-invariant, and $\mu_k\to_{k
  \to \infty} \mu$ in the strict topology (or equivalently, since they
are in $\cP(X)$, in the weak-* topology). Also assume that $\cS',
\cS$ satisfy the following, for $\nu = \mu$ and for $\nu = \mu_k$, for
any $k \in \N$: 
\begin{enumerate}
\item $\cS$ is a $\nu$-reasonable cross-section.
\item $\cS'\subset \cS$ is $\nu_\cS$-JM.
\item $\cS'$ is tempered. 
\end{enumerate} 
Then $(\mu_k)_{\cS}|_{\cS'}\to_{k \to \infty} \mu_{\cS}|_{\cS'}$,
with respect to the strict topology.
\end{proposition}
\begin{remark}
It seems to us that in the above proposition we cannot deduce that $\mu_{k,\cS}\to \mu_\cS$ or even that 
there is convergence after renormalizing the measures to be probability measures. The reason for this is that 
potentially, there exists $\del>0$ such that for any $\vre>0$, $\mu_{k,\cS}(\cS_{<\vre})>\del$ for infinitely many $k$'s. This is reminiscent of escape of mass. To overcome this we restrict attention to the tempered set $\cS'$. Nevertheless, the proof does give that for any $\mu_\cS$-JM set $E\subset \cS$ we have
$$
\liminf_{k\to \infty} \left( \mu_k\right)_{\cS} (E) \geq \mu_\cS(E).
$$

\end{remark}
\begin{proof}
By Lemma~\ref{lem:convergence-equivalence} we need to show 
that for any $\mu_\cS$-JM set $E\subset \cS'$ we have 
\begin{equation}\label{eq:1310}
(\mu_k)_\cS(E)\to_{k\to\infty} \mu_\cS(E).
\end{equation} 
Take $E\subset \cS'$ a $\mu_\cS$-JM set. 
For any $\vre>0$ we can decompose $E = (E\cap \cS_{\ge\vre}) \cup
(E\cap \cS_{<\vre})$ and hence for  
$\nu\in \set{\mu,\mu_k}$ we have 
$$\nu_\cS(E) = \nu_\cS(E\cap \cS_{\ge\vre}) + \nu_\cS(E\cap \cS_{<\vre}).$$
By Lemma~\ref{lem:JM} we have that $(E\cap \cS_{\ge\vre})^{(0,\vre)}$
is $\mu$-JM, and so by 
Theorem~\ref{cor:csmeasure}\eqref{1204-2},
\[\begin{split}
    & \lim_{k\to\infty} (\mu_k)_\cS(E\cap \cS_{\ge\vre}) =\lim_{k\to
      \infty} \frac{1}{\vre} \,\mu_k\left((E\cap 
\cS_{\ge\vre})^{(0,\vre)}\right) \\ =  &
\frac{1}{\vre}\, \mu\left((E\cap \cS_{\ge\vre})^{(0,\vre)} \right) = \mu_\cS\left(E\cap
\cS_{\ge\vre}\right).
\end{split}\]
Moreover $\mu_\cS(E\cap \cS_{<\vre})\to_{\vre\to 0} 0$, because
$\mu_\cS$ is a finite measure, and hence
$$
\liminf_{k\to \infty} \left( \mu_k\right)_{\cS} (E) \geq \mu_\cS(E),
$$
and \eqref{eq:1310} will follow once
we establish that
\begin{equation}\label{eq: establish}
  \limsup_{k\to\infty} \, (\mu_k)_\cS(E\cap
  \cS_{<\vre}) \to_{\vre \to 0}0.
\end{equation}
Clearly, it is
enough to show \eqref{eq: establish} for $E=\cS'$. Assume first that each
$\mu_k$ is ergodic. Since $\cS'$ is $\mu_k$-JM, there
is $x_k\in X$ which is  
$(a_t,\mu_k)$-generic. Choose $M$ so that $\cS'$ is
$M$-tempered, then by \eqref{eq:2245}
we have 
\small{
\begin{equation}\label{eq: for ergodic}
  (\mu_k)_\cS\left(\cS'\cap \cS_{<\vre} \right) = \lim_{T \to \infty}
\frac{1}{T}N \left(x_k,T,
  \cS'\cap \cS_{<\vre} \right) \leq  M 
\mu_k\left((\cS'\cap \cS_{<\vre})^{(0,1)}\right).
\end{equation}
}
Using the ergodic decomposition and Theorem \ref{cor:csmeasure}\eqref{item: ergodic
  decomposition}, we see that \eqref{eq: for ergodic} also holds without assuming
that $\mu_k$ is ergodic. 
This gives
$$
\limsup_{k \to \infty} (\mu_k)_\cS\left(\cS'\cap \cS_{<\vre} \right)
\leq M  \mu\left( \left(\cS'\cap \cS_{<\vre} \right)^{(0,1)}\right),
$$
and taking the limit as $\vre \to 0$ we
obtain \eqref{eq: establish}. 
\end{proof}

\section{Lifting reasonable cross-sections}\label{sec:lifting}
The goal of this section is to prove Proposition \ref{prop:941} which
roughly says that a lift of a reasonable cross-section to a fiber bundle extension
is again reasonable.
We will begin with an elementary lemma about fiber bundles that will be
needed at a certain point. 

\subsection{Fiber bundles}
Let $X, Y, F$ be topological spaces. 
Recall that a continuous map $\pi:Y\to X$ between two topological
spaces is called an {\em $F$-fiber bundle} if $X$ can be covered by a
collection $\set{U_i}$ of 
open sets, called {\em trivial open sets}, 
such that  
for each index $i$, there exists a homeomorphism $\psi_i:U_i\times F
\to \pi^{-1}(U_i)$ satisfying  
$\pi(\psi_i(x,f)) = x$ for all $(x,f)\in U_i\times F$.
A
%
\textit{morphism} between two $F$-fiber bundles $Y_i \to X_i,\ i=1,2$,
is a pair of continuous maps $\phi: Y_1 \to Y_2, \ \bar \phi: Z_1 \to X_2$ for which the following
diagram commutes:
\begin{equation}\label{eq: bundle morphism}
  \xymatrix{
Y_1\ar[d]_{\pi_1}\ar[r]^\phi & Y_2\ar[d]^{\pi_2}\\ X_1\ar[r]_{\bar{\phi}}& X_2
}
\end{equation}

\begin{definition}\label{def: fiber proper}
  The bundle morphism \eqref{eq: bundle morphism} is said to be
  \textit{proper relative to 
    fibers}, if any sequence $\set{y_n}\subset Y_1$ 
for which $\set{\pi_1(y_n)}$ and $\set{\phi(y_n)}$ are bounded
sequences in $X_1,Y_2$ respectively,  is bounded in $Y_1$ (here
bounded means has compact closure). 
\end{definition}

\begin{lemma}\label{lem: about fiber bundles}
  Let $F, X_i, Y_i \ (i=1,2)$ be lcsc spaces such that $\pi_i: Y_i
  \to X_i$ are $F$-bundles. Suppose the bundle map \eqref{eq: bundle
    morphism} is proper relative to fibers, 
  $\bar \phi: X_1 \to X_2$ is open, and  
  the restriction
 $\phi|_{\pi_1^{-1}(x)}:\pi_1^{-1}(x)\to \pi_2^{-1}(\bar{\phi}(x))$
is a homeomorphism between the fibers, for any $x \in X_1$. 
Then
$\phi$ is open. 
\end{lemma}
\begin{proof}
  By
  restricting
  to trivial open sets in $X_1,X_2$ and their preimages in $Y_1,Y_2$,
  we may assume that $Y_1, Y_2$ are trivial bundles. Note that the
  properness relative to fibers holds after such restriction.  

Assume then that 
$Y_i=X_i\times F$ and $\pi_i$ is the projection on the first
coordinate. For any $x\in X_1$ write
$\phi(x,f) = (\bar{\phi}(x), \phi_x(f))$, so that 
$\phi_{x}:F\to F$ is the homeomorphism $\phi$ induces between the
fibers $\set{x}\times F$ and  
$\set{\bar{\phi}(x)}\times F$. 
Since we assumed all spaces involved are lcsc, the topology is induced
by a metric.
 Let 
$(x_0,f_0)\in X_1\times F$ and let $U\times V$ be a basic open
neighborhood of $(x_0, f_0)$ in $X_1\times F$. We will show that
$\phi(x_0,f_0)
$ is in the interior
of  
\begin{equation}\label{eq:image is open}
\phi(U\times V) = \set{(\bar{\phi}(z), \phi_z(f)):z\in U, f\in V}.
\end{equation}
For any $z\in U$, since $\phi_z$ is a homeomorphism of $F$, there
exists $\vre>0$ so that  
 $B_\vre^F(\phi_z(f_0))\subset \phi_z(V)$ (where $B_{\vre}^F(x)$ is
 the ball of radius $\vre$ in $F$ centered at $x$, with respect to a
 metric defining the topology). The following claim says
 that $\vre$ can be taken to be uniform 
 for $z$ close enough to $x_0$. 

 \medskip
 
\noindent \textbf{Claim}: {\em There exists $\vre>0$ such that for any
  $z$ close enough to $x_0$, $$B_\vre^F(\phi_z(f_0))\subset \phi_z(V).$$
  }

We first assume the Claim and complete the proof of the
Lemma. By continuity of $\phi$ at $(x_0,f_0)$, there is a
neighborhood $U'$ of $x_0$ so that for  
$z \in U'$ we have $d(\phi_{x_0}(f_0),
\phi_z(f_0))<\vre/2$. Thus for $z \in U'$, $B_{\vre/2}^F(\phi_{x_0}(f_0))\subset
B_\vre^F(\phi_z(f_0))$. By the Claim,  making $U'$ smaller if
necessary, we have that $B_{\vre/2}^F(\phi_{x_0}(f_0))\subset
\phi_z(V)$ for $z \in U'$, 
so  \eqref{eq:image is open} 
contains 
$\bar{\phi}(U')\times B_{\vre/2}^F(\phi_{x_0}(f_0))$. 
Since  $\bar{\phi}$ is open, we have shown that 
$\phi(x_0,f_0)$ is in the interior of~\eqref{eq:image is open}.

We now prove the claim. Assume to the contrary  that there exists
a sequence 
$z_n\to x_0$ for which the homeomorphisms  
$\phi_{z_n}:F\to F$ map a point $f_n\notin V$ into
$B_{1/n}^F(\phi_{z_n}(f_0))$. By the properness relative to fibers of
$\phi$, we may assume $f_n\to f_*\ne f_0$. Moreover, by the continuity
of $\phi$ we get that  
$\lim \phi(z_n,f_n) =\phi(x_0,f_*)$. On the other hand, this limit is also equal to 
$$\lim(\bar{\phi}(z_n), \phi_{z_n}(f_n)) =(\bar{\phi}(x_0),\phi_{x_0}(f_0))= \phi(x_0,f_0).$$ 
This contradicts the
assumption that $\phi_{x_0}$ is one to one. 
\end{proof}

\subsection{Extensions}
Let $\{a_t\} \curvearrowright (X, \cB_X, \mu)$ and $\{a_t\}
\curvearrowright (\wt{X}, \cB_{\wt{X}}, \tilde{\mu} )$ be two actions as in the
lcsc setup, and assume that $\mu, \tilde{\mu}$ are both probability measures. 
We say that $\pi: \wt{ X} \to
X$ is a {\em continuous factor map} or equivalently a {\em continuous extension map} if it is continuous and satisfies
$\mu = \pi_* \tilde{\mu}$ and
$a_t \circ \pi = \pi \circ a_t$ for any $t \in \R$. 

\begin{proposition}\label{prop:941}
In the above setup, let $\pi : \wt{X}
\to X$ be a continuous factor map.
Assume furthermore that $\pi$ is a fiber bundle. Let $\cS\subset X$ be a
$\mu$-reasonable cross-section and \index{S@$\widetilde{\cS}$ -- lifted
  cross-section} let $\widetilde{\cS} \df \pi^{-1}(\cS)$.  
Then:
\begin{enumerate}[(a)]
\item \label{item: 1 lift} $\wt{\cS}$ is a $\tilde{\mu}$-cross-section. 
\item\label{item: 2 lift} $\pi_*\tilde{\mu}_{\wt{\cS}} = \mu_\cS$.
  \item\label{item: 2.5 lift} $\tilde{\mu}_{\wt{\cS}}$ is finite.
\item \label{item: reasonable lift}
  $\wt{\cS}$ is $\tilde{\mu}$-reasonable.
\item \label{item: 3 lift}
  For any $\mu_\cS$-JM set $E\subset \cS$, $\pi^{-1}(E)\subset
  \wt{\cS}$ is $\mu_{\wt{\cS}}$-JM.  
\item\label{item: 4 lift} If $E\subset \cS$ is $M$-tempered then
  $\pi^{-1}(E)$ is $M$-tempered. 
\end{enumerate}
\end{proposition}
\begin{remark}
In the proof of Proposition~\ref{prop:941}, the only place that we use the assumption that 
$\pi$ is a fiber bundle and not simply a continuous extension is in the proof of item \eqref{item: reasonable lift}.
\end{remark}
\begin{proof}
  Since $a_t \circ \pi = \pi \circ a_t$, we have
  \begin{equation}\label{eq: thickening lift}
 \left(\pi^{-1}(E)\right)^{I} = \pi^{-1}\left( E^I \right)
 \end{equation}
 for the thickened sets as in \eqref{eq: thickening}. In addition, for
 $ x = \pi(\tilde{x})$ we
 have   
\begin{equation}\label{eq:returntimes}
\set{t: a_t\tilde{x}\in\wt{\cS}} = \set{t:a_tx\in \cS}, 
\end{equation} 
and hence
\begin{equation}\label{eq:1029}
\cY_{\tilde{x}} = \cY_x, \  \ \ \tau_{\cS}(x) =
\tau_{\wt{\cS}}(\tilde{x}), \  \ \ \text{ and } \ \
 \wt{\cS}_{\ge \vre}= \pi^{-1}(\cS_{\ge \vre}). 
 \end{equation}
We prove  \eqref{item: 1 lift}: since $\cS$ is a $\mu$-cross-section,
there exists an $\set{a_t}$-invariant set $X_0\in \cB_X$ 
such that $\mu(X_0) = 1$ and $\cS\cap X_0$ is a
Borel cross-section (see Definitions \ref{def: Borel cs},
\ref{def:cs}). It follows that if $\wt{X}_0 \df \pi^{-1}(X_0)$, then
$\wt{X}_0\in \cB_{\wt{X}}$ is an $\set{a_t}$-invariant set 
with $\tilde{\mu}(\wt{X}_0) = 1$. 
Finally, it follows from \eqref{eq:1029} that $\wt{\cS}\cap
\wt{X}_0$ is a Borel cross-section according to Definition \ref{def:
  Borel cs}. This finishes the verification of Definition~\ref{def:cs}
and proves \eqref{item: 1 lift}.

Item~\eqref{item: 2 lift} follows from 
\eqref{eq:1019}, \eqref{eq: thickening lift},  and the assumption that
$\pi_*\tilde{\mu} = \mu$.


Item \eqref{item: 2.5 lift} follows from  \eqref{item: 2 lift} since the $\mu$-reasonability of $\cS$ implies that $\mu_{\cS}$ is finite. 

Since $\pi:\wt{\cS}\to \cS$ is continuous, for any
$E\subset \cS$,
$$\partial_{\wt{\cS}}(\pi^{-1}(E)) \subset
\pi^{-1}(\partial_\cS E).$$
It follows that if $E\subset \cS$ is $\mu_\cS$-JM then 
$$\tilde{\mu}_{\wt{\cS}}(\partial_{\wt{\cS}}(\pi^{-1}(E)) ) \le
\tilde{\mu}_{\wt{\cS}}(\pi^{-1}(\partial_\cS E)) =
\mu_\cS(\partial_\cS E) = 0,$$
and \eqref{item: 3 lift} follows. 

Item \eqref{item: 4 lift} follows from the equivariance of $\pi$. 

It remains to prove \eqref{item: reasonable lift}. We verify Definition~\ref{def:mureasonable}.
We note that $\wt{\cS}$ satisfies \eqref{item: probability}, \eqref{item: B} and
\eqref{item: C}. Indeed, \eqref{item: probability} is part of our assumptions, \eqref{item: B} follows from \eqref{item: 2
  lift} and \eqref{item: C} 
follows easily from the fact that $\wt{X}$
and $\cS$ are lcsc and $\pi$ is continuous.
We now check the further conditions of
Definition~\ref{def:mureasonable}. Condition
\ref{def:mureasonable}\eqref{item: reasonable 1} follows from \eqref{item: 3 lift} and
\eqref{eq:1029}. 
For the technical condition \ref{def:mureasonable}\eqref{item:
  reasonable 2}, let $\cU\subset \cS$ be the 
subset appearing in the definition for $\cS$, and let $\wt{\cU} =
\pi^{-1}(\cU)$. 
First,
$$\on{cl}_{\wt{X}}(\wt{\cS})\smallsetminus \wt{\cU}
\subset \pi^{-1}\left(\on{cl}_X(\cS) 
\smallsetminus \cU \right)$$
and since $\pi_*\tilde{\mu} = \mu$,
$$\tilde{\mu}\pa{\on{cl}_{\wt{X}}(\wt{\cS})\smallsetminus
  \wt{\cU}}\le  
\mu\left( \on{cl}_X(\cS)\smallsetminus \cU \right) = 0.$$
In order to verify that the map $(0,1)\times \wt{\cU}\to
\wt{X}$, $(t,y)\mapsto a_ty$ is open, 
note that $\wt{\cU}^{(0,1)} = \pi^{-1}(\cU^{(0,1)})$ is open in
$\wt{X}$, and thus 
 it is enough to show that the map $(0,1)\times \wt{\cU}\to
 \wt{\cU}^{(0,1)}$ is open. Consider the 
commutative diagram
\begin{equation}\label{eq:diagram of bundles}
\xymatrix{
(0,1)\times \wt{\cU}\ar[rrr]^{(t,\tilde{x}) \ \mapsto \  a_t\tilde{x}} \ar[d]_{(\on{id},\pi)} &&& \wt{\cU}^{(0,1)}\ar[d]^\pi\\
(0,1)\times \cU\ar[rrr]_{(t,x) \ \mapsto \  a_tx} &&& \cU^{(0,1)}
}
\end{equation}
in which the vertical maps are fiber bundles and the horizontal maps
constitute a morphism of fiber bundles. Note that the map between the
bottom spaces is open by assumption and the map between the top spaces
is a homeomorphism when restricted to single fibers.  Moreover, we
claim that this  
morphism of bundles is proper relative to fibers in the sense of
Definition~\ref{def: fiber proper}. An application of 
Lemma~\ref{lem: about fiber bundles} then shows that the upper
horizontal map is open which finishes the proof.  

We verify properness relative to fibers. Let $\set{(t_n,\wt{x}_n)}\subset (0,1)\times \wt{\cU}$ be a sequence such that 
both $\set{(t_n,\pi(\wt{x}_n))}\subset (0,1)\times \cU$ and
$\set{a_{t_n}\wt{x}_n}\subset \wt{\cU}^{(0,1)}$ are bounded. We may
assume without loss of generality that  
$$(t_n,\pi(\wt{x}_n))\to (t_0,x_0)\in (0,1)\times \cU \textrm{ and } 
a_{t_n}\wt{x}_n\to \wt{z}\in \wt{\cU}^{(0,1)}.$$ 
It follows first that $t_0\in (0,1)$, and second,
by applying $\pi$ we get that $\wt{z}\in \pi^{-1}(a_{t_0}x_0)$, which implies that 
$a_{-t_0}\wt{z}\in \pi^{-1}(\cU) = \wt{\cU}.$ Finally,
by the continuity of the action we see that
that  
$$\lim \wt{x}_n = \lim a_{-t_n}a_{t_n}\wt{x}_n  = a_{-t_0}\wt{z} \in \wt{\cU}$$
which shows that $\set{(t_n,\wt{x}_n)}$ converges in $(0,1)\times \wt{\cU}$ to $(t_0,a_{-t_0}\wt{z})$. 
\if
In fact, as $\rho = \vre/2$ this map is 1-1 and onto and so we might as well check the continuity of its inverse. Assume  
that $a(t_k)y_k\to a(t_0)y_0$ where $(t_k,y_k),(t_0,y_0)\in (0,\rho)\times \wt{\cU}_\vre$. 
Applying $\pi$ and using the fact that $(0,\rho)\times \cU_\vre\to \cU_\vre^{(0,\rho)}$ is a homeomorphism,
we deduce that $t_k\to t_0,$ and that any accumulation point of $y_k$ belongs to $\pi^{-1}(\pi(y_0))$. 
Assume for a moment that after passing to a subsequence, $y_k$ converges to some $z_0$. 
Then we get that $a(t_0)z_0 \leftarrow a(t_k)y_k\rightarrow a(t_0)y_0$ and so $z_0=y_0$. That is, the only possible accumulation
point of $y_k$ is $y_0$. Now since $y_k$ is a bounded sequence in the lcsc space $\wt{X}$ 
(as $a(t_k)y_k$ converges and $t_k$ is bounded), we conclude that since it has only one accumulation
point, it actually converges to it. That is, we established that $y_k\to y_0$ which together with $t_k\to t_0$ which was already established,
constitutes the continuity of the inverse map and finishes the proof.
\fi
\end{proof}

\section{Homogeneous spaces and homogeneous measures} \label{sec:
  spaces and measures}
\subsection{Guide to the rest of the paper}
In order to obtain our results, we will apply the theory developed in
\S\S\ref{sec: prelims}-\ref{sec:lifting} in order to get
equidistributed sequences on cross-sections in various spaces, and
then interpret them as being related to best approximations and
$\vre$-approximations. For our results we will have three distinctions
which give rise to $8 =
2^3$ cases. The first distinction is between best approximations and
$\vre$-approximations (compare Theorems \ref{thm: main best} and
\ref{thm: main epsilon}), the second distinction is between Lebesgue
a.e.\,vector and vectors with entries in a totally real field (compare both of the above, with
Theorem \ref{thm: number field}) and the third is between results on
equidistribution in the real spaces $\crly{E}_n \times \R^d$ and in the
larger locally compact space $\crly{E}_n \times \R^d \times \widehat{\Z}^n$
(compare the first two  
coordinates of \eqref{eq: seq bundle best} with all three of them). In all cases
we will specify the dynamical system, define our
cross-section, give an expression for the 
cross-section measures, and check that the axioms used in \S\S\ref{sec:
  prelims}-\ref{sec:lifting} are satisfied. Some of these verifications
will be routine but others will require detailed argumentation.

Since the reader may not be equally interested in all eight cases, we
preface this discussion with a short guide. 
We will first consider the {\em real homogeneous space} $\XX_n =
\SL_n(\R)/\SL_n(\Z)$, where the
action is given by left 
multiplication by the one-parameter group $\{a_t\}$ in \eqref{eq: def
  at}.
%
Understanding the space $\XX_n$ will only give information about the
real components $\crly{E}_n \times \R^d$ in \eqref{eq: seq bundle best}. 
On $\XX_n$ we will
define a cross-section $\sro$ (defined below in \eqref{eq:main
  sets}). We will consider two kinds of measures 
$\mu$ 
on $\XX_n$. The first is the 
Haar-Siegel measure $m_{\XX_n}$. This measure will give information
about the properties of Lebesgue a.e.\;$\theta \in \R^d$. When
discussing this measure we will 
say that we are in {\em Case I}. In the second case, related to
approximation in totally real fields, which we will refer
to as {\em Case II}, \index{Case II (measures from totally real number
  fields)} we will have a homogeneous measure
$m_{\vec{\alpha}}$ (see Proposition \ref{prop: explicit compact
  push}), where by {\em homogeneous} we mean that there is a closed
subgroup $L \subset 
\SL_n(\R)$, such that $m_{\vec{\alpha}}$ is $L$-invariant and $\supp\pa{
m_{\vec{\alpha}}}$ is a closed $L$-orbit. In our case $L$ is a
conjugate of the diagonal group $A$, where the conjugating matrix will
depend on the algebraic vector $\vec{\alpha}$, and $\supp\pa{
m_{\vec{\alpha}}}$ is compact. We will prove that $\sro$ is
$\mu$-reasonable in both Case I and Case II (see \S \ref{subsec: 8.3} and \S
\ref{subsec: 8.4} respectively). 

In order to derive information about best approximations, we will work
with a subset $\cB \subset \sro$, and for $\vre$-approximations, we
will work with a subset 
$\cS_{\vre} \subset \sro$. We will show (see \S \ref{subsec: 9.1}, \S
\ref{subsec: 9.2}) that these sets are
$\mu_{\sro}$-JM in both cases. As remarked in the introduction,
$\cB$ will be a tempered subset (see Proposition \ref{prop: B
  tempered}), but $\cS_{\vre}$ will not be.  

In order to understand all three components in \eqref{eq: seq bundle best},
we will consider the {\em 
adelic homogeneous space }
$\XXnA \df \SL_n(\bA)/\SL_n(\Q)$, and the factor map $\pi: \XXnA \to
\XX_n$. These will 
be defined in \S \ref{subsec:
  adeles case 1}. 
The
group $\{a_t\}$ defined in \eqref{eq: def at} is contained in
$\SL_n(\bR)$ and thus in $\SL_n(\bA)$, and hence acts on
$\XXnA$, and the map $\pi$ is $\{a_t\}$-equivariant. 
We will lift $\sro$ to a
cross-section $\srotilde = \pi^{-1}(\sro)$ in $\XXnA$, and
define relevant $a_t$-invariant measures 
$ \mu$ on $\XXnA$ which descend under $\pi$ to the correct measure on 
$\crly{X}_n$ according to the case at hand. Namely, in Case I, the
measure $\mu = m_{\XXnA}$ is the unique $\SL_n(\bA)$-invariant
probability 
measure on $\XXnA$, and in Case II we will take measures $\mu =
\wt{m}_{\vec{\alpha}} = m_{\tilde{L}_{\vec{\al}} \tilde{y}_{\vec{\al}}}$
corresponding to $\vec{\alpha}$, which are homogeneous measures
supported on a compact adelic torus-orbit. These measures satisfy $\pi_*
m_{\XXnA} = m_{\XX_n}$ and $\pi_* \widetilde m_{\vec{\alpha}} =
m_{\vec{\alpha}}$. Using the results of \S\ref{sec:lifting}, we will
obtain that the
lifted cross-section $\srotilde \df \pi^{-1}(\sro)$ is $\widetilde
\mu$-reasonable in both cases. Similarly we will obtain that the lifted subsets
$\widetilde{\cB} \df \pi^{-1}(\cB), \, 
\widetilde{\cS}_{\vre} \df \pi^{-1}(\cS_{\vre})$ are
$\widetilde{\mu}$-JM.
Throughout the discussion we will take care to obtain explicit
description of the maps and measures that arise. 

In the subsection below we
will introduce the spaces and measures $(X, \mu)$, and in the
subsequent sections we will introduce the cross-sections and their JM-subsets. 


\subsection{The real homogeneous space, Case I}\label{sec: rhc1}
We \index{Case I (Haar-Siegel measure)} will work with the space of lattices $\XX_n \df
\SL_n(\R)/\SL_n(\Z)$, and the Haar-Siegel 
measure $m_{\XX_n}$. Let $\SL_n^{(\pm)}(\R)$ denote the $n \times
n$ real matrices of determinant $\pm 1$. The following simple
observation will be useful.

\begin{proposition}\label{prop: always confused}
$\XX_n$ is isomorphic to the quotient  $
\SL^{(\pm)}_n(\R)/\SL^{(\pm)}_n(\Z)$, via a map which is
$\SL_n(\R)$-equivariant. 
\end{proposition}

\begin{proof}
  Let $\tau: \SL_n(\R) \hookrightarrow \SL^{(\pm)}_n(\R)$ be the
  embedding. Then $\tau(\SL_n(\Z)) \subset \SL^{(\pm)}_n(\Z)$ and
  hence $\tau$ induces an $\SL_n(\R)$-equivariant map $\bar \tau:
  \SL_n(\R)/\SL_n(\Z) \to 
  \SL^{(\pm)}_n(\R)/ \SL^{(\pm)}_n(\Z)$. We leave it to the reader to
  verify that $\bar \tau$ is a bijection.
      \end{proof}

Note that elements of $\SL_n^{(\pm)}(\R)$
act on $\XX_n$ via their action on $\R^n$;
in terms of the isomorphism given in Proposition \ref{prop: always
  confused}, the left-action of
$\SL_n^{(\pm)}(\R)$ by matrix multiplication on lattices, is given by
left multiplication on cosets. 
 
\subsubsection{Contracting horospherical group} Let $G$ be an lcsc
group, $\{a_t\}\subset G$ a one-parameter subgroup, and $\Gamma<G$ a
lattice. Let $\mathscr{X} = G/\Gamma$ and let 
$\mu$ be an $\{a_t\}$-invariant measure on $\mathscr{X}$. The group
\begin{equation}\label{eq: the contracting group}
H^- \df \{g \in G: a_t g a_{-t} \to_{t\to +\infty} e\},
\end{equation}
is known 
as \index{H@$H^-$ -- contracting horospherical subgroup} the {\em
  contracting 
  horospherical subgroup} of $G$, corresponding to 
$\{a_t\}$. Also we denote the centralizer of $\{a_t\}$ by $H^0$, that
is,
\begin{equation}\label{eq: the centralizer}
H^0 \df \{g \in G: \forall t, a_t g = g a_t\}.
\end{equation}
We \index{H@$H^0$ -- central subgroup} will need the following
well-known fact. We leave the proof to the 
reader. 
\begin{proposition}\label{prop: horospherical and generic} The product
  $H^{\leq} =
  H^0H^-$ is a \index{H@$H^\leq$ -- non-expanding subgroup} group. 
If $x_0$ is $(a_t, \mu)$-generic and 
$h^0 \in H^0, \, h^- \in H^-$, then
$h^0h^-x_0$ is $(a_t, h^0_*\mu)$-generic.  
\end{proposition}
\ignore{
\begin{proof}
Let $f \in C_c(\mathscr{X})$ and write $\tilde{f} (x) \df
f(h^0x)$, so that
$\int f \, dh^0_* \mu = \int \tilde{f} \, d\,\mu.$ Let $\vre>0$. Then
$\tilde{f} \in C_c(X)$, and since $x_0$ is $(a_t, \mu)$-generic, there
is $T_0$ such that for all $T>T_0$ we have
$$
\left| \frac{1}{T}\int_0^T \tilde{f}(a_tx_0) dt - \int
\tilde{f} \, d\mu \right| < \frac{\vre}{2}.
$$
Since $\tilde{f}$ is compactly supported, there is a
neighborhood $W$ of $e $ in $G$ such that for all $w\in W$ and all $x
\in \mathscr{X}$, 
$|\tilde{f}(wx) - \tilde{f}(x)|< \frac{\vre}{2}$. Let $T_0> 2
\|f\|_\infty$, then for all $T> T_0$ we have  
  \end{proof}
}
    \subsection{The real homogeneous space, Case II}\label{subsec:
      case II}
    Let $d \geq 2$, and let $\bK =\on{span}_\bQ\set{\al_1,\dots,
      \al_d,1}$ be a totally real field of degree $n$ over $\bQ$. 
      Let
    $\sigma_1, \ldots, \sigma_n: \bK \to \R$ 
denote the distinct field embeddings. Departing
slightly from common 
conventions, we let
$\sigma_{n} = \mathrm{Id}$. Let
\begin{equation}\label{eq: all about alpha}
\vec{\alpha} \df
  \left( \begin{matrix} \alpha_1 \\ \vdots \\ \alpha_d \end{matrix}
  \right) \in \R^d.
\end{equation}
We point out a slight abuse of notation: what we now denote by 
$\vec{\alpha}$ is a column vector whereas in the introduction the same
notation was used for a row vector with the same entries. In this section it will be more
convenient to work with column vectors, and this should cause no
confusion.

  For $\beta \in \bK$, let $\pmb{\sig}: \bK
    \to \R^n$ be the $\Q$-linear map defined by 
  \begin{equation}\label{eq: geometric embedding}    
\pmb{\sig}(\beta) \df \left( \begin{matrix} \sigma_1(\beta)
    \\ 
\vdots \\ \sigma_{n-1}(\beta) \\
    \beta\end{matrix} \right ) \in \R^n. 
\end{equation}
We will refer to $\pmb{\sig}(\beta)$ as \index{G@$\pmb{\sig}(\beta)$
  -- geometric embedding of a totally real number field} the {\em
  geometric embedding} 
of  $\beta$. A lattice $\Lam\in \crly{X}_n$ is called \textit{of type}
$(\pmb{\sig}, \bK)$ \index{type $(\pmb{\sig}, \bK)$ -- lattice arising
  from embeddings of a number field} 
if there exist a basis
$\be_1,\dots, \be_n$ of $\bK$ 
over 
$\bQ$
 such that $\Lam$ is homothetic to the lattice 
 \begin{equation}\label{eq: homothetic to}
  \left\{ \pmb{\sig}(\beta) :  \beta \in
    \on{span}_\bZ\set{\be_1,\dots,\be_n} \right\}  =  
\mat{| & \dots &| \\ \pmb{\sig}(\be_1)& \dots & \pmb{\sig}(\be_n) \\ |
  & \dots & |}\bZ^n.
\end{equation}
Let $A \subset \SL_n(\R)$ \index{A@$A$ -- the diagonal group} denote
the group of diagonal matrices with 
positive diagonal entries. We say that an orbit $A\Lam' \subset \XX_n$ is {\em an orbit of
  type $(\pmb{\sig}, \bK) $} if there is 
  $\Lam\in \XX_n$ of type $(\pmb{\sig}, \bK)$ such
that 
$A\Lam'=A\Lam$.

\begin{lemma}\label{lem: addition compact orbits}
  Orbits of type $(\pmb{\sig}, \bK)$ are compact, and all
  compact $A$-orbits are of type $(\pmb{\sig}, \bK)$, for
  some totally real number field $\bK$ and some $\pmb{\sig}$. 
If 
$\Lam',\Lam\in \XX_n$ are lattices such that 
$\Lam'$ is of
type $(\pmb{\sig}, \bK)$, and
 $\on{Stab}_A(\Lam)  = 
 \on{Stab}_A(\Lam')$, then 
$A\Lam$ is an orbit of type $(\pmb{\sig}, \bK)$ (for the same field $\bK$).
\end{lemma}
\begin{proof}
  For the first two assertions see \cite[\S 6]{LW}. For the
  third assertion, let $\Del \df \on{Stab}_A(
  \Lam')$. It is well-known that the
$\bQ$-linear span of $\Del$, in the linear space of $n \times n$ real matrices, satisfies 
$$\wt{\bK}:= \on{span}_\bQ(\Del) = \set{\diag{\sig_1(\al),\dots,
    \sig_n(\al)}:\al\in \bK}.$$ 
This follows from the fact that any finite index subgroup of the group of units in the ring of integers of $\bK$ spans $\bK$. 
For $\al\in \bK$ let us denote $\tilde{\al} := \diag{\sig_1(\al),\dots, \sig_n(\al)}$. Fix a vector
$\mb{w}\in \Lam$ all of whose coordinates are strictly positive. Consider the map
$$\iota  :\bK\to \bR^n,\quad \iota(\al) \df  \tilde{\al} \mb{w}.$$
Note that this is an injective $\bQ$-linear map. We claim that its
image is the $\bQ$-span of
$\Lam$, which we denote by $\Lam_\bQ$. Indeed,
any $\tilde{\al}\in
\wt{\bK}$ is a linear combination over $\bQ$ of elements  
from $\Del$, and this implies that $\iota(\bK)\subset \Lam_\bQ$. Since
$\bK$ has dimension $n$ over $\bQ$ we conclude that   
$\iota$ is a linear isomorphism between $\bK$ and $\Lam_\bQ$. 

Now let $M = \iota^{-1}(\Lam)$ and let $\be_1,\dots,\be _n$ be a basis
of $M$ as a $\bZ$-module. It is clear from the construction that
$\Lam$ is obtained from $\pmb{\sig}(M)$ by applying the diagonal
matrix 
whose diagonal entries are the coordinates of $\mb{w}$. This shows
that $A \Lam = A \Lam''$, where $\Lam''$ is the lattice of type
$(\pmb{\sig}, \bK)$  obtained from
$\pmb{\sig}(M)$ by normalizing its covolume to 
be one.
\end{proof}

Let $\vec{\al}$ be as in \eqref{eq: all about alpha}. We set
\begin{equation}\label{eq: def g alpha}
g_{\vec{\al}} \df \mat{
\vline &\dots&\vline& \vline \\
\pmb{\sig}(\al_1)&\dots&\pmb{\sig}(\al_d)&\pmb{\sig}(1)\\
\vline &\dots&\vline&\vline \\
}, \end{equation}
so the bottom row of $g_{\vec{\al}}$ is $(\vec{\al}^{\mathbf{t}},1)$, where
$\vec{\al}^{\mathbf{t}}$ denotes the transpose of $\vec{\al}$. It is well-known
that the matrix $g_{\vec{\alpha}}$ is invertible.
Let 
\begin{equation}\label{eq:x alpha}
x_{\vec{\al}}\defi \av{\det g_{\vec{\al}}}^{-1/n} g_{\vec{\al}}\bZ^n\in \XX_n.
\end{equation}
For $M \in \GL_n(\R)$ we denote by $M^* =(M^{-1})^{\mathbf{t}}$ the inverse of
the transpose of $M$. 
For a lattice $x = g\bZ^n \in \XX_n$, the dual lattice is defined by $x^*\defi g^*\bZ^n$.
We then have that 
\begin{equation}\label{eq: def lattice with compact orbit}
x_{\vec{\al}}^* = cg_{\vec{\alpha}}^* \Z^n,
\end{equation}
where $c>0$ is \index{X@$x_{\vec{\al}}^* $ -- lattice
  with compact $A$-orbit} 
chosen so that $x_{\vec{\al}}^* \in \XX_n$.

Note that the one-parameter group $\{a_t\}$ defined in \eqref{eq: def at} is
contained in $A$. 
 We will need the following well-known fact:
\begin{proposition}
  \label{prop: compact A orbit}
The orbit $Ax_{\vec{\al}}^*$ is of type $(\pmb{\sig}, \bK)$, and $\{a_t\}$
acts uniquely ergodically on  
$A x_{\vec{\al}}^*$. For any $\Lam \in Ax_{\vec{\al}}^*,$ and any $v \in \Lam \sm
\{0\}$, all the coordinates of $v$ are nonzero. 
  \end{proposition}
\begin{proof}
Let
$c' \df |\det(g_{\vec{\alpha}})|^{-1/n}>0$, so that $x_{\vec{\al}} \df c' 
g_{\vec{\alpha}} \Z^n \in \XX_n$. Then by Lemma \ref{lem: addition
  compact orbits}, $A x_{\vec{\al}} \subset \XX_n$ is a compact orbit. The
map $M \mapsto M^*$ is a continuous group automorphism of $\SL_n(\R)$ 
which maps the groups $A$ and $\SL_n(\Z)$ to
themselves. Therefore it induces an automorphism $\Psi$ of
$\XX_n$, such that $x_{\vec{\al}}^* = \Psi(x_{\vec{\al}})$. Since $a^* = a^{-1}$
for $a \in A$, the stabilizers
satisfy 
$\on{Stab}_A(x_{\vec{\al}}^*)  = \on{Stab}_A(x_{\vec{\al}})$, and Lemma \ref{lem: addition
  compact orbits} implies that $Ax_{\vec{\al}}^*$ is 
an orbit of type $(\pmb{\sig},\bK)$. 

We now prove the second assertion.
The group $A$ is isomorphic (as a Lie group) to $\R^d$, and we can
realize this group isomorphism explicitly using the exponential map
$\mathrm{Lie}(A) \to A$. The orbit $Ax_{\vec{\al}}^*$ is
isomorphic to $\bT = \R^d/\Delta$ for some lattice $\Delta$ in
$\R^d$, and the action $a_t \curvearrowright Ax_{\vec{\al}}^*$ is mapped by this
isomorphism to a straightline flow on $\bT$, that is, a flow of the
form
$$
t P(x) = P\left(x + t \vec{\ell} \right), \ \ \text{ where } \vec\ell \in \R^d \sm
\{0\} \text{ and }
 P: \R^d \to \bT$$
is the projection. Recall that such a 
straightline flow is uniquely ergodic unless the straightline orbit of
$ 0 \in \bT$ is
contained in a proper subtorus of $\bT$.
Let $\bL$ be the Galois closure of $\bK/\Q$ and let $\mathscr{G}
= \mathrm{Gal}(\bL/\Q)$ denote the corresponding Galois group. Then
$\mathscr{G}$ acts transitively on the field embeddings $\sigma_1, \ldots, \sigma_n$
by post-composition. This gives an identification of $\mathscr{G}$
with a transitive subgroup of the group  
$S_n$ of permutations of $\{1, \ldots,
n\}$. In turn, this allows $\mathscr{G}$ to act on the group of diagonal matrices by permuting 
the coordinates on the diagonal. As the discussion in
\cite[\S6]{LW} shows (in particular, from \cite[Step 6.1]{LW}),  for $N \subset \R^d \cong
\mathrm{Lie}(A)$, $P(N)$ is a compact subtorus
of $\bT$ if and only 
if $A_0 \df \exp(N) \subset A$ is $\mathscr{G}$-invariant. Thus if $\{a_t\} \subset
A_0$, then $A_0$  contains any group obtained from $\{a_t\}$ by 
acting on it with $\mathscr{G}$. Because of the transitivity mentioned above, $A_0$ must contain all
the subgroups
$$\left\{a^i_t: t \in \R \right\}, \  \ \text{ where } a^i_t = \diag{e^{t},
\ldots, e^{t}, \underset{i\text{th position}}{e^{-dt}}, e^t, \ldots, e^t}. $$
Since the groups $\left\{a^i_t\right \}$ generate $A$, we must have $A_0=A$, and
this establishes unique ergodicity.

For the third assertion, note that if $\pmb{\sig}(\beta)$ has one of
its coordinates equal to zero, then $\beta =0$. This observation,
along with \eqref{eq: homothetic to}, 
implies the third assertion for $\Lam = x_{\vec{\al}}^*$. The statement now
follows for general $\Lam =
ax_{\vec{\al}}^*$ by the definition of the $A$-action. 
\end{proof}

We will denote the $A$-invariant probability measure
\index{M@$m_{Ax_{\vec{\al}}^*}$ -- $A$-homogeneous measure on orbit of
  $x_{\vec{\al}}^*$} on
$A x_{\vec{\al}}^*$ by $m_{Ax_{\vec{\al}}^*}$.

\begin{proposition}
  \label{prop: explicit compact push}
  Let \index{L@$\Lam_\theta$ -- lattice corresponding to $\theta \in \R^d$}
  \begin{equation}\label{eq: def lattice to a vector}
\Lambda_{\vec{\alpha}} \df
\left( \begin{matrix} I_d & - \vec{\alpha} \\ \mathbf{0}^{\mathbf{t}} &
    1\end{matrix} \right) \Z^n \in \XX_n,
\end{equation}
where $\mathbf{0} \in \R^d$ is the zero (column) vector,
and let 
 $$ B_{\vec{\al}}  \df (b_{ij})_{i,j =1, \ldots, d} \ \text{ where } \ 
b_{ij} \df \sigma_j(\alpha_i) -\alpha_i. 
$$
Then $B_{\vec{\al}}$ is invertible, and for $c_1 \df \left|\det
  \left(B_{\vec{\al}}\right)\right|^{-1/n},$ the matrix  \index{H@$\bar h_{\vec{\al}}$ -- the
  matrix conjugating $A$ to $\bar A_{\vec{\al}}$}
\begin{equation}\label{eq: def Bar B}\bar{h}_{\vec{\al}} \df c_1
  \left( \begin{matrix} B_{\vec{\al}} & \mathbf{0} \\ \mathbf{0}^{\mathbf{t}} & 
      1\end{matrix} \right) \in \SL_n^{(\pm)}(\R)
  \end{equation}
satisfies that the 
trajectory $\{a_t \Lambda_{\vec{\alpha}}: t> 0\}$ is generic for the
measure\index{M@$m_{\vec{\alpha}}$ -- measure corresponding to
  approximation of algebraic $\vec{\alpha}$}
\begin{equation}\label{eq: def m alpha}
  m_{\vec{\alpha}} \df 
  (\bar{h}_{\vec{\al}})_* m_{Ax_{\vec{\al}}^*}.
  \end{equation}
  \end{proposition} 

  \begin{proof}
    Note that
   $$
   B_{\vec{\al}}^{\mathbf{t}} = \left( \sigma_i(\alpha_j) - 
     \alpha_j 
   \right)_{i,j=1, \ldots, d}.
$$
Also note that the entries of the right-most column of \eqref{eq: def g
  alpha} are all equal to 1, and the first $d$ entries of the bottom
row in  \eqref{eq: def g
  alpha}  are $\vec{\alpha}^{\mathbf{t}}$. Thus, letting $\mathbf{1} \in \R^d $
be the column vector all of whose entries are 1, we find 
\begin{equation}\label{eq: furthermore}
g_{\vec{\alpha}} = \left( \begin{matrix}
    B_{\vec{\al}}^{\mathbf{t}}  & \mathbf{1} \\
    \mathbf{0}^{\mathbf{t}} & 1  \end{matrix} \right) \, \left( \begin{matrix} I_d & \mathbf{0} \\
    \vec{\alpha}^{\mathbf{t}} & 1  \end{matrix} \right) ,
\end{equation}
which implies $\det(B_{\vec{\al}}) \neq 0$, and thus $\det\left(\bar{h}_{\vec{\al}} \right) = \pm 1.$
  
From \eqref{eq: furthermore} we have that
\begin{equation}\label{eq: for det}
  g_{\vec{\alpha}}^* = \left( \begin{matrix}
    B_{\vec{\al}}^{\mathbf{t}}  & \mathbf{1} \\
    \mathbf{0}^{\mathbf{t}} & 1  \end{matrix} \right)^* \, \left( \begin{matrix}
    I_d & -\vec \alpha  \\
   \mathbf{0}^{\mathbf{t}} & 1  \end{matrix} \right) = \left( \begin{matrix}
    B_{\vec{\al}}^{-1} & \mathbf{0} \\
    \mathbf{0}^{\mathbf{t}} & 1  \end{matrix} \right) q \left( \begin{matrix}
    I_d & -\vec \alpha \\
    \mathbf{0}^{\mathbf{t}} & 1  \end{matrix} \right),
\end{equation}
where $q$ is of the form
$$
q = \left( \begin{matrix} I_d & \mathbf{0} \\ \mathbf{x}^{\mathbf{t}} & 1 \end{matrix}
\right), \text{ for some } \mathbf{x}\in \R^d.
$$
In particular we have $\lim_{t \to \infty} a_t q a_{-t} = e$; i.e., $q
\in H^-$, 
the contracting horospherical subgroup of $\{a_t\}$.

Let $x_{\vec{\al}}^* = c_1 g_{\vec{\alpha}}^* \Z^n$. By \eqref{eq: for det}, $c_1 =
\left| \det \left( g_{\vec{\alpha}}^*  
  \right) \right|^{-1/n}$, and thus $x_{\vec{\al}}^*$ is the lattice as in Proposition \ref{prop:
  compact A orbit} and satisfies 
\begin{equation}\label{eq:stable relation}
x_{\vec{\al}}^* =
(\bar{h}_{\vec{\al}})^{-1} 
q \left( \begin{matrix}
    I_d & -\vec \alpha  \\
    \mathbf{0}^{\mathbf{t}} & 1  \end{matrix} \right) \Z^n = (\bar
{h}_{\vec{\al}})^{-1} q \Lambda_{\vec{\alpha}}.
\end{equation}
The lattice $x_{\vec{\al}}^*$ is $(a_t, m_{Ax_{\vec{\al}}^*})$-generic by
Proposition \ref{prop: compact A orbit}. Since 
$\Lambda_{\vec{\alpha}} = q^{-1} \bar{h}_{\vec{\al}} x_{\vec{\al}}^*$, $\bar h
$
commutes with the $\{a_t\}$ action, and $q^{-1}$ belongs to the contracting 
horospherical group for $\{a_t\}$, we have by Proposition \ref{prop:
  horospherical and generic} that $\Lambda_\alpha$ is $(a_t,
m_{\vec{\alpha}})$-generic. 
%
\end{proof}

\subsection{The adelic homogeneous space, Case I}\label{subsec:
  adeles case 1} 
We briefly recall facts and notation regarding the rational adeles.
See \cite{Weil_adeles, Platonov_Rapinchuk} for more details on adeles
and arithmetic groups, and 
see \cite{Guilloux} for a gentle recent introduction. 
Let $\mb{P}$ be the set of (rational) primes. Let \index{A@$\bA_f$ --
  finite adeles} \index{A@$\bA$ --
  adeles} $\bA = \bR\times \bA_f =
\bR\times\prod_{p\in \mb{P}}' \bQ_p$ be the ring of adeles. Here $\prod'$
stands for the restricted product --- that is, 
a sequence
$\underline \beta = (\beta_\infty,\beta_f) =
(\beta_\infty,\beta_2,\beta_3,\dots,\beta_p,\dots)$ belongs to  $\bA$ if and
only if $\beta_p\in \bZ_p$ for all but finitely  
many $p$. As suggested by the notation, we denote the real coordinate
of a sequence $\underline \beta\in \bA$ 
by $\beta_\infty$ and the sequence of $p$-adic coordinates by $\beta_f =
(\beta_p)_{p\in \mb{P}}$. The rational numbers $\Q$ are embedded in $\bA$ diagonally,
that is, $q \in \Q$ is identified with the constant sequence $(q,q,
\cdots).$ We let $\SL_n(\bA) = \SL_n(\bR)\times
\SL_n(\bA_f) = \SL_n(\bR)\times \prod_{p\in \mb{P}}'\SL_n(\bQ_p)$ and
use similar notation $(g_\infty,g_f) = (g_\infty, (g_p)_{p\in
  \mb{P}})$ to denote elements of $\SL_n(\bA)$. It is well-known that
the diagonal embedding of $\SL_n(\bQ)$ in $\SL_n(\bA)$ is a lattice in
$\SL_n(\bA)$.  Let \index{K@$K_f$ -- product of $\SL_d(\bZ_p)$'s, adelic case}
\begin{equation}\label{eq: def K_f}
K_f \df \prod_{p\in
  \mb{P}}\SL_n(\bZ_p)
  \end{equation}
and \index{P@$\pi_f : \SL_n(\bA) \to \SL_n(\bA_f)$ -- projection to
  finite places}
$$\pi_f: \SL_n(\bA) \to \SL_n(\bA_f), \ \ \ 
\pi_f(g_\infty, g_f) \df g_f.$$
Then $K_f$ is a compact open subgroup of
$\SL_n(\bA_f)$. Via the embedding $\SL_n(\bA_f) \cong \{e\} \times
\SL_n(\bA_f)$ we also think of $K_f$ as a subgroup of
$\SL_n(\bA)$. 
We shall use the following  
two basic facts (see \cite[Chap. 7]{Platonov_Rapinchuk}):
\begin{enumerate}[(i)]
\item\label{fact1} The intersection $K_f \cap \pi_f(\SL_n(\bQ))$ is equal
  to $\pi_f(\SL_n(\bZ))$.
\item\label{fact2} The projection $\pi_f(\SL_n(\bQ))$ is dense in $\SL_n(\bA_f)$.
\end{enumerate}
 
Let \index{X@$\XXnA= \SL_n(\bA)/\SL_n(\bQ)$ -- adelic space}
$$\XXnA= \SL_n(\bA)/\SL_n(\bQ),$$
and let \index{M@$m_{\XXnA}$ -- natural measure on $\XXnA$}
$m_{\XXnA}$ denote the $\SL_n(\bA)$-invariant probability 
measure on $\XXnA$. 
There is a natural projection \index{P@$\pi:\XXnA\to\XX_n$ --
  projection from adelic space to real space} $\pi:\XXnA\to\XX_n$
which we now 
describe in two equivalent ways. 

\textbf{First definition of} $\pi$: Given $\tilde{x} =
(g_\infty,g_f)\SL_n(\Q)\in \XXnA$, using \eqref{fact2} and the fact
that $K_f$ is open, 
we may replace the representative
$(g_\infty,g_f)$ by another $(g_\infty\ga,g_f\ga)$, where
$\ga\in\SL_n(\bQ)$ is such that $g_f\ga\in K_f$. We then define 
$\pi(\tilde{x}) = g_\infty\ga \SL_n(\bZ)$. This is well-defined
since, if $g_f\ga_1, g_f\ga_2\in K_f$, then by \eqref{fact1},
$\pi_f(\ga_1^{-1}\ga_2 ) \in K_f\cap
\pi_f(\SL_n(\bQ))= \pi_f(\SL_n(\bZ))$, 
and so $g_\infty\ga_1\SL_n(\bZ) = g_\infty\ga_2\SL_n(\bZ)$.  

\textbf{Second definition of} $\pi$: View $K_f$ as a subgroup of  $\SL_n(\bA)$. 
We claim that 
%
we may 
identify the double coset space $K_f\backslash \SL_n(\bA)/\SL_n(\Q)$ with
$\SL_n(\bR)/\SL_n(\bZ)$. Indeed, by \eqref{fact2} and since $K_f$ is
open (as a subgroup of $\SL_n(\bA_f)$), each double coset $K_f(g_\infty,g_f)\SL_n(\Q)$
contains representatives with $g_f= e_f$ (the identity element in
$\SL_n(\bA_f)$). The real coordinates of
all such  representatives form a single  
left coset of $K_f\cap \SL_n(\bQ) = \SL_n(\bZ)$. With this
identification $\pi$ is simply the projection from the  
coset space $\SL_n(\bA)/\SL_n(\Q)$ to the double coset space $K_f\backslash
\SL_n(\bA)/\SL_n(\Q)$.  

We leave it to the reader to check that these two
definitions agree and that $\pi$ intertwines the actions of
$G_\infty \df \SL_n(\bR)$ on $\XXnA, \XX_n$. Since there is a unique
$G_\infty$-invariant probability measure on $\XX_n$, we have that $\pi_* m_{\XXnA}
= m_{\XX_n}$. In particular, the
1-parameter 
group $\{a_t\} \subset G_\infty$ acts on both of these spaces and
$\pi$ is a factor map for these actions. 
The following standard statement will be important for us:
\begin{lemma}\label{lem: mautner}
The group $\{a_t\}$ acts ergodically on $\left(\XXnA,m_{\XXnA}\right)$. 
\end{lemma}
\begin{proof}
By the 
Mautner phenomenon (see e.g.\,\cite{EW}), it is enough to show that
$G_\infty$ acts ergodically on $\left(\XXnA,m_{\XXnA}\right)$. By
duality, this is equivalent to the ergodicity of the action by right translations, of
$\pi_f(\SL_n(\Q))$ on $G_\infty \backslash \SL_n(\bA) =
\SL_n(\bA_f)$. Since the stabilizer of a measure is a closed 
group, by 
property \eqref{fact2} above, any $\pi_f(\SL_n(\Q))$-invariant measure
must be $\SL_n(\bA_f)$-invariant, and thus is the Haar measure of
$\SL_n(\bA_f)$. In particular the action of
$\pi_f(\SL_n(\Q))$ is uniquely ergodic, and hence ergodic. 
\end{proof} 
\subsubsection{The adelic homogeneous space, Case II} \label{subsec:
  adeles case 2}
Let $A$ be the diagonal group, let $x_{\vec{\al}}^*$ be
as in \eqref{eq: def lattice with compact orbit},  let $\bar h_{\vec{\al}}$ be
as in \eqref{eq: def Bar B} and let \index{A@$\bar A_{\vec{\al}}$ -- conjugate of
  $A$ such that $m_{\vec{\alpha}}$ is $\bar A_{\vec{\al}}$-homogeneous}
  \index{Y@$y_{\vec{\al}} $ -- lattice
  with compact $\bar A_{\vec{\al}}$-orbit}
\begin{equation}\label{eq: def bar A}
  \bar A_{\vec{\al}} \df \bar h_{\vec{\al}} A \bar h_{\vec{\al}}^{-1}, \ \ \ \ 
  y_{\vec{\al}} \df \bar h_{\vec{\al}} x_{\vec{\al}}^*.
  \end{equation}
The measure $m_{\vec{\alpha}}$
defined in \eqref{eq: def m alpha} is $\bar A_{\vec{\al}}$-homogeneous, more precisely, by
Proposition \ref{prop: compact A orbit}, we have that the orbit 
$\bar A_{\vec{\al}}y_{\vec{\al}}$ is compact and 
$m_{\vec{\alpha}}$ is the $\bar A_{\vec{\al}}$-invariant probability measure on $\bar A_{\vec{\al}}
y_{\vec{\al}}$.

Let
\begin{equation}\label{eq: def bar A Lam}
  \on{Stab}_{\bar A_{\vec{\al}}}(y_{\vec{\al}}) \df \set{\bar{a} \in
    \bar A_{\vec{\al}}: \bar{a} y_{\vec{\al}} = y_{\vec{\al}}} 
  \end{equation} 
  be the stabilizer group of $y_{\vec{\al}}$ in $\bar A_{\vec{\al}}$. 
  Write 
  $y_{\vec{\al}} = g_\infty\Z^n$ for some $g_\infty \in \SL_n(\R)$,
  so that 
  $
g_\infty^{-1} \on{Stab}_{\bar{A}_{\vec{\al}}}({y}_{\vec{\al}})
g_\infty \subset \SL_n(\Z)  
$
is cocompact in the conjugated group $g_\infty^{-1} \bar A_{\vec{\al}} g_\infty.$ 
Let $\Delta$ be the diagonal embedding of $\SL_n(\Z)$ in $\SL_n(\bA)$,
 let \index{M@$M_{\vec{\al}}$ -- the finite component of the adelic
   torus $\tilde{L}_{\vec{\al}}$} $M_{\vec{\al}} 
 \subset K_f$ denote the closure
of 
$\pi_f \circ \Delta \left( g_\infty^{-1}
  \on{Stab}_{\bar{A}_{\vec{\al}}}({y}_{\vec{\al}}) g_\infty
\right),$ and 
let \index{L@$\tilde{L}_{\vec{\al}} =\bar A_{\vec{\al}}\times
  M_{\vec{\al}}$ -- the adelic lifted torus} 
$$\tilde{L}_{\vec{\al}} \df \bar A_{\vec{\al}} \times M_{\vec{\al}} \subset \SL_n(\bA).$$ 
For $\bar{a} \in
  \on{Stab}_{\bar{A}_{\vec{\al}}}({y}_{\vec{\al}})$ we let
$$\ga_{\bar{a}} \df 
\pi_f \circ \Delta (g_\infty^{-1} \bar{a} g_\infty ) \in M_{\vec{\al}},$$
and
\begin{equation}\label{eq:defining lambda tilde}
\tilde{y}_{\vec{\al}} \df (g_\infty,e_f)\SL_n(\Q)\in \XXnA.
\end{equation}
\index{Y@$\tilde{y}_{\vec{\al}}$ -- adelic lift of $y_{\vec{\al}}$. Point with compact $\tilde{L}_{\vec{\al}}$-orbit.}

Note that $M_{\vec{\al}}$ is a compact abelian group. Note also that the group
$M_{\vec{\al}}$, the homomorphism $\bar{a} \mapsto \ga_{\bar{a}}$, and the point
$\tilde{y}_{\vec{\al}}$ all depend on the choice of the representative  
$g_\infty$ of $y_{\vec{\al}}$. This dependence will not matter to us and we
suppress it from the notation.

\begin{proposition} \label{prop: lifted measure case II}
  The orbit
  $\tilde{L}_{\vec{\al}}\tilde{y}_{\vec{\al}}\subset \XXnA$ is compact, and supports a
  finite $\tilde{L}_{\vec{\al}}$-invariant
  measure\index{M@$m_{\tilde{L}_{\vec{\al}}\tilde{y}_{\vec{\al}}}$ -- the natural measure
    on the adelic lifted torus orbit $\tilde{L}_{\vec{\al}}\tilde{y}_{\vec{\al}}$}
  $m_{\tilde{L}_{\vec{\al}}\tilde{y}_{\vec{\al}}}$. The 
  action of $\{a_t\}$ on $\left(\tilde{L}_{\vec{\al}}  
   \tilde{y}_{\vec{\al}}, m_{\tilde{L}_{\vec{\al}}\tilde{y}_{\vec{\al}}} \right)$ is 
  uniquely ergodic, and 
$\pi_*\left(m_{\tilde{L}_{\vec{\al}}\tilde{y}_{\vec{\al}}}\right)= m_{\vec{\alpha}}$. 
Moreover, $M_{\vec{\al}}$ acts transitively on the fibers of $\pi|_{\tilde{L}_{\vec{\al}}\tilde{y}_{\vec{\al}}}$.
\end{proposition}
\begin{proof}
In order to show that $\tilde{L}_{\vec{\al}} \tilde{y}_{\vec{\al}}$ is compact we need to
show that $\on{Stab}_{\tilde{L}_{\vec{\al}}}(\tilde{y}_{\vec{\al}})$ is a lattice 
in $\tilde{L}_{\vec{\al}}$. 
We claim that \begin{equation}\label{eq: inclusion}
\set{(\bar{a},\ga_{\bar{a}}):\bar{a}\in 
    \on{Stab}_{\bar{A}_{\vec{\al}}}({y}_{\vec{\al}})} \subset 
\on{Stab}_{\tilde{L}_{\vec{\al}}}(\tilde{y}_{\vec{\al}}).
\end{equation}
Indeed, for any $\bar{a}\in \on{Stab}_{\bar A_{\vec{\al}}}(y_{\vec{\al}})$ we have
\[\begin{split}
    (\bar{a},\ga_{\bar{a}}) \tilde{y}_{\vec{\al}} = &
    (\bar{a},\ga_{\bar{a}})(g_\infty,e_f)\SL_n(\Q) 
    =
(g_\infty g_\infty^{-1} \bar{a} g_\infty, \ga_{\bar{a}})\SL_n(\Q) \\  = &
(g_\infty, e_f) (g^{-1}_\infty \bar{a} g_\infty, \pi_f \circ \Delta
(g^{-1}_\infty \bar{a} g_\infty)) \SL_n(\Q) \\= &
(g_\infty, e_f)\SL_n(\Q) =
\tilde{y}_{\vec{\al}}.
\end{split}\]
We remark that the other inclusion in \eqref{eq: inclusion} is also true,
but we will not need it. 
Now, since $M_{\vec{\al}}$ is
compact and $\on{Stab}_{\bar A_{\vec{\al}}}({y}_{\vec{\al}})$ is a lattice  
in $\bar A_{\vec{\al}}$, the graph $ \set{(\bar{a},\ga_{\bar{a}}):\bar{a}\in \on{Stab}_{\bar A_{\vec{\al}}}({y}_{\vec{\al}})}$ is
a lattice in $\tilde L_{\vec{\al}}$.

Since $M_{\vec{\al}} \subset K_f$  
we have that $\pi(\tilde{L}_{\vec{\al}}\tilde{y}_{\vec{\al}}) = \bar A_{\vec{\al}} y_{\vec{\al}}$ and since $\pi$
intertwines the $\bar A_{\vec{\al}}$-action, we have that  
$\pi_*m_{\tilde{L}_{\vec{\al}}\tilde{y}_{\vec{\al}}}$ is an $\bar A_{\vec{\al}}$-invariant probability measure supported
on $\bar A_{\vec{\al}} y_{\vec{\al}}$. As $m_{\vec{\alpha}}$ is the unique such
measure, we have $\pi_*
m_{\tilde{L}_{\vec{\al}}\tilde{y}_{\vec{\al}}}= m_{\vec{\alpha}}$. 

It remains to establish the unique ergodicity of the action of $a_t$
on $\tilde{L}_{\vec{\al}} \tilde{y}_{\vec{\al}}$.
By standard facts on translation flows on compact abelian groups (see
e.g. \cite[Thm. 4.14]{EW}), this is equivalent to  
showing that any character on
$\tilde{L}_{\vec{\al}}/\on{Stab}_{\tilde{L}_{\vec{\al}}}(\tilde{y}_{\vec{\al}})$ which is trivial on
the 
image of $\{a_t\}$ is trivial. This is in turn equivalent to 
the fact that any character on $\tilde{L}_{\vec{\al}}$ which is trivial on $\{a_t\}$
and on $\on{Stab}_{\tilde{L}_{\vec{\al}}}(\tilde{y}_{\vec{\al}})$ is trivial. To this end, let  
$\chi: \tilde {L}_{\vec{\al}} \to \bS^1$ be a character such that 
$$\chi|_{\{a_t\}} \equiv
\chi|_{\on{Stab}_{\tilde{L}_{\vec{\al}}}(\tilde{y}_{\vec{\al}})} \equiv 1.$$
There are characters  
$\chi_1:\bar A_{\vec{\al}}\to \bS^1$, $\chi_2:M_{\vec{\al}}\to\bS^1$ such that
$$\chi(\bar{a},h) =
\chi_1(\bar{a}) \cdot \chi_2(h).$$

We claim that 
there exists $k \in \N$ such that 
$\chi_2^k$ is trivial. Let $\cW$ be a neighborhood of $1$ in $\bS^1$
which does not contain any nontrivial subgroups. The group $K_f$ has a collection
of clopen subgroups that 
give a basis of the topology at the identity, and thus the same is
true for $M_{\vec{\al}}$. By continuity of $\chi_2$
there is a clopen subgroup $M'$ of $M_{\vec{\al}}$ such that $\chi_2(M') \subset
\cW$, and hence $M' \subset \ker \chi_2$. This implies that $\chi_2$
factors through the finite quotient of $M_{\vec{\al}}/M'$, proving the claim.

It follows
that $\chi_1$ is trivial  
on $\on{Stab}_{\bar A_{\vec{\al}}}(y_{\vec{\al}})^k$ and on $\{a_t\}$ and
therefore induces a character on 
$\bar A_{\vec{\al}}/\on{Stab}_{\bar A_{\vec{\al}}}(y_{\vec{\al}})^k$ which is trivial on the image of $\{a_t\}$.   
Since $\{a_t\}$ acts ergodically on $\bar A_{\vec{\al}}/\on{Stab}_{\bar A_{\vec{\al}}}(y_{\vec{\al}})$
and $\bar A_{\vec{\al}}$ is connected, $\{a_t\}$ also acts 
ergodically on 
$\bar A_{\vec{\al}}/\on{Stab}_{\bar A_{\vec{\al}}}(y_{\vec{\al}})^k$.  Hence $\chi_1$ is trivial. This 
in turn implies that $\chi_2$ is trivial on $\set{\ga_{\bar{a}} : \bar{a} \in
  \on{Stab}_{\tilde{L}_{\vec{\al}}}(\tilde{y}_{\vec{\al}})}$, which is a  
dense subgroup of $M$. Therefore $\chi_2$ is also
trivial. 

We now establish the transitivity of the $M_{\vec{\al}}$-action on the fibers of $\pi|_{\tilde {L}_{\vec{\al}}\tilde{y}_{\vec{\al}}}$. We need to show that if 
$x,y\in \tilde {L}_{\vec{\al}}\tilde{y}_{\vec{\al}}$ are such that $\pi(x) = \pi(y)$, then there exists $m\in M_{\vec{\al}}$ such that $mx = y$. 
Since $\tilde {L}_{\vec{\al}}$ commutes with $\pi$, acts transitively, and is commutative we may assume without loss of generality that $x = \tilde{y}_{\vec{\al}}$.
\ignore{
 We first
claim that each fiber decomposes into finitely many $M_{\vec{\al}}$-orbits. Let $x\in \bar{A}_{\vec{\al}}y_{\vec{\al}}$ and consider the fiber
$F_x\def \pi^{-1}(x)\cap (\bar{A}_{\vec{\al}}\times M_{\vec{\al}})\tilde{y}_{\vec{\al}}$. 
Choose a neighborhood of the identity $e_\infty\in \bar{A}_{\vec{\al}}$, $U$, and consider the open cover
$$F_x\subset \bigcup_{y\in  F_x} (U\times M) y.$$
By compactness of $F_x$, there exists a finite sub-cover 
$$F_x\subset \bigcup_1^k (U\times M)y_i.$$
If we choose $U$ small enough, the map $u\mapsto uy$ is injective for all $y\in (\bar{A}_{\vec{\al}}\times M_{\vec{\al}}\tilde{y}_{\vec{\al}}$, and therefore we conclude that 
$F_x\subset \cup_1^k My_i$ and since $F_x$ is $M_{\vec{\al}}$-invariant we conclude that it decomposes into finitely many $M_{\vec{\al}}$-orbits. 
}
Recall that $\tilde{y}_{\vec{\al}} = (g_\infty,e_f)\Gaa$, where we set
$\Gaa \df \SL_n(\bQ)$. Since $\bar{A}_{\vec{\al}}\times M_{\vec{\al}}$ acts transitively,  we can write 
$$y = (\bar{a},m_f) \tilde{y}_{\vec{\al}} = (\bar{a}, m_f)(g_\infty, e_f) \Gaa =
(\bar{a}g_\infty, m_f)\Gaa.$$ 
Applying $\pi$ we see that 
$$y_{\vec{\al}}  = g_\infty \SL_n(\bZ) = \pi(x) = \pi(y) = \bar{a} g_\infty\SL_n(\bZ),$$
and therefore by definition of $M_{\vec{\al}}$, there exists $\ga\in M\cap
\SL_d(\bZ)$ such that $\bar{a}g_\infty = g_\infty \ga$. 
It follows that 
$$ y = (g_\infty\ga, m_f)\Gaa = (g_\infty, m_f\ga^{-1})\Gaa  =
(e_\infty,m_f\ga^{-1})(g_\infty,e_f)\Gaa \in M\tilde{y}_{\vec{\al}}.$$ 
\end{proof}

\section{Some  reasonable
  cross-sections}\label{sec: applications}
We will now define the cross-sections $\cS$ we will need for
our applications, as well as the cross-section measures
$\mu_{\cS}$. In this section we will work with
the real space $\XX_n$, and the adelic space $\XXnA$ will be discussed
in \S \ref{subsec: adelic lifts}. 

\subsection{The cross-section $\sro$}\label{sec:main example}
We introduce some convenient notation. Recall that
the set of primitive
vectors in $\Lambda$ is denoted by $\Lambda_{\prim}$.
Given a subset $W\subset \bR^n$ and $k\ge 1$ we let
\index{X@$\XX_n(W,k)$ -- lattices with $\#(\Lam_{\prim} \cap W)\geq
  k$}
\begin{align}
\label{def:special lattices}\XX_n(W,k) &\df  \set{\Lambda\in \XX_n:
                                         \#(\Lambda_{\prim} \cap W)\ge
                                         k}. 
\end{align}
For $k=1$ we will omit $k$ and denote 
\begin{align*}
\XX_n(W) &\df \XX_n(W,1).
\end{align*}
We will be interested in the case when there is a unique
primitive vector in $W$, and thus we let \index{X@$\XX_n^{\sharp}(W)$
  -- lattices with $\#(\Lam_{\prim} \cap W) = 1$} 
\begin{align*}
\XX_n^{\sharp}(W) & \df \XX_n(W)\smallsetminus \XX_n(W,2).
\end{align*}
There is a natural map \index{v@$v
\left(\Lambda\right)$ -- unique vector in $\Lambda_{\on{prim}}\cap W$}
\begin{equation}\label{eq: well defined vector}
v:\XX_n^{\sharp}(W)\to W, \  \  \textrm{defined by } \{v
\left(\Lambda\right)\} = \Lambda_{\on{prim}}\cap W. 
\end{equation} 
With this notation we have:
\begin{lemma}\label{lem:continuity of vector}
Let $W\subset \bR^n$ be a compact set, $V\subset W$ a relatively open
subset and $k\ge 1$ an integer. 
\begin{enumerate}
\item\label{eq:12111} The set $\XX_n(W, k)$
is closed in $\XX_n$. 
\item\label{eq:12112} The set $\XX_n^{\sharp}(W)\cap \XX_n(V)$ is open in $\XX_n(W)$.
\item\label{eq:12113} The map $v:\XX_n^{\sharp}(W)\to W$ is continuous.
\end{enumerate}
\end{lemma}
\begin{proof} 
Let $\Lambda_i \in \XX_n(W,k)$ such that $\Lambda_i
\to \Lambda$. We can choose $g_i \to h$ such that $\Lambda_i = g_i
\Z^n$ and $\Lambda = h\Z^n$. 
We need to show that $\Lambda$ contains at least $k$
primitive vectors in $W$. Since $g_i \to h$ and $W$ is compact, there is a compact
subset of $\R^n$ containing all of the sets $g_i^{-1}W$. Since $\Lambda_i
\in \XX_n(W,k)$, each $g_i^{-1}W$ contains at least $k$  distinct
elements of $\Z^n_{\prim}$. 
After passing to a
subsequence if necessary, there are 
distinct $\mb{p}_1, \ldots, \mb{p}_k \in \Z^n_{\prim}$, such
that $\mb{p}_j\in g_i^{-1}W$ for $j=1, \ldots, k$. Since $W$ is closed, 
their limits $ h\mb{p}_j$ belong to
$\Lambda_{\prim} \cap W$, and are distinct. This proves \eqref{eq:12111}.

The complement  of $\XX_n^{\sharp}(W)\cap \XX_n(V)$ in $\XX_n(W)$ consists
of lattices that either contain at least two distinct primitive
vectors in $W$ or else, contain a primitive vector 
in  $W\smallsetminus V$ (these cases are not mutually exclusive). That is, 
$$\XX_n(W)\smallsetminus(\XX_n^{\sharp}(W)\cap \XX_n(V)) = \XX_n(W,2) \cup
\XX_n(W\smallsetminus V),$$
which by 
part~\eqref{eq:12111} of the Lemma, is a union of two closed subsets
of $\XX_n$. This proves \eqref{eq:12112}.


Suppose $\Lambda_i = g_i \Z^n \to \Lambda = h\Z^n $ is a converging sequence
in $\XX_n^{\sharp}(W)$, with $g_i \to h$. As before, the sets $g_i^{-1}W$
are  all contained in  
a fixed compact set in $\bR^n$. Since $\Lambda_i$ and $\Lambda$ belong
to $\XX_n^{\sharp}(W),$
the sets $g_i^{-1}W \cap \bZ^n_{\prim}$ and $h^{-1}W \cap \bZ^n_{\prim}
$ are singletons. Passing to a  
subsequence if needed, we may assume that $g_i^{-1}W = \{\mb{p}\}$ for
a fixed  $\mb{p}$. Therefore
 $v(\Lambda_i) =g_i\mb{p} \to h\mb{p}  \in \Lambda_{\on{prim}}\cap W$ and
$h\mb{p} = v(\Lambda)$, proving \eqref{eq:12113}.
\end{proof}

The following Lemma is proved using similar arguments: 
\begin{lemma}\label{lem:open sets}
Let $W\subset \bR^n$ be an open set. For any $k\ge 1$, $\XX_n(W,k)$ is open in $\XX_n$.
\end{lemma}
\qed

Fix a norm $\| \cdot \|$ on $\R^d$, and denote the Lebesgue measure on
$\R^n$ by $m$. 
For $\bx = (x_1, \ldots, x_n) \in \R^n$, we refer to $\pi_{\R^d}(\bx)$
as the {\em horizontal component} of $\bx$ and to $x_n$ as the {\em
    vertical component} of $\bx$. For positive numbers $r$ and $s$,
  let \index{C@$C_r$ -- cylinder in $\R^n$}
\begin{equation}\label{eq: def C r}
C_r (s) \df \left\{\bx \in \R^n:  \left\| \pi_{\R^d} (\bx) \right\| \leq r,\ |x_n|
  \leq s \right\} \text{ and } C_r \df C_r(1).
\end{equation}
Note that this cylinder depends on the choice of the norm; we consider
the norm as fixed and thus it does not appear in the notation. Choose
\index{R@$r_0$ -- large enough radius for cylinder}
$r_0>0$ large enough so that
\begin{equation}\label{eq: ro satisfies}
  m(C_{r_0}) \geq 2^n.
\end{equation}
By Minkowski's convex body theorem, this
implies that for any $\Lambda \in \XX_n$, $\Lambda_{\prim} \cap
C_{r_0}\neq \varnothing$. 
In other words, 
$\XX_n(C_{r_0}) = \XX_n$. Also set \index{D@$D_r$ -- horizontal disk
  at height one}
\begin{equation}\label{eq: def Dro}
D_r \df \left\{\bx\in \R^n : \left\|\pi_{\R^d}(\bx) \right \| \leq r,
  \ x_n =1 \right\},
\end{equation}
and define  \index{S@$\sro$ -- cross-section }
\begin{align}\label{eq:main sets}
  \sr = \XX_n(D_{r})
  \  \ \text{ and } \ 
  \sr^{\sharp} = \XX_n^{\sharp}(D_{r}).
\end{align}
The set $\sro$ will be our cross-section.
\begin{remark}
  The reader will note that we have not specified $r_0$ explicitly,
  e.g. we did not specify an equality
  $m(C_{r_0})=2^n$. Two situations in which this additional
  flexibility will be useful, are when dealing with $\vre$-approximations for
  large $\vre$ (in which case we will require $r_0 \geq \vre$), and in
  the proof of  Proposition \ref{prop: weak stable 3}.  
  \end{remark}
  \begin{figure}[htbp]
\center{
\tdplotsetmaincoords{80}{30}
\begin{tikzpicture}[tdplot_main_coords,   scale=2]

\draw[black, thick, ->] (-1.3,0,0) -- (1.3,0,0);
\draw[black, thick, ->] (0,-1.5,0) -- (0,1.5,0);
\draw[black, thick, dashed] (0,0,0) -- (0,0,1);
\draw (1,1,0) node[anchor=west]{\!\tiny{$\bR^d$}};

\draw [thick, draw=black, fill=red, opacity=0.2] (0,0,1) circle (1);
\draw [thick, draw=black] (0,0,1) circle (1);

\filldraw [thin, draw=black, fill=gray, opacity=0.8] (0,0,1) circle (0.6);
\filldraw [thin, draw=black, fill=gray, opacity=0.8] (0,0,0) circle (0.6);

\draw[black, thick, ->] (0,0,1) -- (0,0,2) node[anchor=west]{\!\tiny{$\bR\mb{e}_n$}};
\tdplotsetrotatedcoords{0}{0}{0}
\begin{scope}[tdplot_rotated_coords]
\foreach \y in {0,0.02,...,1}
{\tdplotdrawarc[blue, very thin, opacity = 0.3]{(0,0,\y)}{0.6}{0}{360}{}{}{}}
\end{scope}

\filldraw [black] (0,0,1) circle (0.5pt);
\draw (0,0,1.05) node[anchor = west]{\!\tiny{$\mb{e}_n$}};
\filldraw [black] (0.6,0,1) circle (0.5pt);


\draw[black,  thin, ->] (1.5,0,0.7)--(0.71,0,0.5);
\draw (1.5,0,0.7) node[anchor = west]{\!\tiny{Upper half of $C_r$}};

\draw[black, thin] (0,0,1)--(0.6,0,1);
\draw[black, thin, ->] (0.6, 0, 1.5)--(0.3, 0 , 1.02);
\draw (0.65, 0, 1.5) node[anchor = south]{\!\tiny{$r$}};
\draw[black,  thin, ->] (-1.6,0,1.2)--(-1.02,-0.2,1);
\draw (-1.6,0,1.2) node[anchor = south]{\!\tiny{$D_{r_0}$}};

\end{tikzpicture}
}


\end{figure}

\begin{lemma}\label{lem:Sstar}
Let $\mu$ be any $\{a_t\}$-invariant probability measure on $\XX_n$. Then
$\sro$ is a $\mu$-cross-section 
for $(\XX_n,\mu, \{a_t\})$. Furthermore, the cross-section measure
satisfies
\begin{equation}\label{eq: cross section measure small}
  \mu_{\sro}\left((\XX_n(D_{r_0}, 2)
    \right) =
0.
\end{equation}
\end{lemma}
\begin{proof}
We verify Definition~\ref{def:cs}. We take $X_0$ to be the set of
lattices which intersect 
the horizontal space and the vertical axis $\spa\set{\mb{e}_n}$ only at $\set{0}$. 
It is easy to see that this is an $\{a_t\}$-invariant 
$G_\del$-set and therefore  $(X_0,\cB_{X_0})$ is a  
standard Borel space. In order to show that $\cS$ is  a
$\mu$-cross-section we need to check three things:
that $\mu(\XX_n\smallsetminus X_0) =0$, that for any $\Lam\in X_0$ the
set of visit times $\cY_{\Lam}$ is discrete and unbounded  
from below and above, and finally that the return time function
$\tau_\cS:\cS\cap X_0\to \bR_+$ is measurable.  

Let $\Lam\in \XX_n$. We say that $\Lam$ is {\em divergent in positive
  (negative) time} if for any compact $K \subset \XX_n$ there is
$t_0\in \bR$ such that for all $t> t_0$ (respectively, for all $t<t_0$) we
have $a_t \Lam \notin K$. Recall that by Mahler's compactness criterion,
a closed subset $K \subset \XX_n$ is compact if and only if there is
$\vre>0$ such that for all $\Lam \in K$, any nonzero $v \in \Lam$
satisfies $\| v\| \ge \vre$. In particular, if $\Lam\in
\XX_n\smallsetminus X_0$ then either it is divergent 
in positive time or in negative time. By Poincar\'e
recurrence, $\XX_n\smallsetminus X_0$ is a $\mu$-null set. 

Let $\Lam\in X_0$. We verify
\eqref{eq: requirements visit times}, that is, we show
that $\set{t \in \R : a_t\Lam\in \sro}$ is discrete and unbounded from
below and from above. 

Note that $a_t \Lam \in \sro$ if and only if $\Lam_{\prim}$ contains a
vector in $a_{-t}(D_{r_0})$. 
Discreteness of the 
set of visit times readily follows from this, and the fact that $\Lam_{\prim}$ is a
discrete set in $\bR^n$. Now suppose by contradiction that there is
$T>0$ such that for all $s \geq T$,
$a_s\Lam\notin \sro$.  
The set $F \df a_T\Lam\cap C_{r_0}$ is finite and since $\Lam\in X_0$, all vectors in $F$ have non-zero
horizontal component.  
It follows that for all large enough $t>T$, 
$a_{t-T}(F) \cap C_{r_0}=
\varnothing$. By Minkowski's convex body theorem  and the choice of $r_0$,
the lattice  
$a_t\Lam$ contains a primitive vector $v = (v_1, \ldots, v_n)$ in the
cylinder $C_{r_0}$. Since $\Lam\in X_0$, we have $v_n \neq 0$, and we can assume without loss of
generality that $v_n>0.$ Let $v = a_tv_0$ for $v_0 \in \Lam_{\prim}$. By
\eqref{eq: def at} there is a unique $s$ such that $a_s v_0 \in
D_{r_0}$, and since the vertical component of $v$ is at most 1, we
have $s \leq t$. This means that $a_s \Lam \in \sro$ and by choice of
$T$ we have $s \leq T$, so that the vertical component of $a_Tv_0$
is at most 1. On the other hand
$$
\| \pi_{\R^d}(a_Tv_0)\| = \|\pi_{\R^d}(a_{T-t}v)\| =
e^{t-T}\|\pi_{\R^d}(v)\| \leq r_0,
$$
so that $a_Tv_0 \in C_{r_0}$. This shows that $a_Tv_0 \in F$ and
hence $v \in a_{t-T}F$, a contradiction. 
The argument showing
unboundedness from below is similar. 

To complete the verification of Definition~\ref{def:cs}  
we need to show that the return time function is Borel measurable, or
equivalently, that the sub-level sets 
$\sroleps$ defined in \eqref{eq: sublevel sets} are Borel. The
thickened set 
\begin{equation}\label{eq: def Weps}
  W_\vre \df D_{r_0}^{(-\vre,0)}=
  \bigcup_{t\in(0,\vre)} 
  a_{-t} (D_{r_0})
  \end{equation}
is a Borel subset of $\R^n$, and therefore 
\begin{equation}\label{eq:1402}
\sroleps = \set{\Lam\in \sro : \Lam_{\prim}\cap W_\vre\ne\varnothing} =
\XX_n(D_{r_0}) \cap \XX_n(W_\vre) 
\end{equation} 
 is Borel as well.

To finish the proof we show \eqref{eq: cross section measure small} holds. This follows since  $\XX_n(D_{r_0}, 2)^\bR\subset \XX_n\smallsetminus X_0$
(where we use the notation
 in \eqref{eq: thickening}) is $\mu$-null  and so by
 Theorem~\ref{cor:csmeasure}\eqref{1204-4}, $\XX_n(D_{r_0},2)$ 
 is $\mu_{\sro}$-null.
 \end{proof}

\subsection{Parameterizing $\sro$}\label{subsec: parameterizing}
One of the advantages of the cross-section $\sro$ is that it 
has a nice description in terms  
of orbits and groups. Let
\index{U@$U$ -- the expanding horospherical group}
\begin{align}\label{eq: def subgroups}
&H=\set{\smallmat{A&0\\ w^{\mathbf{t}}&1} : A\in \SL_d(\bR), \, w\in \bR^d}
                             , \ \text{ and } \\ 
&U=\set{u(v) : v\in \bR^d}, \ \ \text{ where } u(v) \text{ is as in
                                                    \eqref{eq: def
                                                    u(v)} } 
                                                           ;
\end{align}
that is,  $H $ is the group defined in \eqref{eq: non-expanding
group}, the group generated by $H$ and $\{a_t\}$ is the group
$H^{\leq}$ appearing in Proposition
\ref{prop: horospherical and generic}, the orbit $H\Z^n$ is the space 
$\crly{E}_n= \XX_n(\mathbf{e}_n)$ of lattices which contain $\mathbf{e}_n$ as a
primitive vector (see \eqref{eq: all functionals}), and $U$ is the {\em expanding
  horospherical group of $a_t$ in positive time}.

We let \index{B@$\bar B_r$ -- closed ball in $\R^d$} $\bar{B}_r
\subset \R^d$
denote the closed ball centered at $\mathbf{0} \in \R^d$, with
respect to our chosen norm (note that the norm is suppressed from the
notation).
Consider the map \index{P@$\varphi$ -- the map parameterizing $\sro$}
$$\vphi : \crly{E}_n \times \bar{B}_{r_0}\to \sro,\quad
\varphi(\Lam, v)\df u(v)\Lam.$$ 
Note that the map $v\mapsto
u(v)\mb{e}_n$ is a bijection between  
$\bar{B}_{r_0}$ and $D_{r_0}$. It follows that $\vphi$ is onto
$\sro$, and 
\begin{equation}\label{eq: smaller r}
  \text{ for any } r \in (0, r_0), \ \ 
\varphi(\crly{E}_n \times \bar B_r) = \cS_r.
  \end{equation}
Furthermore, for $\Lambda \in \sro$,  
\begin{equation}\label{eq:fiber size}
\# \,\vphi^{-1}(\Lambda) = \# \left( \Lambda_{\on{prim}}\cap D_{r_0}\right).
\end{equation}
Indeed, for any $v\in \Lam_{\prim} \cap D_{r_0}$,
$$\left(u(\pi_{\bR^d}(v))^{-1}\Lam, \pi_{\R^d}(v) \right) \in
\varphi^{-1}(\Lam),$$ 
and this assignment is easily
seen to be a bijection.
Let \index{P@$\psi$ -- the inverse of $\varphi$} \index{V@$v_\Lam$ --
  the horizontal component of $v(\Lam)$}
\begin{align}
\label{eq: psi coordinates}
  &\psi : \srosharp \to \crly{E}_n \times \bar{B}_{r_0}, \\
 \nonumber  &\psi(\Lam) \df  \left(
              u(
              v_\Lam)^{-1}\Lam, v_\Lam \right),
  \end{align}
  where
  \begin{equation}\label{eq: vLam}
    v_\Lam \df \pi_{\R^d}(v(\Lam)) \ \text{ and } \{v(\Lam)\} = \Lam \cap
    D_{r_0}.
    \end{equation}
\begin{lemma}\label{lem:sro-is-lcsc} The set
$\sro$ is closed in $\XX_n$, $\psi$ as in
\eqref{eq: psi coordinates} 
is the inverse of $\vphi|_{
  \varphi^{-1}(\srosharp)}$, and is a homeomorphism  
between $\srosharp$ and $
\varphi^{-1}(\srosharp)$. 
\end{lemma}
\begin{proof}
  The assertions that $\sro$ is closed, and that $\psi$ is
  continuous, follow from 
Lemma~\ref{lem:continuity of vector}.
It is clear  that $\vphi$ is continuous. To see that $\psi$ and
$\varphi|_{\varphi^{-1}(\srosharp)}$ are mutual inverses, we see easily
that $\varphi \circ \psi = \mathrm{Id}_{\srosharp}$. By
\eqref{eq:fiber size}, $\varphi$ is injective on
$\varphi^{-1}(\srosharp)$, and thus we also have 
$\psi
\circ \varphi|_{\varphi^{-1}(\srosharp)} =
\mathrm{Id}_{\varphi^{-1}(\srosharp)}$.
\end{proof}

In the following subsections we will use the map $\varphi$ to describe
the cross-section measure 
$\mu_{\sro}$, corresponding to certain invariant measures $\mu$.
\subsection{The cross-section measure, real homogeneous space, Case 
  I}\label{subsec: 8.3}
The goal of this subsection is the following result:
\begin{theorem} \label{thm: reasonable case I}
The cross-section $\sro$ is $m_{\XX_n}$-reasonable.  
\end{theorem}

For the proof we will need some preparations. The first is an
explicit description  of the cross-section measure \index{M@$m_{\R^d}$ --
  Lebesgue measure on $\R^d$} $\mu_{\sro}$. 
Let $m_{\R^d}$ be the Lebesgue measure on $\R^d$.

\begin{proposition}\label{prop:description of measures1}
  In Case I, 
  $$\mu_{\sro} =
  \frac{d}{\zeta(n)}\, \varphi_*\left(
    m_{\crly{E}_n} \times m_{\R^d}|_{\bar B_{r_0}} \right)$$
  (where $\zeta(n) = \sum_{k\in \N} k^{-n}$). In
particular, $\mu_{\sro}$ is finite and $\on{supp}(\mu_{\sro}) = \sro$.
\end{proposition}
The proof is a straightforward but lengthy computation which is 
postponed to \S \ref{subsec: proof of prop}. 

Here is another result used in the proof of 
 Theorem \ref{thm: reasonable case I}. We state it in a general form which will
 be useful in the sequel. 

 \begin{lemma}\label{lem:periodic-measures-nullity}
Let $L \subset \SL_n(\bR)$ be a closed subgroup, with left Haar measure
$m_L$. Let $\Lam_0 \in \XX_n$ such that
$L\Lam_0$ is a closed orbit supporting a finite $L$-invariant measure
$m_{L\Lam_0}$. Let $W\subset \bR^n$ such that for any $v \in W \cap
L\Lam_0$ there 
is $\vre>0$ such that
\begin{equation}\label{eq: weak transversality}
  m_L\left(\left\{ \ell \in B^L_\vre: \ell v
      \in W \right\} \right)=0
\end{equation}
(where $B^L_\vre$ denotes the $\vre$-ball around the identity element
of $L$, with respect to some metric inducing the topology).
Then $m_{L\Lam_0} \left(\XX_n(W) \right)=0.$

\end{lemma}
\begin{proof}
By covering $W$ with countably many bounded sets we may assume
that $W$ is bounded.
Since $W$ is bounded, for any $\Lam \in L\Lam_0$, the cardinality of 
$\Lam_{\on{prim}}\cap W$ is finite, and bounded for $\Lam$ in a 
compact subset of $L\Lam_0$. 
Using the hypothesis we
deduce that for each 
$\Lam\in L\Lam_0\cap \XX_n(W)$ there exist $\vre>0$ such that
\eqref{eq: weak transversality} holds for any $v\in
\Lam_{\on{prim}}\cap W$. Since $m_{L\Lam_0}$ is the restriction of
$m_L$ to a fundamental domain for the action of the stabilizer $L_{\Lam_0}$, it follows 
that $m_{L\Lam_0}\left( \XX_n(W) \cap B^L_\vre \Lam \right)=0$. 
We can cover $L\Lam_0$ by countably many  sets $\{B^L_\vre \Lam_i\}_{i
  \in \N}$, and therefore
$m_{L\Lam_0}(\XX_n(W))=0.$
\end{proof}

For positive $r$ and $\vre$, define
$$
D_r^\circ \df \left\{\bx\in \R^n : \left\|\pi_{\R^d}(\bx) \right \| < r,
  \ x_n =1 \right\}
$$
(compare with the set $D_{r}$ defined in \eqref{eq: def
  Dro}), and
$$F_{\vre,r} \df 
\left(D_{r}\smallsetminus
  D_{r}^\circ \right)^{[-\vre,0]}.$$
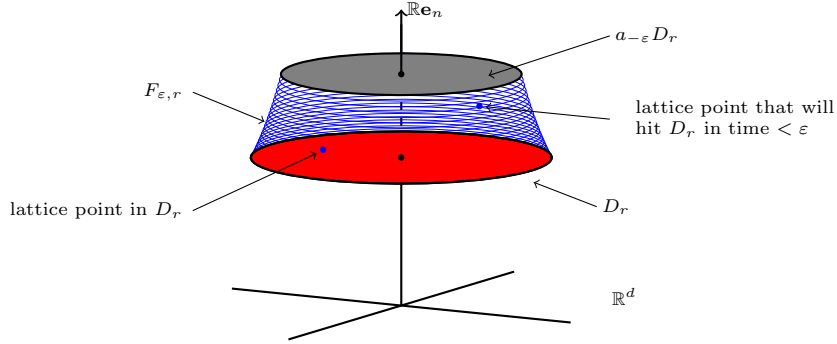
\begin{figure}[htbp]
\center{
\tdplotsetmaincoords{80}{30}
\begin{tikzpicture}[tdplot_main_coords,   scale=2]

\draw[black, thick] (-1.3,0,0) -- (1.3,0,0);
\draw[black, thick] (0,-1.5,0) -- (0,1.5,0);
\draw (1,1,0) node[anchor = west]{\!\tiny{$\bR^d$}};
\draw[black, thick] (0,0,0) -- (0,0,1);
\draw[black, thick, dashed] (0,0,1) -- (0,0,1/0.64);
\draw[black, thick, ->] (0,0,1/0.64) -- (0,0,2) node[anchor=west]{\!\tiny{$\bR \mb{e}_n$}};
\tdplotsetrotatedcoords{0}{0}{0}
\begin{scope}[tdplot_rotated_coords]
 \foreach \y in {0.8,0.81,...,1}
 {\tdplotdrawarc[blue, very thin, opacity = 0.4]{(0,0,1/\y^2)}{\y}{0}{360}{}{}{}}
\end{scope}
\draw [thick, draw=black, fill=red, opacity=0.2] (0,0,1) circle (1);
\draw [thick, draw=black] (0,0,1) circle (1);
\filldraw [thick, draw=black, fill=gray, opacity=0.8] (0,0,1/0.64) circle (0.8);
\draw[black, thick] (0,0,1/0.64) -- (0,0,1.9);
\filldraw [black] (0,0,1) circle (0.5pt);
\filldraw [black] (0,0,1/0.64) circle (0.5pt);

\draw[black,  thin, ->] (1.5,0,0.8)--(0.67,0.6,0.8);
\draw (1.5,0,0.8) node[anchor = west]{\!\tiny{$D_{r}$}};

\draw[black,  thin, ->] (1.6,0,1/0.64 +0.4)--(0.67,0,1/0.64 + 0.05);
\draw (1.6,0,1/0.64 +0.4) node[anchor = west]{\!\tiny{$a_{-\vre}D_{r}$}};

\draw[black,  thin, ->] (-1.6,0,1.3)--(-1.05,0,1.15);
\draw (-1.6,0,1.3) node[anchor = east]{\!\tiny{$F_{\vre,r}$}};

\filldraw [blue] (-0.6,0,1) circle (0.5pt);
\filldraw [blue] (0.6,0,1.4) circle (0.5pt);
\draw[black,  thin, ->] (-1.6,0,0.5)--(-0.65,0,0.97);
\draw (-1.6,0,0.5) node[anchor = east]{\!\tiny{lattice point in $D_r$}};
\draw[black,  thin, ->] (1.6,0,1.4)--(0.65,0,1.4);
\draw (1.6,0,1.4) node[anchor = west]
{\!\tiny{
\begin{tabular}{ll} 
lattice point that will \\ hit  $D_r$ in  time $< \vre$
\end{tabular}}};

\end{tikzpicture}
}

\caption{The set $D_r^{(-\vre,0)}$ bounded by the surface $F_{\vre, r}$ and the disks
  $D_r$ and $a_{\vre}(D_r)$. Lattice points in 
  $D_r^{(0, \vre)}$ correspond to 
  visits to $\cS_{r,<\vre}$.  }
\label{fig: circus tent}
\end{figure}

See Figure \ref{fig: circus tent}. 
We will need the following: 
\begin{lemma}\label{lem: boundary estimate}
For any $r\ge r_0, \, \vre>0$ we have that 
\begin{align}\label{eq: 13211}
\partial_{\sr}( \srleps)&\subset \XX_n(D_{r},2)
                            \cup\pa{\XX_n(D_{r})\cap
                            \XX_n(a_{-\vre} (D_{r}))} \cup \XX_n(F_{\vre,r}),
\end{align}
and the set $\srleps$ is $\mu_{\sr}$-JM. 
\end{lemma}
\begin{proof}
Let $W_\vre$ be as in \eqref{eq: def Weps}  and let $W_\vre^\circ$,
$\overline{W}_\vre$ be its interior and closure in $\bR^n$. Note that 
\begin{equation}\label{eq:11411}
\overline{W}_\vre = W_\vre\cup D_{r} \cup a_{-\vre}(D_{r} )\cup F_{\vre,r}.
\end{equation}
Recall 
from \eqref{eq:1402}
that $\srleps$ consists of the lattices containing a primitive point
in $D_{r}$ and another primitive point in $W_\vre$.  
In particular, $\srleps$ is contained in
$$E_1 \df \XX_n(D_{r})\cap \XX_n(\overline{W}_\vre,2),$$
which is closed by Lemma
\ref{lem:continuity of vector}. By Lemma~\ref{lem:open sets},
$\XX_n(W_\vre^\circ)$ is open and in conjunction with 
Lemma~\ref{lem:continuity of vector} we deduce that the set 
$$E_2 \df \XX_n^{\#}(D_{r})\cap \XX_n(D_{r}^\circ) \cap
\XX_n(W_\vre^\circ) 
$$
is open in
$\XX_n(D_{r})$. It follows that  
$$\partial_{\sr}(\srleps) \subset E_1\smallsetminus E_2.$$
We now show that 
$E_1\smallsetminus E_2$ is contained in the RHS of \eqref{eq:
  13211}. 
Suppose $\Lam \in E_1$ and hence $\Lam_{\on{prim}}$ contains two distinct vectors $v_1,v_2$
such that $v_1\in D_{r}$ and $v_2\in \overline{W}_\vre$. If $\Lam$
does not belong to the  RHS of \eqref{eq: 13211}, then
$v_1$ must be the unique  
primitive vector in $D_{r}$, it cannot lie in $F_{\vre,r}$ and so it
must lie in $D_{r}^\circ$ and $v_2$ cannot lie in $F_{\vre,r}$ or 
$a_{-\vre}(D_{r})$ and so it must lie in $W_\vre^\circ$. That is, $\Lam\in
E_2$. This concludes the proof of \eqref{eq: 13211}. 

In order to show that $\srleps $ is $\mu_{\sr}$-JM, it suffices to show that
the sets
$$\XX_n(D_{r},2),\; \XX_n(F_{\vre,r}),\; \XX_n(D_{r})\cap
\XX_n(a_{-\vre}(D_{r}))$$
are all $\mu_{\sr}$-null.
By Theorem~\ref{cor:csmeasure} it is enough to show that the sets 
$$\XX_n(D_{r},2)^\bR,\; \XX_n(F_{\vre,r})^\bR,\; \pa{\XX_n(D_{r})\cap
  \XX_n(a_{-\vre}(D_{r}))}^\bR$$ are all $m_{\XX_n}$-null.

The set $\XX_n(D_{r},2)^\bR$ is $m_{\XX_n}$-null because of 
\eqref{eq: cross section measure small}. Next, if we set
\begin{equation}\label{eq: def M r}
M_r \df  \set{\mb{x}\in \bR^n: x_n \cdot 
  \norm{\pi_{\R^d}(\mb{x})}^d = r^d} = 
 \bigcup_{t\in \bR} a_t\left(D_{r}\smallsetminus
 D_{r}^\circ \right),
\end{equation}
then 
\begin{equation}
  \label{eq: no waste}
  \XX_n (D_r \sm D_r^\circ)^{\bR} = \XX_n(F_{\vre,r})^\bR =
  \XX_n(M_{r}).
  \end{equation}
Since $\SL_n(\R)$ contains elements which expand the vertical
component $x_n$ without affecting $ \norm{\pi_{\R^d}(\mb{x})}$,
condition \eqref{eq: weak transversality} is satisfied for $W=M_r$, so
applying Lemma~\ref{lem:periodic-measures-nullity} we have that
$\XX_n(M_{r})$ is $m_{\XX_n}$-null.

It remains  to show that 
$$E_\vre \df \pa{\XX_n(D_{r})\cap \XX_n(a_{-\vre}(D_{r}))}^\bR =
\bigcup_{t\in\bR} (\XX_n(a_t(D_{r}))\cap \XX_n(a_{t+\vre}(D_{r}))$$ 
is $m_{\XX_n}$-null. If $\Lam \in \XX_n(a_t(D_{r}))\cap
\XX_n(a_{t+\vre}(D_{r}))$ for some $t$,   
and $v,w$ are primitive vectors in $\Lam$ such that $v\in a_t(D_{r})$
and $w\in a_{t+\vre}(D_{r})$, then $v,w$ are two linearly 
independent vectors in $\Lam$ and the ratio between their vertical
components satisfies $\frac{w_n}{v_n} = e^\vre$. 

We use a result of Siegel \cite{Siegel_mean_value_theorem}, according
to  which for a null set in $\Om\subset \bR^n\times \bR^n$, the function
$$\widehat{\chi_\Om}(\Lam) = \sum_{\substack{v, w\in \Lam \\ \text{
    linearlty  independent}}} \chi_\Om(v,w)$$
has integral zero with respect to $m_{\XX_n}$. We apply this with 
$$\Om \df \set{(v,w)\in \bR^n\times \bR^n : \frac{w_n}{v_n} =e^\vre}.$$
Since  $\widehat{\chi_\Om}$ bounds the 
characteristic function of $E_\vre$ from above, we obtain that $E_\vre$
is $m_{\XX_n}$-null. 
\end{proof}

Let
\begin{equation}\label{eq: def uro}
  \uro \df \XX_n^{\#}(D_{r_0}) \cap \XX_n(D_{r_0}^\circ).
  \end{equation}
\begin{lemma}\label{lem: is open1}
The set $\uro$ is open in $\sro$, the set
$(\on{cl}_{X_n}(\sro)\smallsetminus \uro)^{(0,1)}$ is
$m_{\XX_n}$-null, and the map  $(t,\Lam)\mapsto a_t\Lam$ from $(0,1)\times \uro$ to $\XX_n$ is open. 
\end{lemma}
\begin{proof}
Item \eqref{eq:12112} of Lemma~\ref{lem:continuity of
  vector} shows that $\uro$ is open in $\sro$, and item \eqref{eq:12111} shows
that $\sro = \XX_n(D_{r_0})$ is  
closed in $\XX_n$. Therefore
\begin{equation} \label{eq: therefore}
  \begin{split} \on{cl}_{\XX_n}(\sro)\smallsetminus \uro = & \XX_n(D_{r_0})\smallsetminus
(\XX_n^{\#}(D_{r_0})\cap \XX_n(D_{r_0}^\circ)) \\ \subset &  
\XX_n(D_{r_0},2)\cup \XX_n(D_{r_0}\smallsetminus
D_{r_0}^\circ).\end{split} \end{equation}
In the proof of Lemma~\ref{lem: boundary estimate} we showed that
$\XX_n(D_{r_0},2)^\bR$ as well as 
$\XX_n(D_{r_0}\smallsetminus D_{r_0}^\circ)^\bR$ are
$m_{\XX_n}$-null. This proves the second assertion.

For the third assertion, let
$U$ and $H$ be as in \eqref{eq: def subgroups}, let $Q \df
  \{a_t\} \ltimes U$, and let $\mathfrak{h, q, g}$ denote respectively
  the Lie algebras of $H, \, Q$ and $\SL_n(\R)$.  The product map $(t, u)
  \mapsto a_tu$ is a homeomorphism $\R \times U \to Q$, and since 
 $\mathfrak{g} = \mathfrak{q} \oplus \mathfrak{h}$, the product map 
$\bR\times U\times H\to G$ is open. This implies that the  map
$$\bR\times U\times H\bZ^n\to \XX_n, \ \ (t, u, h\bZ^n) \mapsto a_t u h\bZ^n$$
is open, and the map 
in the statement of the Lemma is its restriction to an open set because
$\XX_n(D_{r_0}^\circ) = B_{r_0}^U\cdot H\bZ^n$ (see \S \ref{subsec: parameterizing}).  
\end{proof}

\begin{proof}[Proof of Theorem \ref{thm: reasonable case I}]
Property \eqref{item: probability} of \S \ref{sec: JM sections} is
immediate, property \eqref{item: B} 
follows from Proposition 
\ref{prop:description of measures1}, and property \eqref{item: C} follows from Lemma
\ref{lem:sro-is-lcsc}. 
Item \eqref{item: reasonable 1} of Definition 
\ref{def:mureasonable} follows from Lemma \ref{lem: boundary estimate}, and item
\eqref{item: reasonable 2} follows from
Lemma \ref{lem: is open1}.
\end{proof}

\subsection{The cross-section measure, real homogeneous space, Case 
  II}\label{subsec: 8.4} 
In this subsection, the notation is as in \S \ref{subsec:
      case II}, and we write $\mu = m_{\vec{\alpha}}$ (see \eqref{eq: def m
      alpha}). We denote by $\mu_{\sro}$ the corresponding measure on
    $\sro$, defined via Theorem \ref{cor:csmeasure}. The goal of this
    subsection is the following result:
    \begin{theorem}
\label{thm: reasonable case II} 
With this choice of $\mu$, the cross-section $\sro$ is $\mu$-reasonable. 
\end{theorem}
\begin{remark}\label{remark: failure for d=1}
Recall that in Case II we always assume $d\ge2$. The reason is that for $d=1$, for some choices of $r_0$, 
$\sro$ may fail to be $\mu$-reasonable. This is because the Jordan measurability of some sets may fail. More 
specifically, the conclusion of Lemma~\ref{lem: using convexity} may
fail for some values of $r$.  
\end{remark}
We will need a detailed description of $\mu_{\sro}$. Let $x_{\vec{\al}}^*$ be
as in \eqref{eq: def lattice with compact orbit},  let $\bar h_\alpha$
be as in \eqref{eq: def Bar B},  let $\bar A_{\vec{\al}}, \,
y_{\vec{\al}}$ be as in \eqref{eq: def bar A}, so that
$$\supp \, m_{\vec{\alpha}} = \bar h_\alpha A x_{\vec{\al}}^* = \bar
A_{\vec{\al}} y_{\vec{\al}}.$$ 
Also let \index{A@$\bar A_{\vec{\al}}^{(1)}$ -- subgroup of $\bar
  A_{\vec{\al}}$ preserving $\mathbf{e}_n$} 
\begin{equation}\label{eq: def bar A 1}
\bar A_{\vec{\al}}^{(1)} \df \{ \bar a \in \bar A_{\vec{\al}} : \bar a \mb{e}_n = \mb{e}_n\},
\end{equation}
and let $m_{\bar A_{\vec{\al}}^{(1)}}$ denote the Haar measure on
$\bar A_{\vec{\al}}^{(1)}$.  Note that both 
$\bar h_\alpha$ and $\bar A_{\vec{\al}}^{(1)}$ act on $\R^n$ 
without changing the vertical component of any vector. 

\begin{proposition}\label{prop:description of measures}
 The following hold: 
 \begin{enumerate}[(a)]
 \item \label{item: 1.1}
The support of $\mu_{\sro}$ is the compact set $\bar A_{\vec{\al}} y_{\vec{\al}} \cap
\sro$;
\item \label{item: 1.2}
  The support of $\mu_{\sro}$ is equal to 
a finite union  $\bigcup_{i=1}^k \cO_i\bar{\Lam}_i$ where each $\cO_i\subset \bar A_{\vec{\al}}^{(1)}$ is homeomorphic to a closed ball,
 and $\bar{\Lam}_i\in \bar{A}_{\vec{\al}}y_{\vec{\al}}$;
\item \label{item: 1.25}
  The restriction of $\mu_{\sro}$ to each of the subsets $\cO_i\bar{\Lam}_i$ 
  in \eqref{item: 1.2} is the pushforward of the restriction of the Haar measure $m_{\bar A_{\vec{\al}}^{(1)}}|_{\cO_i}$
  under the orbit map $\bar{a}\mapsto \bar{a}y_{\vec{\al}}$.
\item \label{item: 1.3} 
$\bar A_{\vec{\al}} y_{\vec{\al}} \cap 
\sro \subset \srogeps$ for some $\vre>0$;
\item\label{item: 1.4}
$\mu_{\sro}$ is finite. 
\end{enumerate}

\end{proposition}

\begin{proof}
  Item \eqref{item: 1.1} follows from
Theorem \ref{cor:csmeasure}\eqref{1204-4}. 
For item \eqref{item: 1.2}, we have from Proposition \ref{prop:
  compact A orbit} that the orbit  $Ax_{\vec{\al}}^*$ is of type
$(\pmb{\sigma}, \bK)$. 
Define $N:\bR^n\to \bR$ by 
\begin{equation}\label{eq: def N}
  N(\bx) \df \prod_{i=1}^n s_i, \ \ \text{ where } (s_1, \ldots, s_n)^\tb{t}
  = \bar h_{\vec{\al}}^{-1}\bx.
  \end{equation}
In other words, the restriction of $N$ to $y_{\vec{\al}}$ is obtained
from the norm $N_{\bK/\Q}$ by the geometric embedding and a 
change of variables. Then $N$ is $\bar A_{\vec{\al}}$-invariant and, since the $N_{\bK/\bQ}$ takes a 
discrete set of  values on $x_{\vec{\al}}^*$, there is a
sequence $\al_i\to\infty$  
such that the vectors comprising $y_{\vec{\al}}$
all lie in the countable union of hypersurfaces 
$\bigcup_{i=1}^\infty
N^{-1}(\set{\al_i})$, and each $N^{-1}(\set{\al_i})$ is a finite union
of $\bar A_{\vec{\al}}$-orbits. Only finitely many of these
hypersurfaces intersect the bounded set $D_{r_0}$ and therefore there
is a finite set $\{v_1, \ldots, v_k\}\subset y_{\vec{\al}}$ such that 
if we denote 
$$\bar{\cO}_i \df \set{\bar{a}\in \bar{A}_{\vec{\al}}:\bar{a} v_i\in D_{r_0}}\quad
\text { and } \quad
W_i \df  \bar{\cO}_i y_{\vec{\al}},$$
then
$$\bar A_{\vec{\al}} y_{\vec{\al}} \cap \sro = \bigcup_{i=1}^k W_i.$$
Let
$
L_1 \df \{\bx \in \R^n: x_n = 1\},$
so that $\bar A_{\vec{\al}}^{(1)}$ is the subgroup of $\bar A_{\vec{\al}}$ leaving $L_1$
invariant, and $D_r \subset L_1$. 
Thus we may write $\bar{\cO}_i = \cO_i \bar{a}_i $ for some $\bar{a}_i\in \bar{A}_{\vec{\al}}$ and with $\cO_i\subset \bar A_{\vec{\al}}^{(1)}$.
Setting $\bar{\Lam}_i = \bar{a}_iy_{\vec{\al}}$ we obtain \eqref{item: 1.2}. 

We show \eqref{item: 1.25}. Note that the orbit map $\bar{a}\mapsto \bar{a}\bar{\Lam}_i$ is injective on each 
$\cO_i$ since otherwise $y_{\vec{\al}}$ would contain two distinct vectors with the same vertical component, contradicting 
the boundedness of the orbit $\set{a_ty_{\vec{\al}}:t\in \bR}$. It follows that $\cO_i$ is contained in a fundamental domain
for the orbit $\bar{A}_{\vec{\al}}y_{\vec{\al}}$. Since 
$m_{\vec{\alpha}}$ can be identified with the restriction of the Haar measure of $\bar A_{\vec{\al}}$
to a fundamental domain via the orbit map, and since 
$$
\bar A_{\vec{\al}} = \bar A_{\vec{\al}}^{(1)} \times \{a_t\},
$$
item \eqref{item: 1.25} follows from \eqref{eq:1019}.

To prove item \eqref{item: 1.3}, suppose by contradiction that for any
$j \in \N$ there exists  $\Lam_j \in \bar  
A_{\vec{\alpha}} y_{\vec{\al}} \cap \mathcal{S}_{r_0, < 1/j}$. Then by
compactness we can take a convergent sequence to conclude that there
is 
$\bar{\Lam}' \in \bar A_{\vec{\al}}y_{\vec{\al}} \cap \crly{X}_n(D_{r_0},2)$. Hence $\bar\Lam'$
contains a nonzero vector with zero vertical component. The map 
$\bar h_{\vec{\al}}$ preserves the horizontal space $\R^d$, and
thus $ \bar h_{\vec{\al}} ^{-1} \bar {\Lam}' \in Ax_{\vec{\al}}^*$ also contains a nonzero
vector with zero vertical component, contradicting Proposition
\ref{prop: compact A orbit}.
Item \eqref{item: 1.4} now follows from
Theorem~\ref{cor:csmeasure}\eqref{1204-2}. 
\end{proof}

\begin{lemma} \label{lem: using convexity}
For 
any norm $\|\cdot \|$, any $r>0$, and any $\vec{\alpha}$ as in \eqref{eq: all 
  about alpha}, the set $M_r$ defined in  \eqref{eq: def M
  r} satisfies
$  m_{\vec{\alpha}} \left( \XX_n(
  M_r )\right)=0. 
$

  \end{lemma}
 
%
 %
%

 \begin{proof}
We apply Lemma \ref{lem:periodic-measures-nullity}, with 
$L = \bar A_{\vec{\al}}$ and $W = M_r$. We need to check \eqref{eq: weak
      transversality}. Let $m_{\bar A_{\vec{\al}}}, \, m_{\bar A_{\vec{\al}}^{(1)}}$ denote
    respectively the Haar measure on $\bar A_{\vec{\al}}$ and $\bar A_{\vec{\al}}^{(1)}$. 
Assume by way of contradiction  that there exists $v\in \bR^n\cap
y_{\vec{\al}} \cap M_r$ such that
$$
m_{\bar A_{\vec{\al}}} \left(\left\{ \bar a \in \bar A_{\vec{\al}}: \bar a v  \in M_r\right \} \right)>0.
$$
Since $M_r$ is $a_t$-invariant, and the action of $\bar A_{\vec{\al}}^{(1)}$ commutes
with the projection $\pi_{\bR^d}$ this implies that for $\bar v \df 
\pi_{\R^d} (v)$ we have 
$$
m_{\bar A_{\vec{\al}}^{(1)}} \left(B\right)>0, \ \ \ \text{ where } B \df \left\{ \bar
  a \in \bar A_{\vec{\al}}^{(1)}: \| \bar a \bar v \| = 
    \|\bar v \| \right \}.
$$
We can make a change of variables to replace  $\bar A_{\vec{\al}}^{(1)} $ with the group of 
diagonal matrices of positive diagonal entries and determinant
1 in $\SL_d(\R)$. We will denote this group by $A_2$. 
By Proposition \ref{prop: compact A orbit},
after this change of variables, the coordinates of $\bar v $ are nonzero.
Thus, by applying another linear change of variables on $\bR^d$ and changing
the norm,  we may assume that $\bar{v} =
\mathbf{1} $ is the vector all of whose coordinates are
equal to $1$.
Furthermore, we can 
replace $B$ with $b_0B$
for some $b_0$, in order to assume that 
the identity is a Lebesgue density point for $B$.
To get a contradiction we will show that 
there is $b\in B$ such that $\norm{b\mathbf{1}}>\norm{\mathbf{1}}$.

Define 
$$\Xi: A_2 \to \R^d, \ \ \ \ \Xi \left(\diag{x_1, \ldots, x_d}\right)
\df (x_1, \ldots, x_d)^{\mathrm{t}},$$ 
and 
$\nu \df \Xi_*
m_{A_2}.$
Then $\Xi$ is 
a diffeomorphism between $A_2$ and the hypersurface $\set{x\in 
  \bR^d: x_i >0,\; 
  \prod_1^d x_i = 1}$. Note
that $\mathbf{1} = \Xi(e)$ and the ray $\bR_+ \, \mathbf{1}$ is
transverse to $\Xi(A)$. Let $B_0$ be the image under $\Xi$ of the
intersection of $B$ with a  
bounded neighborhood of $e$. Then $\mathbf{1}$ is a density point for
$\nu$, i.e. 
\begin{equation}\label{eq: half space}
\lim_{r\to 0+} \frac{\nu(B_0 \cap B(\mathbf{1}, r))}{\nu(B(\mathbf{1},
  r))} =1.
\end{equation}
This implies that the interior of the convex hull of $B_0$ intersects
the ray $\R_+ \mathbf{1}$, since if this did not 
hold there would be a linear functional vanishing on $\mathbf{1}$ and
non-negative on $B_0$, and the left hand side of \eqref{eq: half space}
would be at most $\frac12$. 

It follows that there are distinct $\mb{x}_1,\dots, \mb{x}_{d+1}\in \Xi(B)$ 
and positive scalars $\be_1, \ldots, \be_{d+1}, c$, such that
$$\sum_{i=1}^{d+1}\be_i = 1 \ \ \text{ and } \ \sum_{i=1}^{d+1} \be_i
\mb{x}_i = c\mathbf{1}.$$   
Write $\mb{x}_i = (x_{i1},\dots, x_{id})$. By  the inequality of means
in the $j$th coordinate,
for each $1\le j\le d$ we have  
$$\prod_{i=1}^{d+1} x_{ij}^{\be_i} \leq \sum_{i=1}^{d+1} \be_i x_{ij}
= c,$$
with strict inequality for at least one $j$ since the $\mb{x}_i$ are
distinct. 
Multiplying these inequalities for $1\le j\le d$ and taking into
account that for each $i$, $\prod_{j=1}^d x_{ij} = 1$, we get  
$$1 = \prod_{j=1}^d \prod_{i = 1}^k x_{ij}^{\be_i} <c^n,$$
and therefore $c>1$.

Since $\mb{x}_i \in \Xi(B)$ we have $ \norm{\mb{x}_i} =
\norm{\mathbf{1}}$ for each $i$. By the triangle inequality
$$c\norm{\mathbf{1}} = \left\|\sum_1^k \be_i \mb{x}_i \right\|\le \sum \be_i
\norm{\mb{x}_i} =\sum \be_i \norm{\mathbf{1}} = \norm{\mathbf{1}},$$ 
a contradiction.

\end{proof}

\begin{proof}[Proof of Theorem \ref{thm: reasonable case II}]
Property \eqref{item: probability} of \S \ref{sec: JM sections} is
immediate, property \eqref{item: B}
follows from Proposition 
\ref{prop:description of measures}\eqref{item: 1.4}, and property
\eqref{item: C}  follows from Lemma 
\ref{lem:sro-is-lcsc}. 
Item \eqref{item: reasonable 1} of Definition 
\ref{def:mureasonable} follows from Proposition \ref{prop:description
  of measures}\eqref{item: 1.3}. It remains to prove 
\eqref{item: reasonable 2}. Define $\uro$ using \eqref{eq: def
  uro}. In light of Lemma \ref{lem: is open1}, we only need to show
that 
\begin{equation}\label{eq: only need to show}
 m_{\vec{\alpha}} \pa{(\on{cl}_{\XX_n}(\sro)\smallsetminus
    \uro)^{(0,1)}} =0.
  \end{equation}

For \eqref{eq: only need to show} it is enough to show that the two
sets on the RHS of \eqref{eq: therefore} are $\mu_{\sro}$-null, and
hence, by Theorem \ref{cor:csmeasure}\eqref{1204-4}, that 
$$
m_{\vec{\alpha}} \left(
\XX_n(D_{r_0},2)^\bR
\right ) = 0 \ \ \ \text{ and } \ m_{\vec{\alpha}}\left  (
\XX_n(D_{r_0}\smallsetminus D_{r_0}^\circ)^\bR
\right) =0.
$$
Since any lattice in $\XX_n(D_{r_0},2)^{\bR}$ contains a  nonzero
horizontal vector,  and $\bar h_{\vec{\al}}$ preserves the horizontal
subspace, Proposition \ref{prop: compact A orbit} implies
$\supp \, m_{\vec{\alpha}} \cap \XX_n(D_{r_0},2)^{\bR} = \varnothing.$ 
The second equality follows from \eqref{eq: no waste} and  Lemma
\ref{lem: using convexity}. 
\end{proof}
\subsection{Proof of  Proposition
\ref{prop:description of measures1}}\label{subsec: proof of prop}
We
will need the following: 
\begin{lemma}\label{lem:hic basic}
  Let $L$ be Lie group and let $L_1,L_2$ be closed subgroups such that
  $L_1\cap L_2 = \set{e}$ and $\dim L_1 + \dim L_2 = \dim L$. Then:
  \begin{enumerate}
  \item
    the
map  
$\al:L_1\times L_2\to L$ given by $\al(\ell_1,\ell_2) =
\ell_1\cdot\ell_2$ is a diffeomorphism onto an open subset $\cU\subset
L$.
\item If furthermore $L$ is unimodular, and $m_{L_1}^{\mathrm{left}}, m_{L_2}^{\mathrm{right}}$ denote left
and right Haar measures on $L_1,L_2$ respectively, then
$\al_*\left(m_{L_1}^{\mathrm{left}} \times m_{L_2}^{\mathrm{right}}\right)$ 
is proportional to the restriction to $\cU$ of a Haar measure on $L$.
\end{enumerate}

\end{lemma}

Lemma \ref{lem:hic basic} is standard, see e.g. \cite[Lemma
11.31]{EW}.

\begin{proof}[Proof of Proposition \ref{prop:description of measures1}]
  Let $G = \SL_n(\R)$, $\Gamma_G \df \SL_n(\Z)$, $\Gamma_H \df
  H(\Z)$, $\pi_G: G \to G/\Gamma_G$ and $\pi_H: H \to H/\Gamma_H$ the
  projections, and let $m_{G}$ and $m_H$ denote respectively the Haar
  measures 
  on $G$ and on $H$. Recall that a fundamental domain for $G/\Gamma$
  in $G$ is a Borel subset  $\Omega \subset G$ for which $\pi_G|_{\Omega}$ is a
  bijection, and that every Borel set on which $\pi_G$ is injective is
  contained in a fundamental domain. One can describe the measure $m_{G/\Gamma}$ by
  $m_{G/\Gamma} (A) = m_{G}(\Omega \cap 
  \pi^{-1}_{G}(A))$, and in particular, this
  formula does not depend on the choice of $\Omega$. The same facts
  hold for $H$ in place of $G$.

We will show that there is $c>0$ such that 
\begin{equation}\label{eq: show up to constant}
  \mu_{\sro} =
  c\, \varphi_*\left(
    m_{\crly{E}_n} \times m_{\R^d}|_{\bar B_{r_0}} \right)
\end{equation}
and then we will determine the constant. For the proof of \eqref{eq:
  show up to constant}, we first claim that the  
  product map
  $$\Psi: \bR\times \R^d \times H\to G, \ \ \Psi(t,v,h)\df
  a_tu(v)h$$
  pushes  the measure $m_\bR\times m_{\R^d}\times
    m_H$ to a multiple of 
$m_G|_{\mathbf{P}}$, where $\mathbf{P}$ is the image of $\Psi$. 
To see this, let $m_U$ be the Haar measure \index{M@$m_U$ --
  Haar measure on $U$} on $U$,
i.e., $m_U$ is the image of $m_{\R^d}$
under the map $v \mapsto u(v)$.
Note also that $m_Q \df m_\bR\times m_U$ is a left Haar measure on the
group $Q \df \set{a_t}\ltimes U$, and, since $H$ is unimodular, $m_H$ is a right
Haar measure on $H$. By Lemma~\ref{lem:hic basic}, $\nu \df \Psi_*(m_Q \times
m_H)$ is a multiple of $m_G|_{\mathbf{P}}$. This
proves the claim and moreover shows that $\on{supp}(\mu_{\sro}) = \sro$ is contained in $ \pi_G(\mathbf{P})$. 

In order to show \eqref{eq:
  show up to constant}, by Lemma \ref{lem:Sstar} we can restrict
attention to the restriction 
of $\mu_{\sro}$ to the set $\srosharp$. Recall the explicit formula
for $\mu_{\sro}$ given in Theorem
\ref{cor:csmeasure}\eqref{1204-3}. It suffices to check this
formula for small sets $E$ which generate the $\sigma$-algebra, and
for such small sets, both $E$ and $E^I$ can be lifted bijectively to
subsets of $\Omega \subset G$, where we can apply the description of
$m_G$ given in the previous paragraph. According to this description,
if $E_0 \subset \Psi (\R^d \times H)$ and $E_0^{I}$ is defined via \eqref{eq:
  thickening} and the action of $a_t$ on $G$, then
$$m_G(E_0^{I}) = m_{\R}(I) \cdot \Psi_*(m_{\R^d} \times m_H)(E_0).$$
From this
  we see that up to a constant multiple, 
  $\mu_{\sro}$ is 
  the pushforward under $\Psi$ of the measure $m_{\R^d} \times m_H$,
  restricted to $\psi(\srosharp)$ (where $\psi$ is as in~\eqref{eq: psi coordinates}). But since the measure $m_{\crly{E}_n}$ is
     given by restricting $m_H$ to a fundamental domain, the same
     description applies to the right-hand side of \eqref{eq: show up to constant}.
This implies \eqref{eq: show up to constant}.

We now claim that 
\begin{equation}\label{eq: zeta c}
  c = \frac{d}{\zeta(n)}.
  \end{equation}
We normalize the Haar measures used above by requiring $m_G(\Omega_G)
= m_H(\Omega_H)=1$, and similarly normalize $m_U$ by requiring that a
fundamental domain for $U/U(\Z)$ has measure one.  
The preceding discussion shows that $c$ is the scalar for which $m_G =
c \nu.$ Write 
$$
q(t, x_1, \ldots, x_d) \df \left(\begin{matrix} e^t & 0 & \cdots &
    e^tx_1 \\
    0 & e^t & \cdots & e^tx_2 \\
    & \cdots & \cdots & \\
  0 & 0 & \cdots & e^{-dt} \end{matrix} \right)
\in Q.
$$
Consider the map 
$$
\tau: Q \to \R^n, \ \ q= q(t, x_1, \ldots, x_d) \stackrel{\tau}{\mapsto} q
\mathbf{e}_n  = \left(e^tx_1, \ldots, e^tx_d,
  e^{-dt}\right)^{\mathbf{t}}$$
which is a bijection between $Q$ and
the upper half
  space
  $$\R^n_+ \df \{(y_1, \ldots, y_n)^{\mathbf{t}} \in \R^n : y_n>0\}.$$
By uniqueness of invariant measures on transitive spaces, and since
$\tau$ is $Q$-equivariant and the $Q$ action on $\R^n_+$ preserves
Lebesgue measure, the measures $\tau_* m_Q$  and $
m_{\R^n}|_{\R^n_+}$ are scalar multiples of each other. Considering the set
$$ E \df \left\{ \left (e^tx_1, \ldots,
e^tx_d, e^{-dt} \right)^{\mathbf{t}} : t \in [0,1], \forall i, \,
|x_i| \leq \frac{1}{2}\right \} \subset \R^n_+,$$
we see that $m_Q(\tau^{-1}(E))=1$ (where we have used that $m_U$
induces a probability 
measure on $U/U(\Z)$) and $m_{\R^n}(E)=d$ (by direct computation) and
thus
$m_{\R^n}|_{\R^n_+}  = d \cdot \tau_* m_Q.$

  Also, the map
  $$Q \times H \to G, \ \ (q,h) \mapsto qh$$
  is injective,
since for $g = qh$ we can reconstruct $q$ uniquely from the vector $q
\mathbf{e}_n$. It is easy to verify that the image $Q H \subset G$ is
open.

\ignore{
  For measurable subsets $A, B$ of a measure space $(X, \mu)$, write $A
 \sim_{\mu}B$ if $\mu(A 
 \triangle B)=0$. 
  Let $\sqcup$
denote disjoint union, and for each $h 
\in \Omega_H$, let
$$Q_h \df \{q \in Q: qh\in \Omega_G\}.$$
We claim that
\begin{equation}\label{eq: fundamental domains}
 \bigcup_{h \in \Omega_H} Q_hh = 
   \bigsqcup_{h \in \Omega_H} Q_hh\sim_{m_G} \Omega_G.
 \end{equation}

 Indeed, the fact that $Q_{h_1}h_1 \cap Q_{h_2}h_2 = \varnothing$ for $h_1 \neq h_2$
 follows from the fact that $Q \cap H = \{e\}$, and the inclusion
 $\subset$ in \eqref{eq: fundamental domains} is obvious. To show
 equality up to a set of $m_G$-measure zero, we will show
 that $\Omega_G \cap QH \subset \bigcup_{h \in \Omega_H} Q_hh$. So let
 $\omega = qh \in \Omega_G$ and let $h'  \in \Omega_H,
 \, 
\gamma' \in \Gamma_H$ such that $h = h'\gamma'$. Then $qh'\gamma' \in
\Omega_G$ and $qh \in \Omega_H$ 
\red{complete this step.}
}
Let $f = \mathbf{1}_E : \R^n \to \R$ be the indicator of $E$.
For $\Lam \in \XX_n$, write
  $$
\widehat{f}(\Lam) \df \sum_{v \in \Lam \sm \{0\}} f(v) \ \ \text{ and
} \ \widehat{f}^p(\Lam) \df \sum_{v \in \Lam_{\prim}} f(v).
  $$
By the Siegel summation formula \cite{Siegel_mean_value_theorem}, 
\begin{equation}\label{eq: SSF}
d = m_{\R^n}(E) = \int_{\bR^n} f \, dm_{\bR^n} = \int_{\XX_n}
\widehat{f} \, dm_{\XX_n} 
= \zeta(n) \int_{\XX_n} \widehat{f}^p \, dm_{\XX_n}.
\end{equation}

We define a lift of $f$ to $QH\subset G$ by
$$F : QH\to \R, \ \ F(qh) = f(\tau(q)
) \cdot
\mathbf{1}_{\Omega_H}(h).$$
It is easily checked that this definition implies
\begin{equation}\label{eq: lift and folding}
\sum_{\gamma \in \Gamma} F(g\gamma) = \hat{f^p} (g \Z^n).
  \end{equation}
Then, using that $m_H(\Omega_H)=1$, 
\begin{equation}\label{eq: lifted f}
  \int_G F \,
  d\nu = \int_{\Omega_H}\int_Q f(\tau(q))dm_Q(q) dm_H= \int_{Q}
  f\circ \tau \, dm_{Q}  
  =1.
\end{equation}
Using Fubini and `folding', we have 
\begin{equation}\label{eq: folding}
\int_G F \, dm_G = \sum_{\gamma \in \Gamma} \int_{\Omega_G \gamma} F
\, dm_G = \int_{\XX_n} \sum_{\gamma \in \Gamma} F(g\gamma) \,
 dm_{\XX_n}
\stackrel{(\ref{eq: lift and folding})}{=}
\int_{\XX_n}
\hat{f^p} \, dm_{\XX_n}.
\end{equation}
Recalling that, $m_G = c \nu$, and comparing 
\eqref{eq: SSF}, \eqref{eq:
  lifted f} and \eqref{eq: folding} gives \eqref{eq:
  zeta c}. 
\end{proof} 

\section{Special subsets of the section}\label{sec: special subsets}

The cross-section $\sro$ contains two subsets of Diophantine
significance, related respectively to best
approximations and to $\vre$-approximations. In this section we will
introduce these subsets, 
establish their Jordan measurability under suitable hypotheses, and
discuss their temperedness. 

\subsection{The set $\cB$ for best approximations}\label{subsec: 9.1}
Recall from \eqref{eq: vLam} that for $\Lam\in
\srosharp$ we denote by $v_\Lam\in \bR^d$ the horizontal component of the unique vector  
$v(\Lam)\in \Lam_{\on{prim}}\cap D_{r_0}$. Let \index{R@$r(\Lam)$ -- the length
  of $v_\Lam$} 
\begin{equation}\label{eq: def r Lam}
  r(\Lam) \df  \norm{v_\Lam}
  \end{equation}
be its
distance from the vertical axis, and let \index{B@$\cB$ -- the subset
  of $\sro$ corresponding to best approximations}
\begin{equation}\label{eq: def Best}
  \cB \df \set{\Lam\in  \srosharp : C_{r(\Lam)}\cap \Lam_{\on{prim}}
    = \set{\pm v(\Lam)}},
  \end{equation}
  where $C_r$ is defined in \eqref{eq: def C r}.
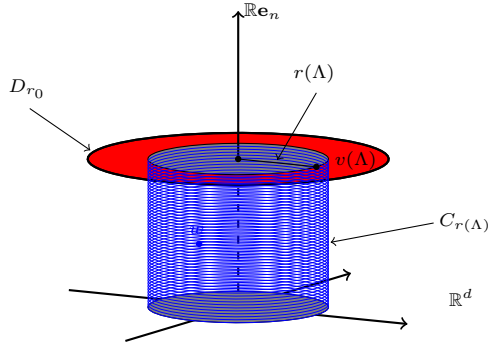
\begin{figure}[htbp]
\center{
\tdplotsetmaincoords{80}{30}
\begin{tikzpicture}[tdplot_main_coords,   scale=2]

\draw[black, thick, ->] (-1.3,0,0) -- (1.3,0,0);
\draw[black, thick, ->] (0,-1.5,0) -- (0,1.5,0);
\draw[black, thick, dashed] (0,0,0) -- (0,0,1);
\draw (1,1,0) node[anchor=west]{\!\tiny{$\bR^d$}};

\draw [thick, draw=black, fill=red, opacity=0.2] (0,0,1) circle (1);
\draw [thick, draw=black] (0,0,1) circle (1);

\filldraw [thin, draw=black, fill=gray, opacity=0.8] (0,0,1) circle (0.6);
\filldraw [thin, draw=black, fill=gray, opacity=0.8] (0,0,0) circle (0.6);

\draw[black, thick, ->] (0,0,1) -- (0,0,2) node[anchor=west]{\!\tiny{$\bR\mb{e}_n$}};
\tdplotsetrotatedcoords{0}{0}{0}
\begin{scope}[tdplot_rotated_coords]
\foreach \y in {0,0.02,...,1}
{\tdplotdrawarc[blue, very thin, opacity = 0.3]{(0,0,\y)}{0.6}{0}{360}{}{}{}}
\end{scope}

\filldraw [black] (0,0,1) circle (0.5pt);

\filldraw [black] (0.6,0,1) circle (0.5pt);
\draw (0.7,0,1.05) node[anchor=west]{\!\tiny{$v(\Lam)$}};

\draw[black,  thin, ->] (1.5,0,0.7)--(0.71,0,0.5);
\draw (1.5,0,0.7) node[anchor = west]{\!\tiny{$C_{r(\Lam)}$}};

\draw[black, thin] (0,0,1)--(0.6,0,1);
\draw[black, thin, ->] (0.6, 0, 1.5)--(0.3, 0 , 1.02);
\draw (0.6, 0, 1.5) node[anchor = south]{\!\tiny{$r(\Lam)$}};
\draw[black,  thin, ->] (-1.6,0,1.2)--(-1.02,-0.2,1);
\draw (-1.6,0,1.2) node[anchor = south]{\!\tiny{$D_{r_0}$}};

\filldraw [blue] (-0.3, 0, 0.4) circle (0.5pt) node[anchor = south]{\!\tiny{$w$}};

\end{tikzpicture}
}

\caption{If the cylinder
  $C_{r(\Lam)}$ defined by the unique vector $v(\Lam) \in \Lam \cap D_{r_0}$
  contains another lattice point $w$, then $\Lam \notin
  \cB$. }
\label{fig: best and not best}
\end{figure}

The set $\cB$ will detect best approximations. It consists of
lattices in the cross-section $\sro$ with a unique vector 
in $D_{r_0}$, such that the cylinder $C_{r(\Lam)}$ they define
contains no lattices points besides $\{0, \pm v(\Lam)\}$ (see Figure
\ref{fig: best and not best}). Note that
$\cB$ depends on the norm, although this is not reflected in the 
notation.

\begin{lemma}\label{lem:boundary of best}
The set $\cB$ is open in $\sro$. The boundary $\partial_{\sro}\cB$ is
contained in the union of $\XX_n(D_{r_0},2)$ and  
{\small
\begin{equation}\label{eq: def E union}
  \crly{Z} \df \set{\Lam\in \XX_n : 
  \begin{array}{lll}
  \exists v,w\in \Lam
  \text{ not co-linear s.t. }\\
  \on{span}_\bZ\set{v,w} = \on{span}_{\bR}\set{v,w}\cap \Lam\\
  \text{ and } 
  \norm{ \pi_{\bR^d}(v)} = \norm{\pi_{\bR^d}(w)}
  \end{array}
  }.
\end{equation}
}
\end{lemma}

\begin{proof}
  We show $\cB$ is open by showing its complement is closed. 
Let
$$K \df
\left\{ \Lam \in \srosharp:
\#(C_{r(\Lam)} \cap \Lam_{\prim} )\geq 3
\right\}
;$$
 that is, $K$ consists of the lattices $\Lam\in
 \sro^{\sharp}$ 
 such that  the cylinder $C_{r(\Lam)}$ contains at least one extra
 primitive vector besides $\pm v(\Lam)$. 
By definition $$\sro\smallsetminus \cB = \XX_n(D_{r_0},2) \cup K .$$
By Lemma
\ref{lem:continuity of vector} $\XX_n(D_{r_0},2)$ is closed, so we let $\Lam_i
\in K$ be a sequence converging  
to some $\Lam$, and show that $\Lam\in
\XX_n(D_{r_0},2) \cup K$.  

By Lemma
\ref{lem:continuity of vector}, $\Lam \in \XX_n(D_{r_0})$, and suppose
that $\Lam\notin \XX_n(D_{r_0},2)$, so that $\Lam \in
\sro^\sharp$.
Since $\Lam_i \in K$, the cylinders
$C_{r(\Lam_i)}$ contain at least three primitive vectors. By
Lemma~\ref{lem:continuity of vector}, $v(\Lam_i)\to v(\Lam)$ and so
$r(\Lam_i)\to r(\Lam)$. Let $r'>r(\Lam)$, then for 
all large enough $i$, $r(\Lam_i)<r'$ and so for all large enough $i$,
$\Lam_i\in \XX_n(C_{r'}, 3)$. Using Lemma~\ref{lem:continuity of
  vector} again we see that  $\Lam\in
\XX_n(C_{r'}(1), 3)$. Since this is true for all $r'>r(\Lam) $ we deduce
that  
$\Lam\in \XX_n(C_{r(\Lam)},3)$, and so $\Lam\in K$. 

We now prove the second assertion. Let 
$$U_0 \df \set{\Lam\in \sro^{\sharp} : \Lam_{\on{prim}}\cap
  C_{r(\Lam)}^\circ \ne \varnothing}$$
(where $C_r^\circ$ denotes the interior of $C_r$). 
We will shortly show that $U_0$ is open in $\sro$. 
From this, and since $\sro\smallsetminus \cB$ is closed, it
follows that 
\begin{equation*}
\partial_{\sro}\cB = \partial_{\sro}(\sro\smallsetminus \cB) \subset
\pa{\XX_n(D_{r_0},2)\cup K} \smallsetminus U_0. 
\end{equation*}
We have that $K\smallsetminus U_0\subset \crly{Z}$ because once $\Lam\in K\smallsetminus  U_0$, apart from the vector $v=v(\Lam)$ which defines the cylinder $C_{r(\Lam)}$, there must exist another primitive vector 
$w\in \Lam$, not co-linear to $v$, on the boundary of $C_{r(\Lam)}$ outside of $D_{r_0}$. It follows that $w$ satisfies 
$\norm{\pi_{\bR^d}(v)}=\norm{\pi_{\bR^d}(w)}$. 
Moreover, since the triangle given by the convex hull of $0,v,w$ is contained in $C_{r(\Lam)}^\circ$, and thus is free of lattice points (apart from its its vertices), 
we deduce that $v,w$ satisfy 
$\on{span}_\bZ\set{v,w} = \on{span}_\bR\set{v,w}\cap \Lam$ and so $\Lam\in \crly{Z}$ as claimed. This implies the second
assertion of the Lemma, once we show that $U_0$ is open. 

To prove that $U_0$ is open, let $\Lam_i$ be a sequence in
$\XX_n(D_{r_0})$ that converges to a lattice $\Lam\in U_0$. We need to
show 
that for all large enough $i$, $\Lam_i\in U_0$. Since $\Lam$ does
not belong to the closed set $\XX_n(D_{r_0},2)$,  $\Lam_i \notin
\XX_n(D_{r_0},2)$ for all large $i$. By Lemma~\ref{lem:continuity
  of vector}, $v(\Lam_i)\to v(\Lam)$ and  
thus $r(\Lam_i)\to r(\Lam)$. 
Since there exists a primitive vector $w\in \Lam_{\on{prim}}$ in the
open set $C_{r(\Lam)}^\circ$, there is some $r'<r(\Lam)$ 
such that $w\in C_{r'}^\circ$. In other words, $\Lam\in
\XX_n(C_{r'}^\circ)$. The latter set is open 
in $\XX_n$ by Lemma~\ref{lem:open
  sets}, and hence  $\Lam_i\in \XX_n(C_{r'}^\circ)$ for all large
$i$. But for all large $i$ we also have  
$r(\Lam_i)>r'$, so $\Lam_i$ must contain a
primitive vector in the open cylinder $C_{r(\Lam_i)}^\circ$, i.e., 
$\Lam_i\in U_0$.
\end{proof}

\begin{lemma}\label{lem: cB null 1}
In Case I, for any norm on $\R^d$, the set $\cB$ is $\mu_{\sro}$-JM and 
$\mu_{\sro}(\cB)>0$. 
\end{lemma}

\begin{proof}
The fact that $\mu_{\sro}(\cB)>0$ follows from the openness of $\cB$ in $\sro$ (see Lemma~\ref{lem:boundary of best}) and the fact that $\mu_{\sro}$ has full support by Proposition \ref{prop:description of measures1}.

By Lemmas \ref{lem:Sstar} and \ref{lem:boundary of best}, it is enough
to show that $\mu_{\sro}(\crly{Z})=0,$ for $\crly{Z}$ as in \eqref{eq: def E
  union}. Since $\crly{Z}$ is $\{a_t\}$-invariant, 
this is equivalent 
to showing that $m_{\XX_n}(\crly{Z})=0$.
Let $G = \SL_n(\R)$ and let \index{M@$m_G$ -- Haar measure on $G =
  \SL_n(\R)$} $m_G$ denote the Haar measure on $G$, so
that $m_{\XX_n}$ is the restriction of $m_G$ to a fundamental
domain. For $g\in G $ denote the
column vectors of $g$ by $g_1,\dots, g_n$, and for  
$v\in \bR^n$, denote $\bar{v} = \pi_{\bR^d}(v)$. Let 
$$\widetilde{\crly{Z}} \df \set{g\in G: \norm{\bar{g}_1} =
  \norm{\bar{g}_n}}.$$
Because
the image of $\widetilde{\crly{Z}}$ in $\XX_n$ under the projection $g\mapsto
g\bZ^n$ contains $\crly{Z}$, it is enough to show
  that $m_G(\widetilde{\crly{Z}})=0$. To see this, note that $m_G$ is
  invariant under right-multiplication by elements 
of $\{a_t\}$, so any set whose $m_G$ measure is positive must
intersect some $\{a_t\}$-orbits for this action along a positive measure's
worth of parameters $t$. On the other hand, the orbit of each $g \in G$ under
this right multiplication only intersects 
$\widetilde{\crly{Z}}$ once. 
\end{proof}
\begin{lemma}\label{lem: positivity of B case 2}
In case II with $\mu = m_{\vec{\al}}$, $\mu_{\sro}(\cB)>0$.
\end{lemma}
\begin{proof}
Assume by way of contradiction that $\mu_{\sro}(\cB) = 0$.
 Since $\cB$ is open in 
$\sro$ by Lemma~\ref{lem:boundary of best}, it follows from Proposition~\ref{prop:description of measures}\eqref{item: 1.2},
that $\bar{A}_{\vec{\al}}y_{\vec{\al}}\cap \cB = \varnothing$. The contradiction follows since it is easy to see that for any
lattice $\Lam$ without non-zero vectors on the vertical axis, the trajectory $\set{a_t\Lam:t>0}$ always visits $\cB$.

\end{proof}
Let $\bar{A}_{\vec{\al}} $ be the group defined in \eqref{eq: def bar
  A}. \index{A@$\bar{A}_{\vec{\al}}$-analytic norm}
\begin{definition}\label{def: a bar analytic}
We say that the norm $\| \cdot \|$ on $\R^d$ is {\em $\bar{A}_{\vec{\al}}$-analytic}, if
for any $v, w \in \R^n$, the set
$$
\bar
A_{v,w} \df \left\{ \bar a \in \bar{A}_{\vec{\al}}: \| \pi_{\R^d}(\bar a v)\| = \|
  \pi_{\R^d}(\bar a w)\|  \right \}
$$
is an analytic subset of $\bar{A}_{\vec{\al}} \cong \R^d; $ that is, the  zero-set of an analytic
function.
\end{definition}

For example, the Euclidean norm is analytic. 
\begin{lemma}\label{lem: cB null 2} 
In Case II, if the norm is $\bar{A}_{\vec{\al}}$-analytic then the set $\cB$ is
$\mu_{\sro}$-JM. 
\end{lemma}

\begin{proof}
As in the proof of Lemma \ref{lem: cB null 1}, we need to show that
$m_{\vec{\alpha}}(\crly{Z})=0.$
Let $\bar h_{\vec{\al}}, \, y_{\vec{\al}}$ and $\bar A_{\vec{\al}}^{(1)}$ be as in \eqref{eq:
  def Bar B}, \eqref{eq: def bar A} and 
\eqref{eq: def bar A 1}. Let $m_{\bar{A}_{\vec{\al}}}$ denote the Haar measure
on $\bar{A}_{\vec{\al}}$, so that $m_{\vec{\alpha}}$ is the pushforward under the
map $\bar a \mapsto \bar a y_{\vec{\al}}$, of the
restriction of $m_{\bar{A}_{\vec{\al}}}$ to a measurable set.

Suppose by way of contradiction that $m_{\vec{\alpha}}(\crly{Z})>0$. Then
there are fixed $v, w \in (y_{\vec{\al}})_{\prim}$ such that $v \neq \pm w$, and
$$
m_{\bar{A}_{\vec{\al}}} \left( \bar A_{v,w} \right)> 0.
$$
Since the norm $\| \cdot \|$ is $\bar{A}_{\vec{\al}}$-analytic, we must have $\bar A_{\vec{\al}}= \bar A_{v,w}$. 
    Let 
\begin{equation*}\label{eq: bar A eigenbasis}    \mathbf{t}_i = \bar h_{\vec{\al}} (\mathbf{e}_i), \ \
  i=1, \ldots, n
  \end{equation*}
    be a basis of $\R^n$ consisting of simultaneous eigenvectors for
    $\bar{A}_{\vec{\al}}$. 
%
%
Write
$$v = \sum_{i=1}^n r_i
\mathbf{t}_i, \quad w = \sum_{i=1}^n s_i\mathbf{t}_i.$$
By acting with
$$\bar a_t \df \bar h_{\vec{\al}} \diag{e^{dt}, e^{-t}, \ldots, e^{-t}}
\bar h_{\vec{\al}}^{-1} \subset \bar{A}_{\vec{\al}}
= \bar A_{v,w}$$
we see that 
that $|r_1| =
|s_1|$; indeed,
$$\frac{|r_1|}{|s_1|}   = \frac{\|\bar a
  _t (r_1\mathbf{t}_1)\|}{ \| \bar a_t (s_1\mathbf{t}_1)\|} =
\frac{\| \pi_{\R^d}(\bar a_t v)\| + o(t)}{\|\pi_{\R^d}(\bar a_t w)
  \| + o(t)}  = 1+o(t).$$

In particular, the first coordinate (with respect to the eigenbasis
$\mathbf{t}_i$) 
of $v \pm w \in y_{\vec{\al}}$ is zero. Therefore the nonzero vector 
$\bar h_{\vec{\al}}^{-1}(v\pm w) \in (x_{\vec{\al}}^*)_{\prim}$ has one of its
coordinates equal to zero, which contradicts Proposition \ref{prop:
  compact A orbit}.
\end{proof}

\begin{lemma}\label{lem: cB null 3} 
In Case II, if the norm on $\R^d$ is the sup-norm, then the set $\cB$ is
$\mu_{\sro}$-JM. 
\end{lemma}

\begin{proof}
  Once again  we need to show that
$m_{\vec{\alpha}}(\crly{Z})=0.$ Supposing by contradiction that
$m_{\vec{\alpha}}(\crly{Z})>0$, and using the explicit description of
$m_{\vec{\alpha}}$ given in Propositions \ref{prop: compact A orbit} and \ref{prop:
  explicit compact push}, we see that there is a lattice $\Lam_1$ of
type $\pmb{\sig}$, and two linearly independent
primitive vectors $v, 
w \in \Lam_1$, such that for a set of
positive measure of $\mathbf{t} = (t_1, \ldots, t_d) \in \R^d$,
$$
\left\| \pi_{\R^d}(\bar h_{\vec{\al}} a_{\mathbf{t}} 
  v )\right \| = \left  \| \pi_{\R^d}( \bar h_{\vec{\al}}
  a_{\mathbf{t}} 
  w )\right \|,  
\ \ \text{ where } a_{\mathbf{t}} \df \diag{e^{t_1}, \ldots, e^{t_d},
  e^{-\sum_{j=1}^d t_j}}.$$
Using \eqref{eq: homothetic to} and \eqref{eq: def Bar B}, we obtain
that for some $\beta_v, \beta_w 
\in \bK$, which are linearly independent
over $\Q$, we have  
  $$
  \left\| \left(
    \begin{matrix} \sum_{j=1}^d e^{t_j} (\sigma_j(\alpha_1)-\alpha_1)
      \sigma_j(\beta_v) \\ \vdots \\  \sum_{j=1}^d   e^{t_j} 
      (\sigma_j(\alpha_d)-\alpha_d)\sigma_j(\beta_v)
    \end{matrix} \right ) \right\| = \left \| 
  \left( \begin{matrix} \sum_{j=1}^d e^{t_j}  (\sigma_j(\alpha_1)-\alpha_1)
      \sigma_j(\beta_w) \\ \vdots \\ \sum_{j=1}^d e^{t_j} 
      (\sigma_j(\alpha_d)-\alpha_d)\sigma_j(\beta_w)
      \end{matrix} \right) \right \|. 
  $$
  Since we are working with the sup-norm, this implies that there are
  indices $1 \leq k_1, k_2\leq d$ 
and $\omega = \pm 1$ such that for a set of positive measure of
$\mathbf{t}$,
$$
\sum_{j=1}^d e^{t_j} (\sigma_j(\alpha_{k_1})-\alpha_{k_1}) \sigma_j(\beta_v) 
= \omega \sum_{j=1}^d e^{t_j }(\sigma_j(\alpha_{k_2})-\alpha_{k_2})
\sigma_j(\beta_w) . 
$$
This is an equality of analytic expressions, which holds for
$\mathbf{t}$ in a set of positive measure, and thus it must hold for
all $\mathbf{t}$. In particular, we can take partial derivatives
$\frac{\partial}{\partial t_j}|_{\mathbf{t}=0}$ to obtain that for any
$j = 1, \ldots, d$  we have
\begin{equation}\label{eq: 68}
(\sigma_j(\alpha_{k_1})- \alpha_{k_1})\sigma_j(\beta_v) =\omega
(\sigma_j(\alpha_{k_2})- \alpha_{k_2})\sigma_j(\beta_w). 
\end{equation}
If $k_1 = k_2$ it follows that
$\beta_v = \omega \beta_w$, which is a contradiction. 
Assume therefore that $k_1 \neq  k_2$. Multiplying \eqref{eq: 68} by
$\sigma_j(\beta_w^{-1})$ 
and letting $\beta \df \beta_v/\beta_w$ we get
\begin{equation}\label{eq: 69}
(\sigma_j(\alpha_{k_1})-\alpha_{k_1}) \sigma_j(\beta) = \omega
(\sigma_j(\alpha_{k_2}) - \alpha_{k_2}).
\end{equation}
Recall our convention $\sigma_n = \mathrm{Id}$, which implies that for
$\gamma \in \bK$, the trace $\mathrm{Tr}$ satisfies 
$$
  \sum_{j=1}^d \sigma_j (\gamma) = \mathrm{Tr}(\gamma)-\gamma.
$$
Summing \eqref{eq: 69} over $j$ gives 
$$
\mathrm{Tr}(\alpha_{k_1}\beta) - \alpha_{k_1}\beta - \alpha_{k_1}
(\mathrm{Tr}(\beta) - \beta) = \omega (\mathrm{Tr}(\alpha_{k_2}) -
\alpha_{k_2} - d \alpha_{k_2}),
$$
and thus
$$
\mathrm{Tr}(\alpha_{k_1}\beta - \omega \alpha_{k_2}) -
\mathrm{Tr}(\beta) \alpha_{k_1} + n \alpha_{k_2}=0
$$
and this is a nontrivial linear dependence over $\Q$ between $1,
\alpha_{k_1}, \alpha_{k_2}$. This is a contradiction. 
  \end{proof}
\begin{remark}  There are other  norms for which the 
conclusion of Lemmas \ref{lem: cB null 2} and \ref{lem: cB null 3}
fails. Here is a sketch of 
how one can build an example. Let $d=2$. 
We define a norm, a
lattice $\Lam$, two  vectors $v,w \in \Lam$, and positive $r,
\vre$, so that:
\begin{enumerate}[(a)]
  \item $v, w$ are primitive in $\Lam$ and linearly independent;
\item \label{item: assume a}
$\Lam$ arises from a number field as in \S
\ref{subsec: case II};
\item \label{item: assume b}
  $v, w$ are both in 
  $\partial C_r$, with $C_r^\circ \cap \Lam = \{0\}$;
  \item \label{item: assume c}
for all $t \in (-\vre, \vre)$, 
$$
\|g_t \pi_{\R^2}(v) \| = \|g_t \pi_{\R^2}(w)\|, \ \ \text{ where }
g_t \df \diag{e^t,
  e^{-t}}. 
$$
\end{enumerate}
One can then see from Proposition
\ref{prop:description of measures}\eqref{item: 1.25}, that in this
case the boundary of $\cB$  
must have positive $\mu_{\sro}$-measure.

To this end, we first construct  
a norm $\| \cdot \|$, two linearly independent vectors $u_1, u_2$
in $\R^2$, and positive $\vre, r$, such that $r = \|g_t u_1 \| =
\|g_t u_2 \| $ for all $t \in (-\vre, \vre)$. 
Let 
 $$u_1 = \left( \begin{matrix} 1 \\  1/4\end{matrix}
\right), \ \  u_2 = \left(\begin{matrix} 2/3 \\
    2/3 \end{matrix} \right).$$
We will define the norm in $\R^2$ by specifying  a symmetric convex
body $\mathbf{B}$ which is its unit ball. The 
boundary of $\mathbf{B}$ 
consists of four small smooth arcs $\pm \gamma_i, i=1,2$, where
$\gamma_i$ passes through $ u_i,$ and four line segments
connecting the ends of these arcs. See Figure \ref{fig: for Uri 1}.
\begin{figure}[htbp]
\center{

\begin{tikzpicture}[ scale=2]
\draw[black, thick, ->] (-1.2,0)--(1.2,0) node[anchor = west]{$x$};
\draw[black, thick,->] (0,-1)--(0,1)node[anchor = west]{$y$};
\draw[blue, thin, ->] (0,0)--(1,1/4) node[anchor = west]{\!\tiny{$u_1$}};
\draw[blue, thin, ->] (0,0)--(-1,-1/4) node[anchor = east]{\!\tiny{$-u_1$}};
\filldraw[blue] (1,1/4) circle (0.5pt);
\filldraw[blue] (-1,-1/4) circle (0.5pt);
\draw[red, thin, ->] (0,0)--(2/3,2/3) node[anchor = west]{\!\tiny{$u_2$}};
\draw[red, thin, ->] (0,0)--(-2/3,-2/3) node[anchor = east]{\!\tiny{$-u_2$}};
\filldraw[red] (2/3,2/3) circle (0.5pt);
\filldraw[red] (-2/3,-2/3) circle (0.5pt);


\draw[blue] (0.97, {sqrt(17/16 - (0.97)^2)}) .. controls (1,1/4) .. (1.03,{sqrt(17/16 - (1.03)^2)});

\draw[blue] (-0.97, {-sqrt(17/16 - (0.97)^2)}) .. controls (-1,-1/4) .. (-1.03,{-sqrt(17/16 - (1.03)^2)});

\draw[red] (0.5, {sqrt(8/9 - (0.5)^2)}) .. controls (2/3,2/3) .. (0.73,{sqrt(8/9 - (0.73)^2)});

\draw[red] (-0.5, {-sqrt(8/9 - (0.5)^2)}) .. controls (-2/3,-2/3) .. (-0.73,{-sqrt(8/9 - (0.73)^2)});

\filldraw[black] (0.97, {sqrt(17/16 - (0.97)^2)}) circle (0.2pt);
\filldraw[black] (1.03,{sqrt(17/16 - (1.03)^2)}) circle (0.2pt);

\filldraw[black] (-0.97, {-sqrt(17/16 - (0.97)^2)}) circle (0.2pt);
\filldraw[black] (-1.03,{-sqrt(17/16 - (1.03)^2)}) circle (0.2pt);

\filldraw[black] (0.5, {sqrt(8/9 - (0.5)^2)}) circle (0.2pt);
\filldraw[black] (0.73,{sqrt(8/9 - (0.73)^2)}) circle (0.2pt);

\filldraw[black] (-0.5, {-sqrt(8/9 - (0.5)^2)}) circle (0.2pt);
\filldraw[black] (-0.73,{-sqrt(8/9 - (0.73)^2)}) circle (0.2pt);

\draw[black] (0.97, {sqrt(17/16 - (0.97)^2)}) -- (0.73,{sqrt(8/9 - (0.73)^2)}) ; 
\draw[black] (0.5, {sqrt(8/9 - (0.5)^2)}) -- (-1.03,{-sqrt(17/16 - (1.03)^2)}) ; 

\draw[black] (-0.97, {-sqrt(17/16 - (0.97)^2)}) -- (-0.73,{-sqrt(8/9 - (0.73)^2)}) ;
\draw[black] (-0.5, {-sqrt(8/9 - (0.5)^2)}) -- (1.03,{sqrt(17/16 - (1.03)^2)}) ;

\end{tikzpicture}
}

\caption{Defining the norm by a carefully chosen convex set in $\R^2$. }
\label{fig: for Uri 1}
\end{figure}
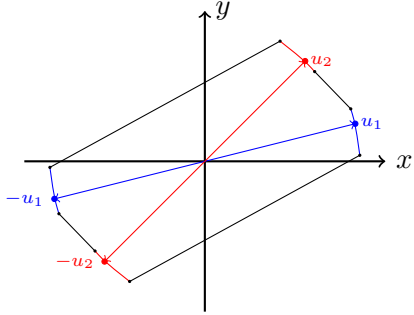
The
$\gamma_i$ are carefully chosen 
so that \eqref{item: assume c} is satisfied, for some small $\vre>0$,
and for $u_1 = \pi_{\R^2}(v), \ u_2 = \pi_{\R^2}(w)$. Since $\gamma_i$
passes through $u_i$ and determines the norm, the 
requirement \eqref{item: assume c} implies that $\gamma_1$ uniquely
determines $\gamma_2$, and we have to choose $\gamma_1$ so that 
the resulting figure is convex. This 
can be shown using an explicit computation in polar
coordinates. 
Moreover the computation shows that in this construction, we have
freedom to vary $u_1, u_2$ in some open set in 
$\R^2 \times \R^2$. 

It is not hard to construct a lattice $\Lam \subset \R^3$ containing
two primitive linearly independent vectors
$v,w$ such that $\pi_{\R^2}(v) = u_1, \pi_{\R^2}(w) = u_2$. Moreover
it is not hard to choose $\Lam$ so that \eqref{item: assume b} holds,
and again this can be carried out for $u_1, u_2$ in some open set in 
$\R^2 \times \R^2$. 
Using the fact that
number field lattices 
are dense (up to a rescaling) in the space of 
lattices, one sees that one can also arrange that \eqref{item: assume
  a} holds.

 \end{remark}

Recall Definition~\ref{def: tempered}. The following temperedness result will be important for our
analysis. For related results, see \cite{Lagarias2}, 
\cite{Cheung_Chevallier} and \cite{Chevallier}. 
\begin{proposition}\label{prop: B tempered}
The set $\cB$ is tempered. 
\end{proposition}

\begin{proof}
  Assume otherwise. Recalling the notation \eqref{eq: def
C r}, let $M$ be large enough so that any set of cardinality $M+1$
in $C_{1}(e^d) = \set{\smallmat{v\\c}\in \bR^n : \norm{v}\le 1, \av{c}
  \leq e^d}$ contains distinct points $\smallmat{v\\c}, 
\smallmat{v'\\c'}$ such that 
$\norm{v-v'}<\frac{1}{3}$ and $ \av{c-c'}<1$. Such a number $M$ exists by
the compactness of $C_{1}(e^d)$. By applying a linear transformation which dilates the
horizontal subspace, one sees that for any $r>0$, in any subset of
$C_r(e^d)$ of cardinality $M+1$, there are distinct points $\smallmat{v\\c}, 
\smallmat{v'\\c'}$ such that 
$\norm{v-v'}<\frac{r}{3}$ and $ \av{c-c'}<1.$

We claim that $\cB$ is
$M$-tempered.
Indeed, suppose by contradiction that one can find 
$\Lam\in \cB$ and 
$0=t_0<t_1<\dots<t_{M} \leq 1$ such that $a_{t_j}\Lam\in \cB$ for $0\le j\le
M$. Then, for each $j$ the vector $w_j \df v(a_{t_j}\Lam)$  (see \eqref{eq: well defined vector}) 
satisfies 
$$w_j \in a_{t_j} \Lam_{\mathrm{prim}} \cap D_{r_0} \  \text{ and } a_{t_j}
\Lam \cap C_{\|\pi_{\R^d}(w_j)\|}^\circ = \{0\}.$$
Applying $a_{-t_j}$, we find 
 vectors $a_{-t_j}w_j = \smallmat{v_j\\e^{dt_j}}\in \Lam$ such that
 
 $$\Lam \cap C^{(j)} =\{0\}, \ \ \text{ where } C^{(j)} \df C^\circ_{\norm{v_j}}( e^{dt_j}). 
 $$
In particular, this implies 
$\norm{v_0} \geq \norm{v_1} \geq  \cdots  \geq \norm{v_M},$
and thus
$$\smallmat{v_j\\e^{dt_j}}
\in C_{\|v_0\|}(e^d) \ \ \text{ for } j=0, \ldots, M.
$$
By the property of $M$, there are indices 
$j_1<j_2$ such that $\norm{v_{j_1}-v_{j_2}}<\frac{\|v_0\|}{3} $ and
$|e^{dt_{j_2}}-e^{dt_{j_1}}| < 1$. The difference  
$$w = \smallmat{v_{j_2}\\e^{t_{j_2}}} - \smallmat{v_{j_1}\\e^{t_{j_1}}}$$
then belongs to the lattice $\Lam$ but also to the interior of the
cylinder 
$C_{\norm{v_0}}( 1)$. If $w$ is primitive we get a contradiction to
the assumption  $\Lam\in \cB$, and if $w$ is not primitive there is a
multiple $w' = tw, t \in (0,1)$ which is primitive and belongs to
$\Lam$,  again giving the desired contradiction.  
\end{proof}

\subsection{The set $\cS_\vre$ for
  $\vre$-approximations}\label{subsec: 9.2}
\begin{lemma}\label{lem: S eps null 1}
 In both Cases I and II,  for any
$0<r \le r_0$, the sets $\cS_r$ are $\mu_{\sro}$-JM. 
\end{lemma}

\begin{proof}
By Lemma~\ref{lem:continuity of vector}, $\XX_n(D_r)$ is closed in
$\XX_n$ and hence in $\XX_n(D_{r_0})$.  
Again by Lemma~\ref{lem:continuity of vector},
$\sro^{\sharp}\cap \XX_n(D_r^\circ) $ is open in
$\XX_n(D_{r_0})$. It follows that
$$\partial_{\sro}\cS_r \subset \XX_n(D_r)\smallsetminus
\pa{\sro^{\sharp}\cap \XX_n(D_r^\circ)}
\subset 
\XX_n(D_{r_0},2) \cup \XX_n(D_r\smallsetminus D_r^\circ)
.$$
 By Lemma~\ref{lem:Sstar}, $\mu_{\sro} \left(\XX_n(D_{r_0},2)
 \right)=0$, and hence, by 
Theorem~\ref{cor:csmeasure}, it suffices to show that
$$ \mu \left(
  \XX_n(M_r
  )\right) =0,$$
where $M_r = (D_r\smallsetminus D_r^\circ)^\bR$ is as in \eqref{eq:
  def M r}. 
In Case I, this follows from Lemma
\ref{lem:periodic-measures-nullity} as in the proof of Lemma \ref{lem:
  boundary estimate}; in Case II, 
this follows from Lemma \ref{lem: using convexity}. 
\end{proof}

\begin{proposition}\label{prof: not tempered}
For any $d>1$, for any $\vre>0$ and any norm, $\cS_\vre$ is not
tempered. For $d=1$, and any $\vre>0$, $\cS_\vre$ is tempered.  
\end{proposition}

\begin{proof}
Suppose $d>1$. Given $M$ and $\vre$, let $\Lam \in \XX_n$ be a lattice 
containing the primitive pair $\bu \df \frac{1}{2M} \mathbf{e}_n$ and
$\bv \df (v,1)^{\mathrm{t}}$, where $v \in \R^d$ satisfies
$\|v\|<\frac{\vre}{2}$ (recall $\bu, \bv$ are a {\em primitive pair} if
$\Lam \cap \mathrm{span}_{\R} (\bu, \bv) = \mathrm{span}_{\Z}(\bu,
\bv)$). Such a lattice exists because $d>1$. Then 
$\Lam$ contains the (necessarily primitive) vectors 
$$ \bv_j \df \bv +j \bu = 
\left( \begin{matrix} v \\ 1+ \frac{j}{2M} \end{matrix} \right) , \ \ \
j=0, \ldots, M.$$
Let $t_j\df  \frac{1}{d}\log\left(1+ \frac{j}{2M} \right) \in [0,1)$. This choice
ensures that $a_{t_j} \Lam$ contains the
vector 
$$a_{t_j} \bv_j = \left( \begin{matrix} e^{t_j}v \\ 1\end{matrix} \right)$$ 
which shows that $a_{t_j}\Lam \in
\cS_\vre$. Since $M$ is arbitrary, $\cS_\vre$ is not tempered.

Now suppose $d=1$. By  Minkowski's second theorem, there is 
$\kappa >0 $ satisfying the following. For any lattice $\Lam \in \XX_2$, denote by
$v_1$ a shortest nonzero vector of $\Lam$, and by $v_2$ a shortest
vector such that $v_1, v_2$ are linearly independent. Then
$\|v_2\|\geq \frac{\kappa}{\|v_1\|}.$ Thus, for any $C>0$ there is
$c>0$ so that if $\|v_1\| < c$, 
then there is no $v_2 \in \Lam_{\prim} \sm \{\pm v_1\}$ with $\|v_2\|
< C. $ Now, given $\vre>0$, choose $C$ large enough so that the ball of
radius $ C$ around the origin contains the rectangle $\mathbf{R} \df
[-\vre, \vre] \times [1, e]$.  Let $c$ be the corresponding constant,
and choose $M$ large enough so that any 
$M$ points in $\mathbf{R}$ contain a pair of distinct points of
distance less than $c$. With this choice, $\cS_\vre$ is
$M$-tempered. Indeed, if this were not the case, there would be a
lattice $\Lam \in \XX_2$ and $0= t_0  \leq t_1 < \cdots < t_M
\leq 1$ such that $a_{t_j}\Lam \in \cS_\vre$. This implies that $\Lam$
contains $M+1$ primitive vectors in $\mathbf{R}$, and hence one of
their differences, which is a nonzero vector of $\Lam$, has length
less than $c$. But by the choice of $c$ and $C$ we get that
$\mathbf{R} \cap \Lam_{\prim}$ contains only two vectors $\pm v_1$. 
\end{proof}

\subsection{The adelic space, cross-section and cross-section
  measure}\label{subsec: adelic lifts}
  We now lift the cross-section $\sro$ to the adelic space, and, using
  the theory developed so far,  
  derive properties crucial to our discussion.
  
Let $\XXnA$ be the adelic space, let $\pi: \XXnA \to \XX_n$ be the
projection, and let $m_{\XXnA}$ and 
$m_{\tilde{L}_{\vec{\al}} \tilde{y}_{\vec{\al}}} $ be the
measures introduced in \S \ref{subsec: adeles 
  case 1} and \S \ref{subsec: adeles case 2} respectively. Let
\index{S@$\wt{\cS}_\vre$ -- the lifted subset for
  $\vre$-approximations} \index{B@$\wt{\cB}$ -- the lifted subset for
  best approximations}
\begin{equation}\label{eq: lifted cross section}
  \widetilde{\cS}_{r_0} \df \pi^{-1}(\sro),
 \ \ \widetilde{\cB} \df \pi^{-1}(\cB), \ \ \text{ and } \widetilde{\cS}_\vre
 \df \pi^{-1}(\cS_\vre). 
 \end{equation}

\begin{theorem} \label{thm: MKf}
Let $\mu$ be $m_{\XX_n^{\bA}}$ (Case I) or $\mu = m_{\tilde{L}_{\vec{\al}}\tilde{y}_{\vec{\al}}}$ (Case II). Then, 
$\srotilde$ is a $\mu$ reasonable cross-section. The sets $\widetilde{\cS}_\vre$ are
$\mu_{\srotilde}$-JM, and the set $\widetilde{\cB}$ is tempered. In Case I, $\wt{\cB}$ is $\mu_{\srotilde}$-JM, and the same 
statement holds in Case II provided the norm on $\bR^d$ is either the sup-norm or an $\bar A_{\vec{\al}}$-analytic norm (in particular, if it is the Euclidean norm). In Case I
the measure $\mu_{\srotilde}$ is $K_f$-invariant (where $K_f$ is as in
\eqref{eq: def K_f}), and in Case II, $\mu_{\srotilde}$ is $M_{\vec{\al}}$-invariant
(where $M_{\vec{\al}}$ is as in \S \ref{subsec: adeles case 2}).
\end{theorem}
\begin{proof}
The fact that $\srotilde$ is a $\mu$-cross-section, its
$\mu$-reasonability, and the statements about Jordan 
measurability and temperedness are immediate consequences of Proposition  
\ref{prop:941}, the fact that $\pi:\XXnA\to \XX_n $ is a $K_f$-fiber bundle, and the results in
\S\ref{sec: applications} and \S\ref{sec: special subsets} proved
above. 

 For the invariance of $\mu_{\srotilde}$ under $K_f$ in
Case I, we use Theorem 
\ref{cor:csmeasure}, item \eqref{item: one more}, the fact that $K_f$ preserves
$m_{\XXnA}$, and the fact that $K_f$ acts transitively on the fibers
of $\pi$ and hence leaves $\srotilde$ invariant (see the discussion of
$\pi$ in \S \ref{subsec: adeles case 
  1}). For the invariance of
$\mu_{\srotilde}$ under $M_{\vec{\al}}$ in Case II, we use the same argument,
noting that $\mu$ is 
$M_{\vec{\al}}$-invariant and $M_{\vec{\al}}$ preserves the fibers of $\pi$. 
  \end{proof}

  Given $\theta \in \R^d$ we define $\Lam_\theta = u(-\theta)\Z^n$ as
in \eqref{eq: def lattice to a vector} and \eqref{eq: def u(v)}, and 
define its lift to the adelic space by
\index{L@$\widetilde{\Lam}_\theta$ -- adelic lift of $\Lam_\theta$}
\begin{equation}\label{eq: adelic lift lattice}
  \widetilde{ \Lam}_\theta \df
  (u(-\theta), e_f) \SL_n(\Q) \in \XXnA.
  \end{equation}

\begin{proposition}\label{cor: adelic genericity alpha}
Let $\vec{\alpha}$ be as in Case II, and let $\mu =
m_{\tilde{L}_{\vec{\al}} \tilde{y}_{\vec{\al}}}$ (see Proposition~\ref{prop: lifted measure case II} for notation). 
Suppose the norm on $\bR^d$ is either the sup norm or $\bar{A}_{\vec{\al}}$-analytic (see Definition~\ref{def: a bar analytic}).
Let $\vre \in (0,r_0)$ be large enough so that $\mu_{\sro}(\cS_\vre)>0$, 
and let $\mu'$ be equal to any one of 
$\mu_{\srotilde}, \, \mu_{\srotilde}|_{\wt{\cS}_\vre}, \,
\mu_{\srotilde}|_{\wt{\cB}}$. Then the point $\wt{\Lam}_{\vec{\alpha}}$
is $(a_t, \mu')$-generic (see Definition~\ref{def: genericity}). Moreover, if $t_k\to\infty$ are such that  
$a_{t_k}\wt{\Lam}_{\vec{\al}}\in \srotilde$ and $\lim_k a_{t_k}\wt{\Lam}_{\vec{\al}} = x$, then $x \in\on{supp}\mu_{\srotilde}$.
\end{proposition}
\begin{proof}
  By Proposition \ref{prop: lifted measure 
  case II}, any point in $\tilde{L}_{\vec{\al}} \tilde{y}_{\vec{\al}}$ is
$(a_t, \mu)$-generic, and thus, 
  by Propositions \ref{prop: horospherical and generic} and \ref{prop: explicit compact push}, $\wt{\Lam}_{\vec{\alpha}}$
is $(a_t, \mu)$-generic. 

The statement will follow by applying
Theorem~\ref{thm:Sgenericity}. We verify its conditions.   
First observe that
$\wt{\Lam}_{\vec{\al}}\notin\Del_{\srotilde}^\bR$. To see this note
that if  $\wt{\Lam}_{\vec{\al}}\in \Del_{\srotilde}^\bR$, the orbit
$\set{a_t\wt{\Lam}_{\vec{\al}}: t>0}$ intersects
$\wt{\cS}_{r_0,<\vre'}$ for arbitrarily small $\vre'>0$. This leads to
a contradiction using a similar argument to the one giving the  
proof of Proposition~\ref{prop:description of measures}\eqref{item: 1.3}.

Next, we apply Theorem~\ref{thm:Sgenericity} to $\cS'$ being one of
the three sets under consideration, $\srotilde, \wt{\cS}_\vre,
\wt{\cB}$. We need to check that $\cS'$ is $\mu_{\srotilde}$-JM and
that $\mu_{\srotilde}(\cS')>0$. 
The Jordan measurability of $\cS'$ follows from Lemmas~\ref{lem: cB
  null 2},~\ref{lem: cB null 3},~\ref{lem: S eps null 1}.  The
positivity of 
$\mu_{\srotilde}(\wt{\cS}_\vre)$ holds by choice of $\vre$ and the
positivity of $\mu_{\srotilde}(\wt{\cB})$ follows from Lemma~\ref{lem: positivity of B case 2} and the fact that 
$\mu_{\srotilde}(\wt{\cB}) = \mu_{\sro}(\cB) $ (by Proposition~\ref{prop:941}\eqref{item: 2 lift} and Proposition~\ref{prop: lifted measure case II}).

We prove the last assertion in the statement. Let $x$ be an accumulation point in positive time of 
$a_t\wt{\Lam}_{\vec{\al}}$. Note that by~\eqref{eq: adelic lift lattice}, \eqref{eq:stable relation}, \eqref{eq: def bar A}, and 
\eqref{eq:defining lambda tilde} we have  that $\wt{\Lam}_{\vec{\al}} = q\wt{\Lam}$ where 
$q$ satisfies $a_tqa_{-t}\to e_\infty$. Thus, the accumulation points of 
$a_t\wt{\Lam}_{\vec{\al}}$ in positive time are equal to those of $ a_t\wt{\Lam}$. By 
Proposition~\ref{prop: lifted measure case II}, $a_t$ acts uniquely ergodically on the orbit $\tilde{L}_{\vec{\al}} \tilde{y}_{\vec{\al}}$, and hence $x\in \tilde{L}_{\vec{\al}} \tilde{y}_{\vec{\al}}= \on{supp}\mu$. It is clear that $\on{supp} \mu_{\srotilde} = \on{supp}\mu \cap \srotilde$ and so $x\in \on{supp}\mu_{\srotilde}$ as claimed.
 
  \end{proof}

\subsection{The special case $d=1$}
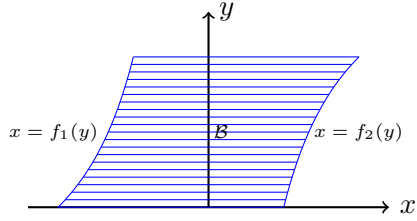
\begin{figure}[htbp]
\center{

\begin{tikzpicture}[ scale=2]
\draw[black, thick, ->] (-1.2,0)--(1.2,0) node[anchor = west]{$x$};

\draw[black, thick,->] (0,0)--(0,1.3)node[anchor = west]{$y$};

\draw[scale=1, domain=0:1, smooth, variable=\y, blue] plot ({-(1/(1+\y))}, {\y});
\draw[scale=1, domain=0:1, smooth, variable=\y, blue] plot ({1/(2-\y)}, {\y});

\foreach \y in {0,0.05,...,1.05}
{ 
\draw[blue, thin, opacity = 0.4] ({-(1/(1+\y))}, {\y}) -- ({1/(2-\y)}, {\y});
}

\draw ({-(1/(1+0.5))}, {0.5}) node[anchor = east]{\!\tiny{$ x = f_1(y)$}};
\draw ({1/(2-0.5)}, {0.5}) node[anchor = west]{\!\tiny{$ x = f_2(y)$}};
\draw (0, {0.5}) node[anchor = west]{\!\tiny{$\cB$}};

\end{tikzpicture}
}
\caption{The set $\{(x,y): y\in [0,1),\; f_1(y) \leq x \leq f_2(y) \}$ parameterizing $\cB$. }
\label{fig: for uri 2}
\end{figure}

When $d=1$ we can fully describe the cross-section $\sro$ and the sets
$\cS_\vre$ and $\cB$ (see Figure \ref{fig: for uri 2}). 

\begin{proposition}\label{prop: explicit d=1}
  Let $d=1$ and let
  $$
u_x \df \left( \begin{matrix}  1 & x \\ 0 & 1 \end{matrix} \right), \
\ h_y \df \left( \begin{matrix} 1 & 0 \\ y & 1 \end{matrix} \right), \
\ \Lam_{x,y} \df u_x h_y \Z^2.$$
Then we can choose $r_0 = 1$, and for $\vre \in (0, 1)$, 
\begin{equation}\label{eq: first assertion Seps}
\cS_\vre = \left\{\Lam_{x,y}: |x|\leq \vre, \, y \in [0,1)
\right\}.
\end{equation}
Setting
$$
 f_1(t) \df - \frac{1}{1+t} \ \ \ \text{ and
} \ \ f_2(t) \df \frac{1}{2-t}, 
$$
we have
\begin{equation}\label{eq: explicit Best}
\cB =\left \{\Lam_{x,y} \in \srosharp : y\in [0,1), \, f_1(y) \leq x \leq
f_2(y)\right\} , 
\end{equation}
and the map $(x,y) \mapsto \Lam_{x,y}$ is injective on
$\left\{(x,y) : y \in (0,1), \, |x| < 1 \right\}.$
  \end{proposition}
  \begin{proof}
In case $d=1$ the choice $r_0=1$ satisfies \eqref{eq: ro satisfies}, and the
group $\{h_y: y \in \R\}$ coincides with the group $H$ in \eqref{eq:
  def subgroups}. The orbit $H\Z^2$ is a periodic orbit
consisting of all lattices that contain $\mathbf{e}_2$ as a primitive
vector, and $\{h_y: y \in [0,1) \}$ is a fundamental domain for the
quotient $H/H(\Z) \cong H\Z^2$. Clearly $(x,1)^{\on{t}} \in \Lam_{x,y} \cap
D_1$ when $|x| \leq \vre$, and conversely, for every $\Lam \in \cS_\vre$ there is
$x \in [-\vre, \vre]$ such that $(x,1)^{\on{t}} \in \Lam_{\mathrm{prim}}$,
and hence $u_{-x} \in H\Z^2$. This proves \eqref{eq: first assertion
  Seps}.

To see the injectivity, note that 
\begin{equation} \label{eq: formula Lam x y}
\Lam_{x,y} = \left(\begin{matrix} 1+ xy & x \\ y &
    1 \end{matrix} \right) \Z^2 = \left\{\left(\begin{matrix}m(1+xy)
      + nx \\ my+n\end{matrix} \right) : m,n \in \Z\right\}.
\end{equation}
Thus,
when $\Lam_{x,y} \in \srosharp$, the vector $(x,1)^{\on{t}}$  is
the unique vector in $\Lam_{x,y} \cap D_1$. That is, for $\Lam_{x,y}
\in \srosharp$, $x$ is uniquely determined, and since $u_{-x} \Lam \in
H\Z^2$, it follows that $y \in [0,1)$ is also uniquely
determined. When $\Lam_{x,y} \notin 
\srosharp$ there must be a horizontal  vector in $\Lam_{x,y}$
of length less than $2$, that is,
integers $(m,n)$ in \eqref{eq: formula Lam x y} satisfying
$$my+n
=0 \ \text{  and } \ |m(1+xy )+ nx |<2.$$
Plugging the first of these equations into the second implies  $m \in \{-1,
0, 1\}$. If $y \in (0,1)$ and $(m,n) \neq (0,0)$, $my+n=0$ is now
impossible. This proves the
injectivity. We note in passing that 
the map we have defined is essentially the map
$\psi$ of \eqref{eq: psi coordinates}.

We now determine the set of $(x, y)\in (-1,1) \times [0,1)
$ for which $\Lam_{x, y} \in \cB$. Using
\eqref{eq: formula Lam x y} and the definition of $\cB$ we see that  $\Lam_{x,y}
\notin \cB$ if and only if $\Lam_{x,y} \notin \srosharp$ or there is
$(m,n) \notin \{ (0,0), (0,1)\}$ such that  
\begin{equation} \label{eq: first condition} 
my +
n \in [0,1)
\end{equation} and 
\begin{equation}\label{eq: second condition}
|m(1+xy)+ nx |\leq |x|.
\end{equation}
 The condition \eqref{eq: first
  condition} implies that 
$ n = -\lfloor my \rfloor$.
Suppose first that $x \in
[0,1)$. In this case, \eqref{eq: second condition} becomes 
$|x(my+n) +m| \leq x$, and 
by plugging \eqref{eq: first condition} into
\eqref{eq: second condition} and examining the possibilities for $m,n$ one
sees that the only possibility is $m=- 1, n= 1$.  This gives the
inequality 
$1+xy -  x \leq x
$, which is equivalent to $x \geq f_2(y)$. In 
the second case $x \in (-1,0]$ the only possibilities become $m= 
1, n=0$, leading to  
$1+ xy \leq  -x$ or $x \leq -f_1(y)$. 
\end{proof}

  \begin{remark}
For $d>1$, our description of $\sro$ is not as explicit since it depends
on the choice of an explicit fundamental domain for $H\Z^n$ in
$H$. It would be interesting to completely describe the set $\cB$, for
$d=2$ and some
fixed norm on $\R^2$. 
    \end{remark}

\section{Interpreting the visits to the
  cross-section}\label{sec: interpreting}
Given $\theta \in \R^d$, let $\wt{\Lam}_\theta$ be as in \eqref{eq:
  adelic lift lattice}. 
The goal of this section is to read off
Diophantine properties of $\theta$, from the successive times $t_i$
for which $a_{t_i} {\widetilde \Lam}_\theta \in \widetilde{ \sro}$. More precisely,
for $\widetilde{\cB}, \, \widetilde{\cS}_\vre$ as in \eqref{eq: lifted cross section}, 
we will relate the successive visits to $\widetilde{\cB}$ to best approximations,
and the successive visits to $\widetilde{\cS}_\vre$, to $\vre$-approximations. 

\subsection{The adelic cross-section as a Cartesian product}
Let
$\psi:\sro^{\sharp}\to  \crly{E}_n\times \bar B_{r_0}$ be as in
\eqref{eq: psi 
  coordinates}, and 
$\srotilde^{\sharp} \df
\pi^{-1}(\sro^{\sharp}).$
We now augment $\psi$ and define a map \index{P@$\widetilde{\psi}$ --
  the adelic lift of $\psi$}
\begin{align}\label{eq: augmented map} 
\widetilde{\psi} 
: \srotilde^{\sharp} \to 
  \crly{E}_n\times \bar B_{r_0}\times \widehat{\bZ}^n \ \ \text{ by } \ \ \widetilde{\psi} =
  (\psi\circ \pi,\psi_f); 
\end{align}
that is, the first two coordinates of $\widetilde{\psi}$ are given by
$\psi \circ \pi$, and the third coordinate  
is given by a map $\psi_f:\srotilde^{\sharp}\to
\widehat{\bZ}^n$, which  
we now define.

Given $\widehat{\Lam}\in \srotilde^{\sharp}$ we may write $\widehat{\Lam} =
(g_\infty, g_f)\Gaa$ with $g_f\in K_f$, and then  
$\Lam \df \pi(\widehat{\Lam}) = g_\infty\Z^n$. Since $\Lam\in
\sro^{\sharp},$  the vector $v(\Lam)$ defined in \eqref{eq: vLam} is the
unique primitive vector of the lattice $\Lam$ in $D_{r_0}$, and the
columns of $g_\infty$ 
are a basis of $\Lam$. Replacing $g_\infty$ by another
representative of the coset $g_\infty \Ga$, we may assume that the
$n$th column of $g_\infty$ is $v(\Lam)$. Moreover, the
uniqueness of $v(\Lam)$ implies that if $\ga\in \Ga$ satisfies that
$g_\infty\ga$ is another 
representative of $\Lam$ having this property, then

$$g_\infty \mathbf{e}_n = v(\Lam) = g_\infty \gamma \mathbf{e}_n,$$
and hence  
$\ga\in H(\bZ)$. With these choices
we define  \index{P@$\psi_f$ -- the finite coordinate of $\widetilde{\psi}$}
\begin{equation}\label{eq: choices define}
  \psi_f(\widehat{\Lam}) \df g_f\mb{e}_n;
  \end{equation}
that is, if we write $g_f= (g_p)_{p\in \mathbf{P}}$, then
$\psi_f(\widehat{\Lam}) $ is 
the element of $\widehat{\bZ}^n$ whose $p$ coordinate is the $n$-th
column of $g_p$. The above discussion implies that $\psi_f$ is
well-defined. 

\begin{lemma}\label{lem: psi tilde continuous}
The map $\widetilde{\psi}$ is $K_f$-equivariant and continuous.
\end{lemma} 
\begin{proof}
The fact that $\wt{\psi}$ commutes the $K_f$-actions on $\srotilde^{\sharp}$ and on $\hat{\bZ}^n_{\on{prim}}$ follows 
directly from the procedure defining $\wt{\psi}$ discussed above.

We prove its continuity.
We have already seen in Lemma \ref{lem:sro-is-lcsc} that $\psi\circ \pi$ is
continuous, and it remains to establish the continuity of
$\psi_f$. 
Assume $\widehat{\Lam}_k\to \widehat{\Lam}$ in $\srotilde^{\sharp}$ and write
$\widehat{\Lam}_k=\left(g_\infty^{(k)}, g_f^{(k)}\right)\Gaa, \ \widehat{\Lam} =
(g_\infty,g_f)\Gaa$, where $\Gaa \df \SL_n(\Q)$
and where the
representatives are chosen so that
$g_f,g_f^{(k)}\in K_f$ and $g_\infty\mb{e}_n,
g_\infty^{(k)}\mb{e}_n\in D_{r_0}$. There are $\ga_k\in \Gaa$ such that  
$$\left(g_\infty^{(k)}\ga_k, g_f^{(k)}\ga_k \right)\to
(g_\infty,g_f),$$
and we need to show that
$g_f^{(k)}\mb{e}_n\to g_f\mb{e}_n$. It follows from our choice of
representatives that
$$\ga_k\in \Gaa\cap 
K_f = \SL_n(\Z).$$

We claim that $\ga_k \mb{e}_n = \mb{e_n}$ for all large enough $k$;
 this will imply $g_f^{(k)}\mb{e}_n = g_f^{(k)}\ga_k\mb{e}_n \to
g_f\mb{e}_n$ and conclude the proof. To prove the claim, suppose by
contradiction that 
$\ga_k\mb{e}_n\ne \mb{e}_n$ for infinitely many $k$, then along a
subsequence, the lattices   
$\Lam_k  = \pi(\widehat{\Lam}_k)$ contain the two distinct primitive vectors $v_k
= g_\infty^{(k)}\mb{e}_n, \, w_k = g_\infty^{(k)}\ga_k \mb{e}_n$, and the
sequence $\Lam_k$ converges in $\XX_n$.  
Note that by assumption $v_k, w_k \in D_{r_0}$. By passing to a further subsequence
we may assume that the sequences $(v_k)$ and $(w_k)$ converge to
limits in $D_{r_0}$. Since the sequence of lattices
$(\Lam_k)$ is bounded in $\XX_n$, and since $v_k, w_k$ are 
distinct, the limits $\lim v_k, \lim w_k$ are distinct vectors  
in $D_{r_0}$ which belong to the limit lattice $\pi(\widehat{\Lam})
= g_\infty\bZ^n$. This contradicts the assumption that $\widehat{\Lam}\in
\srotilde^{\sharp}$.
\end{proof}

We now describe the image of the cross-section measure under
$\wt{\psi}$. We first set up some notation. 
Denote the natural projections by  \index{P@$P_1, P_2, P_3, P_{12}$ --
  projections to factors of $\crly{E}_n \times \R^d \times \widehat{\Z}^n_{\prim}$}
  \begin{equation}\label{eq: definition of lots of projections}
    \begin{split}
      P_{1}: & \ \crly{E}_n \times \R^d
      \times \widehat{\Z}^n_{\prim} \to \crly{E}_n \\
      P_{2}: &  \ \crly{E}_n \times \R^d
      \times \widehat{\Z}^n_{\prim} \to \R^d\\
      P_{3}: & \ \crly{E}_n \times \R^d
      \times \widehat{\Z}^n_{\prim} \to \widehat{\Z}^n_{\prim}\\
P_{12}: & \  \crly{E}_n \times \R^d
\times \widehat{\Z}^n_{\prim} \to \crly{E}_n \times \R^d. 
\end{split}
  \end{equation}
Also let $\mathrm{Proj}: \R^d \to \bS^{d-1}$ be the radial projection
as in \eqref{eq: def Proj}. Given an $\set{a_t}$-invariant measure $\mu$ on $\XXnA$, let
$\mu_{\srotilde}$ 
be the cross-section measure,  
and define a measure on $\crly{E}_n \times \R^d
\times \widehat{\Z}^n_{\prim}$ by 
\begin{equation}\label{eq: definition of projection nu}
  \nu \df \widetilde{\psi}_* \mu_{\srotilde}.
  \end{equation}
Let $\nu^{(\crly{E}_n)}, \, \nu^{(\R^d)}, \nu^{(f)}, \nu^{(\infty)},
\, \nu^{(\bS^{d-1})}, \nu^{(\crly{X}_d)}$
denote the projection of $\nu$ under the maps $P_1, P_2, P_3, P_{12},
\mathrm{Proj} \circ P_2, \pi_{\crly{X}_d}\circ P_1$
respectively.  Let
$\bar A_{\vec{\al}}^{(1)}, M_{\vec{\al}}$ be the groups  in \eqref{eq: 
    def bar A 1} and \S \ref{subsec: adeles case 2} respectively, and
  note that $\bar A_{\vec{\al}}^{(1)}$  and acts linearly on $\R^d$ and on $\crly{E}_n$ via
  the embedding $\bar A_{\vec{\al}}^{(1)}\subset H$. 
  
We have:
\begin{proposition}\label{prop: we have 10.2}
In Case I, with $\mu = m_{\XXnA}$, we have 
\begin{equation}\label{eq: conclusion prop case I}
  \nu = \frac{1}{\zeta(n)} \left(
  m_{\crly{E}_n} \times m_{\R^d}|_{\bar B_{r_0}}\times 
  m_{\widehat{\bZ}^n_{\mathrm{prim}}} \right);
\end{equation}
in particular, the measures $\nu^{(\crly{E}_n)}, \, \nu^{(\R^d)}, \,
\nu^{(f)}, \, \nu^{(\XX_d)}$ are 
scalar multiples of the measures $m_{\crly{E}_n}, \, m_{\R^d}|_{B_{r_0}}, \,
m_{\widehat{\Z}^n_{\prim}}, \, m_{\XX_d}$, and the
measures $\nu^{(\R^d)}$ and $\nu^{(\bS^{d-1})}$ are invariant under
any linear transformations of $\R^d$ preserving the norm $\| \cdot
\|$. 

In Case II, with $\mu = m_{\tilde{L}_{\vec{\al}} \tilde{y}_{\vec{\al}}}$, we have $\nu =
\nu^{(\infty)} \times \nu^{(f)}$, and there exist finitely many subsets $\cO_i\subset \bar{A}_1$ homeomorphic to closed balls, 
lattices $\bar{\Lam}_i\in \bar{A}_{\vec{\al}}y_{\vec{\al}}\cap \crly{E}_n$, and vectors $v_i\in \bar{\Lam}_i$, $i=1,\dots,k$, such that 
the following hold:
\begin{enumerate}
\item \label{item: projection 1}
The measure $\nu^{(\crly{E}_n)}$ is the sum of the push-forwards of the restriction to each $\cO_i$ of 
a Haar measure on $\bar{A}_1$, via the orbit map $\bar{a}\mapsto \bar{a}\bar{\Lam}_i$. In particular, $\nu^{(\crly{E}_n)}$ is supported on a $(d-1)$-dimensional submanifold of $\crly{E}_n$.

\item\label{item: projection 2} 
The measure $\nu^{(\R^d)}$ is the sum of the pushforwards of the restriction to each $\cO_i$ of 
a Haar measure on $\bar{A}_1$, via the orbit map $\bar{a}\mapsto \bar{a}v_i$.
  In particular, $\nu^{(\bR^d)}$ is supported on a $(d-1)$-dimensional submanifold of $\bR^d$.
  \item \label{item: projection 3}
  The measure $\nu^{(f)}$ is the unique $M_{\vec{\al}}$-invariant probability measure supported on the orbit $M_{\vec{\al}}\mb{e}_n \subset \hat{\bZ}_{\on{prim}}$.
\item \label{item: projection 4}
  The measure
  $\nu^{(\infty)} $  is the sum of the pushforwards of the restriction to each $\cO_i$ of 
a Haar measure on $\bar{A}_1$, via the orbit map $\bar{a}\mapsto \bar{a}(\bar{\Lam}_i,v_i)$. In particular, it is supported on
a $(d-1)$ dimensional submanifold and is singular with respect to $\nu^{(\crly{E}_n)} \times 
  \nu^{(\R^d)}$.
\item \label{item: projection 5}
  The measure $\nu^{(\bS^{d-1})}$ is supported on a proper subset of
  $\bS^{d-1}$, and in particular, is not invariant under the group of
  orthogonal transformations of $\R^d$.
  \end{enumerate}

\end{proposition}
\begin{proof}
  We have the following commutative diagram:
  \begin{equation}\label{eq:diagram1160}
\xymatrix{
 \XXnA\ar[r]^{\widetilde{\psi}} \supset \widetilde{\cS}_{r_0}^\sharp \ar@<-3ex>[d]^{\pi}
 & \crly{E}_n \times \R^d \times \widehat{\Z}^n_{\prim} \ar[d]^{P_{12}}
 \\
 \XX_n 
\supset \srosharp
\ar[r]_{\psi} 
& \crly{E}_n \times \R^d 
}
\end{equation}
Recall from Lemma~\ref{lem: psi tilde continuous}, that $\wt{\psi}$ is $K_f$-equivariant. The map $\psi$ has the following 
weak equivariance property: Let $H_0\df \set{\smallmat{g&0\\0&1}:g\in \SL_d(\bR)} \subset H$. Note that $\bar{A}_1\subset H_0$
and that $H_0$ acts diagonally on $\crly{E}_n\times \bR^d$. If $\Lam\in \srosharp$ and $h\in H_0$ is such that $h\Lam\in \srosharp$ and $h v_\Lam = v_{h\Lam}$, then $\psi(h\Lam) = h\psi(\Lam)$.

Suppose we are in Case I. 
By Theorem \ref{thm: MKf}, the measure
$\mu_{\srotilde}$ is $K_f$-invariant.  By the $K_f$-equivariance of $\wt{\psi}$, $\nu$ is $K_f$-invariant as well. Since 
$K_f$ acts transitively on $\hat{\bZ}_{\on{prim}}$ and $m_{\hat{\bZ}_{\on{prim}}}$ is the unique $K_f$-invariant probability measure there, we see that $\nu$ must be the product of $\nu^{(\infty)}$ and $m_{ \hat{\bZ}_{\on{prim}}}$.
Also, by Proposition
\ref{prop:description of measures1} and the commutativity of the diagram, $(P_{12})_*\nu =
\frac{1}{\zeta(n)} \left(m_{\crly{E}_n} \times m_{\R^d}|_{\bar B_{r_0}} \right)$,
  and \eqref{eq: conclusion prop case I} follows. The fact that $\nu^{(\XX_d)}$ is proportional to $m_{\crly{X}_d}$ 
  is immediate
  from \eqref{eq: measures pushed}. Finally, the additional assertion about
  the invariance properties of $\nu^{(\R^d)}$ and $\nu^{(\bS^{d-1})}$  follows from the fact that a linear
  transformation preserving $\| \cdot\|$ also preserves $\sro$ (which
  can be readily seen from \eqref{eq: def Dro} and \eqref{eq:main
    sets}), Theorem \ref{cor:csmeasure}\eqref{item: one more}, and the weak equivariance property of $\psi$ discussed above.

  Now suppose we are in Case II. By Theorem \ref{cor:csmeasure}\eqref{item: one more} the group $M_{\vec{\al}}$ preserves
  $\mu_{\srotilde}$. It also preserves 
  the fibers of $P_{12}$ and using Proposition~\ref{prop: lifted measure case II} and the fact that $\srotilde$ is $M_{\vec{\al}}$-invariant we conclude that it acts transitively on the
  intersection of  each fiber with $\supp \mu_{\srotilde}$. 
Therefore, $\nu = \nu^{(\infty)} \times \nu^{(f)}$. Assertion \eqref{item: projection  
    3} now follows from equations \eqref{eq: choices define} and \eqref{eq:defining lambda tilde}. 
 Assertion \eqref{item: projection 4}, and consequently \eqref{item: projection 1} and 
  \eqref{item: projection 2}, follows from
  Proposition \ref{prop:description of measures}
  and the weak equivariance property of $\psi$ discussed above.
  Assertion \eqref{item:
    projection 5} follows from Proposition \ref{prop: compact A
    orbit};
  indeed, if $\supp \,
  \nu^{(\bS^{d-1})}$ contained an eigendirection for the group $\bar
  A$ then some  $\Lam \in A \Lam_0$ would contain 
a vector whose horizontal component  was in the
 direction of one of the axes.
\end{proof}

\subsection{Hitting the subsets $\cB$ and $\cS_\vre$}
\begin{definition}
We say that two infinite sequences $(a_k)_{k=1}^\infty, \,
(b_k)_{k=1}^\infty$ are {\em prefix-equivalent} if there 
exist $k_0,\ell_0$ such that $a_{k_0+i} = b_{\ell_0+i}$ for all $i\ge
0$. By convention, any two finite sequences are prefix-equivalent, and
a finite sequence is not prefix-equivalent to an infinite sequence. 
\end{definition}
For $\theta\in \bR^d,$ let 
$$\cY_\theta \df \set{t\ge 0 : a_t\Lam_\theta \in \sro} \ \ \ \text{
  and } \ \ \ \cY_\theta^\sharp \df \set{t\in \cY_\theta : a_t\Lam_\theta\in
  \srosharp}.$$
If $t \in \cY_\theta \smallsetminus \cY_\theta^\sharp$ then
$\Lam_\theta$ contains two vectors with the same vertical component
$e^t$, and with horizontal components differing by a vector of size
$O\left(e^{-t/d} \right)$. Thus, by discreteness of $\Lam_\theta$, $\cY_\theta
\smallsetminus \cY_\theta^\sharp$ is finite, and thus $\cY_\theta$ and
$\cY_\theta^\sharp$ are prefix-equivalent. Observe that in the
notation of \eqref{eq: vLam}, for 
each $t\in \cY_\theta^\sharp$, there is a unique primitive vector 
$v(a_t\Lam_\theta)\in (a_t\Lam_\theta)_{\on{prim}} \cap D_{r_0}$. In particular, 
$$u(\theta)a_{-t} v(a_t\Lam_\theta) \in \bZ^n_{\on{prim}}.$$
The next proposition is key to 
our analysis: 
\begin{proposition}\label{prop: hitting the subsets}
Let $\| \cdot \|$ be a norm on $\R^d$, $\vre>0,$ and $ \theta \in
\R^d$. 
\begin{enumerate}
\item Let $\set{t_k}$ be the ordering of $\set{ t\in \cY_\theta^\sharp
    : a_t\Lam_{\theta}\in \cB}$ as an increasing sequence. Then the sequence
  $\mb{v}_k = (\mb{p}_k,q_k)  \df u(\theta)a_{-t_k}
  v(a_t\Lam_\theta)\in \bZ^n_{\prim}$ is prefix-equivalent to the
  sequence of best-approximations of $\theta$. 
  \item Let $\set{t_k}$ be the ordering of $\set{ t\in
      \cY_\theta^\sharp : a_t\Lam_{\theta}\in \cS_\vre}$ as an increasing
    sequence. Then the sequence $\mb{w}_k = (\mb{p}_k,q_k)  \df
    u(\theta)a_{-t_k} v(a_t\Lam_\theta)\in \bZ^n_{\prim}$ is
    prefix-equivalent to the sequence of $\vre$-approximations of
    $\theta$.  
\item For $t\in \cY_\theta^\sharp$, if we let $\mb{u}  =  u(\theta) a_{-t}v(a_t\Lam_\theta)\in \bZ^n_{\on{prim}},$
then (with the notation \eqref{eq: def rho En}, 
  \eqref{eq: def displacement vector}), 
\begin{equation}\label{eq: item iii equivalence}
\widetilde \psi \left(a_t \widetilde \Lam_\theta \right) =
  \left(\rho_{\crly{E}_n}(\bu), \disp(\theta, \bu), 
    \bu\right). 
\end{equation}

\end{enumerate}

%
  \end{proposition}
  \begin{proof}
    Since $\cY_\theta \sm
    \cY_\theta^\sharp$ is finite, we can restrict attention to $t \in
    \cY_\theta^\sharp$.  We first treat 
  $\vre$-approximations. Note that for $\mb{u} = \smallmat{\mb{p}\\q} \in
  \Z^n$, 
  $$\mb{u} \textrm{ is an $\vre$-approximation of $\theta$} \ \
  \Longleftrightarrow \ \ 
  q^{1/d} \norm{\mb{p} - q\theta}\le \vre \ \text{ and } \mb{u} \in \Z^n_{\prim}.$$
   This in turn implies that for $t = \frac{1}{d}\log q$, the lattice
   $a_t\Lam_\theta$ satisfies  
 \begin{equation}\label{eq: lattice satisfies} \begin{split}
(a_t \Lam_\theta)_{\on{prim}}\ni & \diag{q^{1/d},\dots,
                                   q^{1/d},q^{-1}}u(-\theta) \mat{\mb{p}\\q}  \\
   = &  \mat{q^{1/d}(\mb{p}
     -q\theta)\\1} = \mat{\disp (\theta, \mb{u})\\1}\in D_{\vre}.
   \end{split}
  \end{equation}
 Reversing this computation shows that if $t\ge0$ is such that
 $a_t\Lam_\theta\in \cS_\vre$, i.e. 
  $(a_t\Lam_\theta)_{\on{prim}} \cap D_{\vre}$ is not empty, then
  there exist $q\in \bN$ and  
  $\mb{p}\in \bZ^d$ such that $\mb{u} = (\mb{p},q)^{\on{t}}\in
  \bZ^n_{\on{prim}}$ and for $t = \frac{1}{d} \log q$ we have
  $a_tu(-\theta)\smallmat{\mb{p}\\ q}\in 
  D_{\vre}$.  In other words, $\mb{u}$ is an $\vre$-approximation of
  $\theta$. 
  
  The bijection between visits of $a_t\Lam_\theta$ to $\cB$ and best
  approximations described in the statement follow the same
  computation coupled with the following observation: a primitive
  vector $\smallmat{\mb{p}\\ q}\in \bZ^n_{\on{prim}}$ with $q\in \bN$
  is a best approximation of $\theta$ if and only if  
  $$C_{\norm{\mb{p}-q\theta}}(q)\cap \Lam_\theta = \set{0, \pm
    u(-\theta)\mat{\mb{p}\\ q}}$$ 
  where
   $$C_{\norm{\mb{p}-q\theta}}(q) = \set{\mat{\mb{x}\\y}\in\bR^n:
     \norm{\mb{x}}\le \norm{\mb{p}-q\theta}; \; \av{y}\le q}.$$ 
 Acting with $a_t$ on $\Lam_{\theta} = u(-\theta)\bZ^n$ with $t =
 \frac{1}{d}\log q$ to bring $u(-\theta)\smallmat{\mb{p}\\q}$ to
 $D_{r_0}$, we see that $\mb{u}$ is a best approximation for $\theta$
 if and only if $a_t\Lam_\theta\in \cB$ with $t$ as above and  
 $v(a_t\Lam_\theta) = a_tu(-\theta)\mb{u}$.

The third assertion is a straightforward computation combining
\eqref{eq: lattice satisfies} with 
\eqref{eq: def displacement vector}, \eqref{eq: def rho En}, \eqref{eq: psi
  coordinates} and \eqref{eq: augmented map}.  
 \end{proof} 

Using Proposition \ref{prop: hitting the subsets}, we now show that for  $\vre > \vre_0$ which appears in \eqref{eq: def epsilon zero} the condition appearing in Proposition \ref{cor: adelic genericity alpha} is satisfied.
\begin{lemma}\label{lem: epsilon zero}
Let $\vec{\al}$ be as in Case II and let $\mu = m_{\tilde{L}_{\vec{\al}} \tilde{y}_{\vec{\al}}}$ (see Proposition~\ref{prop: lifted measure case II} for notation).  Let 
 $\vre_0$ be as in \eqref{eq: def epsilon zero}. Then, for any $\vre>\vre_0$ we have that $\mu_{\srotilde}(\wt{\cS}_\vre)>0$.
\end{lemma}
\begin{proof}
Let $\vre>\vre'>\vre_0$. By
definition of $\vre_0$ and by Proposition~\ref{prop: hitting the subsets}, there exists a sequence $t_k\to\infty$ such that 
$a_{t_k}\wt{\Lam}_{\vec{\al}}\in \wt{\cS}_{\vre'}$. Since the orbit $\set{a_t\wt{\Lam}_{\vec{\al}}}_{t>0}$ is bounded in $\XXnA$
 and since $\wt{\cS}_{\vre'}$ is closed, we may assume without loss of generality that $a_{t_k}\wt{\Lam}_{\vec{\al}}$ converges to some point $x$. By  Proposition~\ref{cor: adelic genericity alpha} we deduce that
$x\in \wt{\cS}_{\vre'} \cap \on{supp}\mu_{\srotilde}$. Furthermore, since $x$ is an accumulation point of a bounded forward orbit
of $a_t$, it follows that the two-sided orbit of $x$ is bounded. Therefore $x\in \srotilde^\sharp$ (since if $\pi(x)$ contains two vectors with the same vertical coordinate, the orbit $a_tx$ diverges in negative time).
By Lemma~\ref{lem:continuity of vector}, 
$\XX_n(D_\vre^\circ) \cap \srosharp$ is open in $\sro$, and hence
its preimage in $\XXnA$ is open in $\srotilde$. 
By choice of $\vre'$ we have  $\cS_{\vre'}\subset \XX_n( D_\vre^\circ)$,  therefore
$$x\in  \pi^{-1}(\XX_n(D_\vre^\circ) \cap \srosharp)\subset \wt{\cS}_{\vre},$$
and we deduce that $\mu_{\srotilde}(\wt{\cS}_\vre)>0$.
\end{proof}

\section{Properties of the cross-section measures}\label{sec:
  properties measures}
In this section we will describe some properties of 
measures on the product space $\crly{E}_n \times \bR^d \times 
\widehat{\bZ}^n$. The measures we will consider arise as follows: in
Case I let $\mu = m_{\XXnA}$, and in Case II let $
\mu = m_{\tilde{L}_{\vec{\al}} \tilde{y}_{\vec{\al}}}$ be as in \S\ref{subsec: adeles case 2}. Let $\mu_{\srotilde}$ be the cross-section
measure, let $\wt{\psi}, \wt{\cB}$ and $\wt{\cS}_\vre$ be as in
\eqref{eq: augmented map} and
\eqref{eq: lifted cross section}, and let 
\begin{equation}\label{eq: def lambda}
\lambda \df
\left\{ \begin{matrix}
\wt{\psi}_*\left(\mu_{\srotilde}|_{\wt{\cS}_\vre} \right) & \text{
  ($\vre$-approximations)} \\
\wt{\psi}_*\left(\mu_{\srotilde}|_{\wt{\cB}}
    \right) & \text{ (best approximations).} \end{matrix} \right. 
\end{equation}
Note that although our notation does not reflect this, both measures
depend on the norm,  and in the case of 
$\vre$-approximations, $\lambda$ depends on 
$\vre$.  For such a measure $\lambda$, 
let $\lambda^{(\crly{E}_n)}, \lambda^{(\R^d)}, \lambda^{(f)},
\lambda^{(\infty)}$ denote the image of $\lambda$ under the
projections $P_1, P_2, P_3, P_{12}$ as in \eqref{eq: definition of
  lots of projections},
and let $\lambda^{(\XX_d)}, \, \lambda^{(\bS^{d-1})}$ denote the image
under the maps $\pi_{\XX_d} \circ P_1$ and $\mathrm{Proj} \circ P_2$
(see \eqref{eq: def pi 
  XXd}, \eqref{eq: def Proj}). 
We will go through the
different cases and give some properties of these measures.

The case of $\vre$-approximations is much simpler. The following two
propositions are immediate from \eqref{eq: smaller r} and Proposition 
\ref{prop: we have 10.2}.  

\begin{proposition}[Case I, $\vre$-approximations]\label{prop:
    description of measures I eps}
Let $\vre \in (0, r_0)$, where $r_0$ is large enough to satisfy
\eqref{eq: ro satisfies}. Let $\mu= m_{\XXnA}$, and let
$\lambda$ be as in \eqref{eq: def lambda} ($\vre$-approximations). 
Let $B_\vre$ denote the 
  ball of radius $\vre$ around the origin in $\R^d$ and let $V_d = m_{\bR^d}(B_1)$. Then 
  $$\lambda^{(\crly{E}_n)} = \frac{\vre^dV_d}{\zeta(n)}\, m_{\crly{E}_n},  \ \
  \ \ \ \lambda^{(\R^d)} =
  \frac{1}{\zeta(n)}\, m_{\R^d}|_{B_\vre}, $$
  $$ \lambda^{\left(f
    \right)} = \frac{\vre^dV_d}{\zeta(n)}\,
  m_{\widehat{\Z}^n_{\prim}}, \ \ \ \ \ \lambda^{(\XX_d)} = \frac{\vre^dV_d}{\zeta(n)}\,
  m_{\XX_d}. $$
  The measures $\lambda^{(\R^d)}, \, \lambda^{(\bS^{d-1})}$ are
  preserved by any linear transformations of $\R^d$ preserving the
  norm $\| \cdot \|$. 
\end{proposition}

\begin{proposition}[Case II, $\vre$-approximations]\label{prop:
    description of measure II eps} 
Let $d \geq 2$ and let $\vre \in (0, r_0)$, where $r_0$ is large enough to satisfy
\eqref{eq: ro satisfies}. Let $\mu = m_{\tilde{L}_{\vec{\al}} \tilde{y}_{\vec{\al}}}$ be as
in \S \ref{subsec: adeles case 2} and let $\lambda$ be as in \eqref{eq:
  def lambda}.  Then $\lambda$ is proportional to $\lambda^{(\infty)} \times
\lambda^{(f)}$, $\lambda^{(\infty)}$ is not proportional to
$\lambda^{(\crly{E}_n)} \times 
\lambda^{(\R^d)}$, and the projections $\lambda^{(\crly{E}_n)}, \,
\lambda^{(\R^d)}$  satisfy statements \eqref{item: projection 1},
     \eqref{item: projection 2}, \eqref{item: projection 3}
     of Proposition \ref{prop: we have 10.2}.
    The measures $\lambda^{(\R^d)}, \, \lambda^{(\bS^{d-1})}$ are not
    invariant under the group of orthogonal transformations.
  \end{proposition}
  We now turn to best approximations. Let \index{B@$\widehat{\cB}$ -- the
    pushforward of $\cB$ to $\crly{E}_n \times \R^d$}
  \index{I@$\mathbf{1}_{\widehat{\cB}}$ -- the indicator of $\widehat{\cB}$}
  $$\widehat{\cB} \df \psi(\cB) \subset \crly{E}_n \times \R^d,$$
  so that
  \begin{equation}\label{eq: if and only if}
    (\Lam, v) \in \widehat{\cB} \ \iff \ \varphi (\Lam, v) \in
    \cB,
  \end{equation}
  and denote the indicator function of $ \widehat{\cB}$ by $\mathbf{1}_{\widehat{\cB}}$. 
  \begin{proposition}\label{prop: splits as a product best}
Let $\lambda $ be the measure as in \eqref{eq: def lambda}
  (best approximations). Then $\lambda = \lambda^{(\infty)} \times
  \lambda^{(f)}$, where in Case I, $\lambda^{(f)}$ is a 
  multiple of $
  m_{\widehat{\Z}^n_{\prim}}$, and in Case II, using the notation of \S\ref{subsec: adeles case 2}, 
   $\lambda^{(f)}$ is a
  multiple of the 
  $M_{\vec{\al}}$-invariant measure on an $M_{\vec{\al}}$-orbit on
  $\widehat{\Z}^n_{\prim}$. The measure $\lambda^{(\infty)}$ is
  absolutely continuous with respect to the measure $\nu^{(\infty)}$
   as in \eqref{eq: definition of projection nu},  and the
  Radon-Nikodym derivative is given by
  \begin{equation}\label{eq: RN derivative}
\frac{d\lambda^{(\infty)}}{d\nu^{(\infty)}}(\Lam, v) =
\mathbf{1}_{\widehat{\cB}}(\Lam, v). 
\end{equation}
    \end{proposition}
    \begin{proof}
      The proof that $\lambda$ can be written as a product
$\lambda^{(\infty)} \times \lambda^{(f)}$, and that $\lambda^{(f)}$ is
either the unique (up to scaling) $K_f$-invariant measure, or an $M_{\vec{\al}}$-invariant measure
on an $M_{\vec{\al}}$-orbit, is identical to the one used
for $\nu$, in the proof of Proposition \ref{prop: we have
  10.2}. Formula \eqref{eq: RN derivative} is clear from \eqref{eq:
  lifted cross section}, 
\eqref{eq:diagram1160} and \eqref{eq:
  def lambda}. 
      \end{proof}
    Formula \eqref{eq: RN derivative} is an explicit formula for
    $\lambda^{(\infty)}$. We will now use it to describe
    $\lambda^{(\crly{E}_n)}, \, \lambda^{(\XX_d)}, \, \lambda^{(\R^d)}, \,
    \lambda^{(\bS^{d-1})}$ in more detail. For this we need to
    consider 
    Case I and Case II separately. We begin with the simpler Case II.

\begin{proposition}\label{prop: description of measures II best}
Let $d \geq 2$, and suppose the norm on
$\R^d$ is either $\bar{A}_{\vec{\al}}$-analytic (see Definition \ref{def: a bar
  analytic}), or is the sup-norm.  
 Let $\mu = m_{\tilde{L}_{\vec{\al}} \tilde{y}_{\vec{\al}}}$ be as
 in Proposition \ref{prop: lifted measure case II},
 let $\lambda$ be given by
\eqref{eq: def lambda} (best approximations), and let
$\lambda^{(\crly{E}_n)}, \lambda^{(\R^d)}, \lambda^{(f)},
\lambda^{(\infty)}, \lambda^{(\XX_d)}$ be as above.

  We have
  \begin{enumerate}
  \item \label{item: measures simple case 1}
  The measures $\lambda^{(\crly{E}_n)}, \, \lambda^{(\R^d)}, \,
  \lambda^{(\XX_d)}$ are singular with respect to the measures
  $m_{\crly{E}_n}, \, m_{\R^d}, \, 
  m_{\XX_d}$.
\item \label{item: measures simple case 2}
  The measure $\lambda^{(\bS^{d-1})}$ is not globally supported. In
  particular $\lambda^{(\R^d)}$ is not invariant under all orthogonal
  transformations.
\item \label{item: measures simple case 3}
  The measure $\lambda^{(\infty)}$ is not a scalar multiple of
  $\lambda^{(\crly{E}_n)} \times \lambda^{(\R^d)}$. 
\end{enumerate}
\end{proposition}
 \begin{proof}
Assertions \eqref{item: measures simple case 1} and \eqref{item:
  measures simple case 2} follow from Propositions \ref{prop: we have 10.2} and
\ref{prop: splits as a product best}. Furthermore, suppose that
$(\Lam, v) \in \supp \, \lambda^{(\infty)}$. In particular
$\Lam_1 \df \varphi(\Lam, v) \in \cB$. Since $\cB$ is open in $\sro$ (see
Lemma \ref{lem:boundary of best}) and $\psi$ is a continuous inverse
of $\varphi$ on $\cB$, there is a 
neighborhood $\cU$ of the identity in $\bar A_{\vec{\al}}^{(1)}$, and a neighborhood
$\cV$ of $(\Lam, v) \in \crly{E}_n \times \R^d$, such that
$$
\supp \, \lambda^{(\infty)} \cap \cV = \{\psi(\bar a \Lam_1 ) :
\bar a \in \cU\}. 
$$
Restricting the
measures $\lambda^{(\infty)}, \lambda^{(\crly{E}_n)},
\lambda^{(\R^d)}$ to the image of $\cU$ under each of the maps 
$$\bar a \mapsto \psi( \bar a \Lam_1) , \ \ \ \ \bar a \mapsto \bar a \Lam,
\ \ \ \ \bar a \mapsto \bar a v, 
$$ 
we see that in these open sets, each of these measures is supported on
a $d-1$ dimensional manifold, and this implies \eqref{item: measures
  simple case 3}. 
    \end{proof}

    \subsection{Further properties of the measures for best
      approximations, Case I}
  \begin{proposition}\label{prop: description of measures I best}
    In Case I, for best approximations, the measures
     $$\lambda^{(\crly{E}_n)}, \,
  \lambda^{(\XX_d)} , \, \lambda^{(\R^d)} , \,
  \lambda^{(\bS^{d-1})}$$ are absolutely continuous with respect to $m_{\crly{E}_n}, \,
  m_{\XX_d} , \, m_{\R^d} , \,
  m_{\bS^{d-1}}$, and we have the following formulae for the Radon-Nikodym
  derivatives:
  \begin{align} \label{eq: RN derivatives Case Ia}
          \frac{d\lambda^{(\crly{E}_n)}}{dm_{\crly{E}_n}} (\Lam) = &  \frac{1}{\zeta(n)}\,
      \int_{B_{r_0}} \mathbf{1}_{\widehat{\cB}}(\Lam, v) \,
      dm_{\R^d}(v) \\
                    \frac{d\lambda^{(\R^d)}}{dm_{\R^d}}
                  (v) =  & \frac{1}{\zeta(n)}\,  \int_{\crly{E}_n}
                  \mathbf{1}_{\widehat{\cB}}(\Lam, v) \,
                  dm_{\crly{E}_n}(\Lam) \label{eq: RN derivatives Case Ib}
  \\
            \frac{d\lambda^{(\crly{X}_d)}}{dm_{\XX_d}} (\Lam') = &
            \frac{1}{\zeta(n)}\,  \int_{B_{r_0}} \, \int_{\bT_{\Lam'}} 
            \mathbf{1}_{\widehat{\cB}}(\Lam(\Lam',f, \mathbf{e}_n), v) \, d
                                                                   m_{\bT_{\Lam'}}(f)
                                                                   \,
                                                                   dm_{\R^d}(v)
                                                                   \label{eq: RN derivatives Case Ic}
  \end{align}
             where $\bT_{\Lam'}$ is as in \eqref{eq: torus of
               functionals} and $\Lam( \cdot) $ is as in                                                                      \eqref{eq:
                                                                       recovering
                                                                       from
                                                                       lift},
                                                                     and,
                                                                     for
                                                                     some
                                                                     $c>0$, 
                                                                     \begin{equation}
                            \frac{d\lambda^{(\bS^{d-1})}}{dm_{\bS^{d-1}}
                        }(\omega) = c\,  \int_0^{r_0} t^{d-1} \, \int_{\crly{E}_n}
                        \mathbf{1}_{\widehat{\cB}}(\Lam, t\omega) \,
                        dm_{\crly{E}_n}(\Lam) \, dt.\label{eq: RN derivatives Case
                          Id}
                             \end{equation}
                           \end{proposition}

                           \begin{proof}
   Equations \eqref{eq: RN derivatives Case Ia},
    \eqref{eq: RN derivatives Case Ib} both 
    follow immediately from \eqref{eq: 
      conclusion prop case I} and \eqref{eq: RN
      derivative}.
    For \eqref{eq: RN derivatives Case Ic}, we give an expression for
    $m_{\XX_d}$ in terms of $m_{\crly{E}_n}$. Let $H$ be the group as in
    \eqref{eq: non-expanding group}, let $G_0$ and $W$ denote
    respectively the subgroup of $H$ obtained by setting $\bx =0$ and
    $A=0$ in \eqref{eq:
      non-expanding group}. Thus $G_0 \cong \SL_d(\R), \, W \cong
    \R^d$, and $H = G_0 \ltimes W$. By Lemma \ref{lem:hic basic}, the
    Haar measure on 
    $H$ can be written as $dm_H(A\bx) = dm_{G_0}(A) dm_W(\bx)$. 
The discussion in the last paragraph of \S \ref{sec: lift functionals} shows that the
matrix $A$ gives the projected lattice $\Lam' = A\Z^d$, and the vector
$\bx$ gives the lift functional $f \in \bT_{\Lam'}$. Therefore, for
fixed $\Lam' \in \XX_d$, 
$$
dm_{\XX_d}(\Lam') = \int_{\bT_{\Lam'}}  dm_{\crly{E}_n}(\Lam(\Lam', f,
\mathbf{e}_n) ) \, dm_{\bT_{\Lam'}}(f),
$$
and using this, in combination with \eqref{eq: RN derivatives Case
  Ia} and the uniqueness of disintegration of measures, we obtain
\eqref{eq: RN derivatives Case Ic}. 
By polar coordinates, there is some $c>0$ such that
$dm_{\R^d}(t\omega) = ct^{d-1}dm_{\bS^{d-1}}(\omega) $, and thus \eqref{eq:
  RN derivatives Case Id} follows from \eqref{eq: RN 
  derivatives Case Ib}. 
                           \end{proof}
                           \begin{proposition}\label{prop: properties
                               of the measures}
In Case I, for best approximations, the measures $$\lambda^{(\crly{E}_n)}, \, \lambda^{(\R^d)},
\, \lambda^{(\XX_d)}, \, \lambda^{(\bS^{d-1})}$$ satisfy:
  \begin{enumerate}[(a)]
  \item\label{item: measures b}
$\lambda^{(\infty)} $ is not a scalar multiple of
$\lambda^{(\crly{E}_n)} \times 
\lambda^{(\R^d)}$.
\item\label{item: measures d}
  The measures
      $\lambda^{(\crly{E}_n)}, \, \lambda^{(\XX_d)}, \, \lambda^{(\bS^{d-1})}$
      have full support, and the support of $\lambda^{(\R^d)}$
      contains a neighborhood of the origin. 
     \item\label{item: measures e}
       For $d>1$, there is $c>0$ such that for any $\Lam \in
       \crly{E}_n$ and any $\Lam' \in \XX_d$, 
       \begin{equation}\label{eq: eqref follows}
         \frac{d\lam^{(\crly{E}_n)}}{dm_{\crly{E}_n}}(\Lam) \leq c \cdot
         \mathrm{sys}(\pi_{\R^d}(\Lam))^d \ \ \text{ and } \ \ 
         \frac{d\lam^{(\XX_d)}}{dm_{\XX_d}}(\Lam') \leq c \cdot 
         \mathrm{sys}(\Lam')^d, 
         \end{equation}
where $\mathrm{sys}(\Lam')$ is the length
of the shortest nonzero vector of $\Lam' \in \XX_d$.
In particular $\lambda^{(\crly{E}_n)} $
    and $\lambda^{(\XX_d)} $ are not scalar multiples of $m_{\crly{E}_n}$ and
    $m_{\XX_d} $. 

  \item \label{item: measures concavity}
  For any $\Lam \in \crly{E}_n$ with no nonzero horizontal vectors, and any
  $\omega \in \bS^{d-1}$, the set  
\begin{equation}\label{eq: actually an interval}
\{t \in \R: \varphi(\Lam,
t\omega) \in \cB\}
\end{equation}
is an interval containing $0$. In particular, for any $\omega \in
\bS^{d-1}$, the function
\begin{equation}\label{eq: the function}
t \mapsto \frac{d\lambda^{(\R^d)}}{dm_{\R^d}}(t \omega)
\end{equation}
is monotone non-increasing, is not a.e.\,an indicator function, and
$\supp \, \lambda^{(\R^d)}$ is 
star-shaped around the origin. 
  \item\label{item: measures f}
      The measures $\lambda^{(\crly{E}_n)}, \lambda^{(\R^d)}, \,
    \lambda^{(\bS^{d-1})}$ are invariant under any linear
    transformation of $\R^d$ preserving the norm $\| \cdot \|$ (where
    in case of $\crly{E}_n$, the action is via the linear action on
    the first $d$ coordinates in $\R^n$). In
    particular, for the Euclidean norm, $\lambda^{(\crly{E}_n)}, \, \lambda^{(\R^d)}$ and $
    \lambda^{(\bS^{d-1})}$ are $\SO_d(\R)$-invariant.
   
    \end{enumerate}

  \end{proposition}

  \begin{proof}
 Assume  first that $d=1$. In this case, the assertions all follow
 easily
 from Proposition
\ref{prop: explicit d=1}. Indeed, in this case there is only one
norm on $\R^d = \R$, the bundle $\crly{E}_n$ is isomorphic to $\R/\Z$
via the map $y \mapsto h_y \Z^2$, and the measures
$\lam^{(\crly{E}_n)}, \, \lam^{(\R^d)}$ are simply the pushforwards to
the $x$ and $y$ axes, of the set in \eqref{eq: explicit Best}.

We now assume that $d > 1$,
and
note that for any $\Lam \in \crly{E}_n$ there is $\vre =
 \vre(\Lam) >0$ such 
 that for all $v \in \R^d$ with $\|v\|<\vre$, we have $\varphi(\Lam,
 v) \in \cB$. Indeed, this follows from the discreteness of $\Lam$ and
 the fact that the cylinders $C_r$ get smaller and smaller as $r \to
 0$. This implies that the integrands in  \eqref{eq: RN derivatives
   Case Ia}, \eqref{eq: RN derivatives Case Ib}, \eqref{eq: RN
   derivatives Case Ic}, \eqref{eq: RN derivatives Case Id} are
 positive on sets of positive measure, and \eqref{item: measures d}
 follows.  

 For \eqref{item: measures e}, let $\Lam \in \crly{E}_n$, let $\Lam'
 = \pi_{\R^d}(\Lam) \in \XX_d$, and let $\bx' \in
 \Lam'$ with $\|\bx'\| = \mathrm{sys}(\Lam')$. Then
 $\Lam$ contains a
 vector $\bx \in \pi_{\R^d}^{-1}(\bx')$ whose horizontal component is 
 $\bx'$ and whose vertical component is in $\left[-\frac{1}{2},
   \frac12\right]$. Let $v \in \R^d$ with $\|v\|> 2 \, 
 \mathrm{sys}(\Lam')$ and let $u(v)$ be as in \eqref{eq: def u(v)}. By
 the triangle inequality, 
 $$\|
 \pi_{\R^d}(u(v)\bx)\|   \leq \|\bx'\|+ \frac{\|v\|}{2} <  \|v\| =
 \|\pi_{\R^d}(u(v) \mathbf{e}_n)\|,$$ 
 and hence $u(v)\Lam \notin \cB$. Using \eqref{eq: RN derivatives
   Case Ia} we obtain 
$$\frac{d\lambda^{(\crly{E}_n)}}{dm_{\crly{E}_n}}(\Lam) \leq
\frac{1}{\zeta(n)}\, m_{\R^d}(B(0, 
2 \, \mathrm{sys}(\Lam'))),$$
and the first inequality in \eqref{eq: eqref follows} follows. For the
second inequality, repeat this argument for every $\Lam$ of the form
$\Lam = \Lam(\Lam', f, \mathbf{e}_n)$ as in \eqref{eq: recovering from
  lift}, and use \eqref{eq: RN derivatives Case Ic}.

For assertion \eqref{item: measures b}, note from Proposition
\ref{prop: we have 10.2} that $\nu^{(\infty)}$ is a  
scalar multiple of
$m_{\crly{E}_n} \times m_{\R^d}|_{B_{r_0}}$.  By \eqref{eq: RN derivative} and the
fact that
 $\lambda^{(\crly{E}_n)}$ and $\lambda^{(\R^d)}$ are of full support, if 
$\lambda^{(\infty)}$ were a scalar 
multiple of $\lambda^{(\crly{E}_n)} \times \lambda^{(\R^d)}$ then its
Radon-Nikodym derivative would be  
constant on $\crly{E}_n \times B_{r_0}$. This contradicts the first inequality in \eqref{eq: eqref follows}.

For \eqref{item: measures concavity}, let $\Lam \in \crly{E}_n$ have
no horizontal vectors and
let $\omega\in \R^d, \ \|\omega\|=1$. Then $\varphi(\Lam,v) \in \cB$ if and only if
for every nonzero $\bx \in \Lam$ with vertical
component $x_n \in (-1,1)$, we have $\|v \|< \| \bx'+ x_n v\|$, where
$\bx' \df \pi_{\R^d}(\bx) \neq 0$. In
other words, the set in 
\eqref{eq: actually an interval}, which we denote by $\bar I(\Lam,
\omega)$, can be written as 
$$
\bigcap_{\substack{\bx \in \Lam \sm \{0\}\\ x_n \in (-1,1)}} I(\bx, \omega), \ \
\mathrm{ where } \ I(\bx, \omega) \df \left\{t \in \R: \|t
\omega \|< \| \bx'+ x_n t\omega\| \right \}.
$$ 
In order to show that $\supp \, \lambda^{\left(\R^d \right)}$ is
star-shaped, it suffices to show that $\bar I(\Lam, \omega)$ is an
interval containing 0 for each $\Lam$, and for this it is enough to
show that $I(\bx, 
\omega)$ is an interval containing 0 for every 
$\bx$. Clearly
$I(\bx, \omega)$ is bounded and contains 0, and if it 
were not an interval, there would be $0< t_1< t_2$ such
that $\|t_i
\omega \|= \| \bx'+ x_n t_i\omega\|$ for $i=1,2$. This implies
\[\begin{split}
& t_2 - t_1 = \|t_2 \omega\| - \|t_1 \omega\| = \| \bx'+ x_n t_2\omega\|
- \| \bx'+ x_n t_1\omega\| \\ = \| & \bx'+ x_n t_2\omega\| - \| \bx'+ x_n t_2\omega +
x_n(t_1-t_2)\omega \| \leq \|x_n(t_1-t_2)\omega \| < t_2 - t_1,
\end{split}\]
a contradiction. The assertion about monotonicity in \eqref{item: measures
  concavity} now follows using \eqref{eq: RN derivatives Case Ib}. If
the function in \eqref{eq: the function} were an indicator function  for some $\omega$,
by \eqref{eq: RN derivatives Case Ib}, the interval
$\bar I(\Lam, \omega)$ would actually be the same for almost
all $\Lam$. To see that this cannot be the case, let $\Lam$ be some
lattice, and let $\bx \in \Lam \sm \{0\}$ so that 
$\sup I(\bx, \omega) = \sup \bar I (\Lam, \omega) $. For almost
every $\Lam$, $\bx$ is unique with this property. Now for almost all nearby
lattices $\Lam' \in \crly{E}_n$, containing a  small perturbation
$\bx'$ of $\bx$, we 
will have $\sup \bar I(\Lam', \omega) = \sup I(\bx', \omega) \neq \sup
I (\bx, \omega) = \sup \bar I(\Lam, \omega)$.

Assertion \eqref{item: measures f} follows immediately from Proposition \ref{prop: we
  have 10.2} and the fact that $\cB$ is invariant under any linear
transformation of the horizontal plane preserving the norm.  
    \end{proof}

\subsection{Cut-and-project structure of approximations to algebraic
  vectors} \label{subsec: c and p structure}
  This section is not needed for the proof of our main results. Its purpose is to highlight a certain structure
  that the set of approximations to algebraic vectors posses.

A {\em cut-and-project set } is a subset $X_0 \subset \R^r$ for which
there are $s \in \N$, a lattice $\Delta \subset \R^{r+s}$, and $W \subset
\R^s$ for which
$$
X_0  = \{ \bx \in \R^r: \exists \by \in W \text{ such that } (\bx,
\by)^{\mathrm{t}} \in \Delta \}.
$$
Here $(\bx, \by)^{\mathrm{t}}$ denotes the vector whose first $r$
entries are those of $\bx$ and whose last $s$ entries are those of
$\by$. In some cases
the vector spaces $\R^r, \R^s, \R^{r+s}$ appearing in the above definition are 
taken to be more general locally compact abelian groups. We will refer
to the more general case of cut-and-project sets arising in this way
as {\em generalized cut-and-project sets}. 
Cut-and-project sets are widely studied in the field of
mathematical quasicrystals, see \cite{BaakeGrimm} for a
comprehensive introduction.

For discrete sets $X' \subset X \subset \R$, we say that $X'$ {\em has full density
  in $X$} if 
$
\frac{\#(X' \cap [0,T])}{\#(X \cap [0,T])} \to_{T \to \infty} 1,
$
and say that discrete $X, Y \subset \R$ are {\em strongly
  asymptotic} if there are subsets $X', Y'$ of full density and a
bijection $\tau: X' \to Y'$ such that
$$|x' - \tau(x')| \longrightarrow_{x' \to \infty} 0\  \ \ (x' \in X').$$ 
Let $K$ be a compact abelian group, let $\bG = K \times \R$, let
$\pi_{\R} : \bG \to \R$ be the projection,  and let $d_{\bG}$ be a
translation invariant metric on $\bG$. For $X \subset \bG$ and $T \in
\R$ we set
$X^{\geq T} \df \{(k, x) \in \bG: x \geq T\}$, and say that $X, Y
\subset \bG$ are {\em strongly asymptotic} if there are subsets $X',
Y'$ of $X$ and $Y$, such that $\pi_{\R}(X'),
\pi_{\R}(Y')$ are of full density in $\pi_{\R}(X)$ and $\pi_{\R}(Y)$,
and there is a bijection $\tau: X' \to Y'$ such that $d_{\bG}(x, \tau(x)) \to 0$
as $\pi_{\R}(x) \to \infty$. 

Our analysis  yields the following:
\begin{proposition}\label{prop: c and P}
Let $\vec{\alpha}$ be as in Case II, let $\| \cdot \|$ be some norm on
$\R^d$, let $\vre>0$, and let $\bu_k = (\bp_k, q_k) \in \Z^n$ be either the sequence
of best approximation of $\vec{\alpha}$, or its sequence of
$\vre$-approximations. Then the sequence $\left(\log(q_k)\right)_{k
  \in \N} \subset 
\R$ is strongly asymptotic to a one-dimensional cut-and-project set,
and the sequence $(\mb{p}_k, \log(q_k)) \subset \widehat{\Z}^n \times \R$,
is strongly asymptotic to a generalized cut-and-project set. 
\end{proposition}
\begin{proof}[Sketch of proof]
  We will not be using this statement in the sequel, and we only sketch the proof. 
  Let $\bT=\R^d/\Delta$ be a torus, where $\Delta$ is a lattice in
  $\R^d$, let $\pi_{\bT}: \R^d \to \bT$ be the projection, let $\bz
  \in \R^d \sm \{0\}$, and let $\alpha_t \curvearrowright \bT$  
  be the corresponding {\em straightline flow}  defined by $\alpha_t
  \pi_{\bT}(\bx) = \pi_{\bT}(\bx + t \bz)$. 
  A cut-and-project set in $\R$ can also be described as the set of visit
  times as in \eqref{eq: visit times}, to a section $\cS$ which is a
  {\em bounded linear section}; i.e., the image under $\pi_{\bT}$ of a
  bounded subset of an affine subspace of
  dimension $d-1$; see \cite[Prop. 2.3]{ASW} for a proof.

  Let $\bar{A}_{\vec{\al}}, y_{\vec{\al}}$ be as in \eqref{eq: def bar A}, let $\bar A_{y_{\vec{\al}}}$ be the
  stabilizer group as in \eqref{eq: def bar A Lam}, let $\bar{\mathfrak{a}} \cong
  \R^d$ be the Lie algebra of $\bar{A}_{\vec{\al}}$, let $\exp : \bar{\mathfrak{a}} \to
  \bar{A}_{\vec{\al}}$ be the exponential map and denote its inverse by
  $\log$. Then the map $\exp$ induces an isomorphism between the
  compact orbit $\bar{A}_{\vec{\al}}
  y_{\vec{\al}}$ and the torus $\bT \df \bar {\mathfrak{a}} /
  \log(\bar A_{y_{\vec{\al}}})$, and 
  this map conjugates the $a_t$-action on $\bar{A}_{\vec{\al}} y_{\vec{\al}}$ to a
  straightline flow on $\bT$. Moreover, as in Proposition
  \ref{prop:description of measures}, 
  this isomorphism maps the cross-section $\sro$ to a bounded linear
  section. Therefore, for $\cS' = \cB$ or $\cS' = \cS_\vre$, 
   $\{t \in \R:   a_t y_{\vec{\al}}\in \cS'\}$ is a cut-and-project set. Recall from
  Proposition \ref{prop: hitting 
    the subsets} that the set of denominators $q_k$ of convergents to
  $\vec{\alpha}$, and the set of visit times $\{t \geq 0: a_t
  \Lam_{\vec{\alpha}} \in \sro\}$ are related by $t_k = \frac{1}{d}
  \log (q_k)$, and from Proposition \ref{prop: explicit compact push}
  that $y_{\vec{\al}} = q \Lam_{\vec{\alpha}}$ for some $q \in H^-$, the
  contracting horospherical subgroup of $\{a_t\}$. Using this, and the
  fact that $\cS'$ is $\mu_{\sro}$-JM, one can
  prove that for any $\vre>0$, for all large enough $t$, there is a
  bijection between the visit times of the two trajectories $\{a_t
  y_{\vec{\al}}\}, \, \{a_t \Lam_{\vec{\alpha}}\}$  to points
  which are of distance at least $\vre$ from the boundary of $\cS'$; and
  from this, one can deduce that these sets are strongly asymptotic. 

  For the second assertion, we use the notation of 
  \S \ref{subsec: adeles case 2}.
  We let $\Delta \df \on{Stab}_{\tilde{L}_{\vec{\al}}}(\tilde{y}_{\vec{\al}})$.
  In the
  cut-and-project construction, we now replace 
  $\R^r$ with $ \{a_t\} \times M_{\vec{\al}}$, and $R^{s}$ with $\bar A_{\vec{\al}}^{(1)}$, and use
  Proposition \ref{prop: hitting the subsets} and a generalization 
  of \cite[Prop. 2.3]{ASW}.   
  \end{proof}
  \section{Invariance under the weak stable foliation}\label{sec: weak
    stable invariance}
The goal of this section is to explain how our results translate to statements about Lebesgue almost 
every $\theta\in \bR^d$ as they appear in the introduction.
    For $\theta \in \R^d$, let $\wt{\Lam}_\theta$ be as in \eqref{eq: adelic
  lift lattice}. 
Let $\XXnA
$ and $\mu = m_{\XXnA}$ be as in \S
\ref{subsec: adeles case 1},
let $ \widetilde{\cS}_{r_0}$ be as in \eqref{eq: lifted
  cross section}, 
and let $\mu_{\srotilde}$ be the 
cross-section measure. 

It follows easily from the ergodicity of the $a_t$-action with respect 
to $\mu$ that for Lebesgue almost any $\theta$, $\wt{\Lam}_\theta$ is $(a_t,\mu)$-generic. It does not follow
that such points are also $(a_t,\mu_{\srotilde})$-generic. In this section we address this issue and prove the following:
\begin{proposition}\label{prop: weak stable weak}
For Lebesgue a.e. $\theta \in \R^d$,  $\wt{\Lam}_\theta$
is 
$(a_t, \mu_{\srotilde} )$-generic. Moreover, it is $(a_t, \mu_{\srotilde}|_{\wt{\cB}})$-generic as well as
$(a_t, \mu_{\srotilde}|_{\wt{\cS}_\vre} )$-generic for any $\vre\in (0, r_0)$.
\end{proposition}

Recalling the distinction in Definition \ref{def: genericity}, we
define
$$\Omega_1 \df \left \{\wt{\Lam} \in \XXnA: \wt{\Lam} \text{ is } (a_t, 
\mu)\text{-generic} \right \}$$
and
$$
\Omega_2 \df \left \{\wt{\Lam} \in \XXnA: \wt{\Lam} \text{ is } (a_t, 
\mu_{\srotilde})\text{-generic} \right \}.
$$
\begin{remark}
Recall from Proposition \ref{prop:JMtempered} that if $\cS$ were
tempered and reasonable, or if we were interested in $\mu_{\cS}|_E$
for some tempered subset, then we 
would have $\Omega_1 = \Omega_2$, and the arguments of this section
could be avoided. 
\end{remark}

Let $G = \SL_n(\bA)$ and let $H^-$ and $H^0$ be the groups defined in \eqref{eq: the
  contracting group} and \eqref{eq: the centralizer}, and let
$H^{\leq} \df H^-H^0$ be the group as in Proposition \ref{prop:
  horospherical and generic}. 
We refer to
$H^\leq $ as the
{\em weak stable subgroup} corresponding to $\{a_t\}$. Also write
$$
H^+ \df \left \{g \in G: a_t g a_{-t}
  \to_{t \to -\infty} \mathrm{Id} \right\} = \left\{
u(v): v \in 
\R^d \right \}, $$
where
$u(v)$ is as in \eqref{eq: def u(v)}. 
Note that $H^+ \cong \R^d$ is unimodular. We
denote its Haar measure by $m_{H^+}$. This is simply the image of
Lebesgue measure under the map $v \mapsto u(v)$.

We will derive Proposition \ref{prop: weak stable weak} from  the
following three statements.  
\begin{proposition}\label{prop: weak stable 1}
  For any open subgroup $W \subset H^{\leq}$,
  any $W$-invariant
subset $\Omega \subset \XXnA$ of full 
$\mu$-measure, and any $\wt{\Lam} \in \XXnA$,
the set
$
\left\{ h \in H^+ : h\wt{\Lam} \in \Omega \right\}$ is of full $m_{H^+}$-measure. 
\end{proposition}

\begin{proposition}\label{prop: weak stable 2}
The set $\Omega_1$ is $H^{\leq}$-invariant and of full
$\mu$-measure. 
\end{proposition}

\begin{proposition}\label{prop: weak stable 3}
There is an open subgroup $W \subset H^{\leq}$ and a $W$-invariant set
$\Omega_3 \subset \XXnA$ such that the  set
$\Delta_{\srotilde}^{\R}$ defined via \eqref{def:badset} and 
\eqref{eq: thickening} is contained in $\Omega_3$, and
$\mu\left(\Omega_3 \right)=0.$ 
\end{proposition}

\begin{proof}[Proof of Proposition \ref{prop: weak stable weak} assuming Propositions
  \ref{prop: weak stable 1}, \ref{prop: weak stable 2}, \ref{prop:
    weak stable 3}]
  Let $W$ and $\Omega_3$ be as in Proposition \ref{prop: weak stable
    3}.
  By   Theorem \ref{thm:Sgenericity}\eqref{item: tempered i},  
  $$
\Omega_1 \sm \Omega_3 \subset \Omega_1 \sm \Delta_{\srotilde}^{\R} \subset
\Omega_2.
  $$
By Proposition \ref{prop: weak stable 2}, $\Omega_1 \sm \Omega_3$ is
a $W$-invariant set of full measure, and thus by Proposition
\ref{prop: weak stable 1}, $\{\theta \in \R^d: \wt{\Lam}_\theta \in
\Omega_1 \sm \Omega_3\}$ is a set of full measure, proving the first assertion. Moreover, applying Theorem~\ref{thm:Sgenericity} and taking into account the fact that $\Omega_3$ contains $\Del_{\srotilde}$, we see that any $\theta$ in this set is also $(a_t,\mu_{\srotilde}|_{\wt{\cB}})$-generic as well as  $(a_t,\mu_{\srotilde}|_{\wt{\cS}_\vre})$-generic.
  \end{proof}

  \begin{proof}[Proof of Proposition \ref{prop: weak stable 1}]
   Let $G_f = \SL_n(\bA_f), G_\infty =
\SL_n(\R)$ be as in \S 
\ref{subsec: adeles case 1}, and let
\begin{equation}\label{eq: def real leq}
H^\leq_\infty
  \df H^\leq \cap G_\infty.
  \end{equation}
  The groups  $H^+$ and $H^\leq$ are complementary, in the following
  sense: 
  $$H^{\leq } \cap H^+ = \{e\}, \ \  H^\leq = H^\leq_\infty \times G_f, \ \ \dim H^\leq_\infty +
  \dim H^+ = \dim G_\infty.$$
  It follow that there  are neighborhoods $\mathcal{U}_1,
\mathcal{U}_2$ of the 
identity in $H^\leq $ and $H^+$ respectively, such that the map
$$\mathcal{U}_1 \times \mathcal{U}_2 \to G, \ \ \ \  \left(h^\leq
  , h^+ \right)  \mapsto
h^\leq h^+$$
is a homeomorphism onto its image in $G$. For any any $\wt{\Lam} \in \XXnA$, using
the fact that $W$ is open in $H^\leq$ and the stabilizer of $\wt{\Lam}$
in $G$ is discrete, we can replace $\mathcal{U}_1,
\mathcal{U}_2$ with smaller open sets around the identity, to  obtain
that $\cU_1 \subset W$ and the map
$$\mathcal{U}_1 \times \mathcal{U}_2 \to \XXnA, \ \ \ \  \left(w, h^+
\right)  \mapsto w h^+ \wt{\Lam}$$
is a homeomorphism onto its image, which is a neighborhood of $\wt{\Lam}$
in $\XXnA$. Furthermore, by  Lemma \ref{lem:hic basic},
in this neighborhood of $\wt{\Lam}$,
the measure $\mu$ can be written as $d\mu(w h^+ \wt{\Lam}) =
dm^{\mathrm{left}}_{W}(w) dm_{H^+}(h^+)$. Since $\Omega$ is
$W$-invariant, this implies that 
$m_{H^+}$-a.e. $h^+ \in H^+$ satisfies $h^+\wt{\Lam}\in
\Omega$. 
  \end{proof}

  \begin{proof}[Proof of Proposition \ref{prop: weak stable 2}]
By Lemma \ref{lem: mautner}, 
the flow $\{a_t\} \curvearrowright \XXnA$ is ergodic and hence
$\Omega_1$ is of full
$\mu$-measure.
To see that $\Omega_1$ is $H^{\leq}$-invariant, use Proposition
\ref{prop: horospherical and generic}
and the fact that $\mu$ is $H^0$-invariant.
    \end{proof}

For the proof of Proposition \ref{prop: weak stable 3} we
record an observation made previously, in the
  paragraph preceding Proposition \ref{prop: hitting the subsets}:

	\begin{lemma}\label{lem: new lemma about return times}
	Let $\Lam$ be a lattice and let
        $\cY_{\Lam}(\cS_r)\df\set{t\ge 0: a_t\Lam\in \cS_r}$. Then  
	$\set{t\in \cY_{\Lam}(\cS_r):a_t\Lam\notin \cS_r^\sharp}$ is finite.
	\end{lemma}


    \begin{proof}[Proof of Proposition \ref{prop: weak stable 3}]
      Take a sequence $r_i\nearrow \infty$ and define $
      \cS_{r_i}$
      by \eqref{eq:main sets}. Using the projection $\pi:
      \XXnA \to \XX_n$, let 
      $$\widetilde{\cS}_i \df  \pi^{-1}(\cS_{r_i}), \ \ \ \ \Del_i \df
      \Del_{\widetilde{\cS}_i}^\bR.$$ 
      The sets $\widetilde{\cS}_i$ are 
      reasonable cross-sections by Theorem \ref{thm: reasonable case
        I} and Proposition \ref{prop:941}, and the sets $\Del_i$
      satisfy $\mu(\Del_i) =0$ by  Lemma \ref{lem:bad-is-null}.

      Let 
      $$\Om_3\df  \bigcup_{i=1}^\infty \Del_i.$$
      It is clear that $\mu(\Om_3) = 0$. Let 
      $W\df H_0\times K_f$ where $H_0$ is the connected component of the identity in $H_\infty^{\le}$. We 
      claim that $\Om_3$ is $W$-invariant. The invariance under $K_f$ is automatic since each $\Del_i$ is a preimage under $\pi$. The invariance under $a_t$ is clear from the definition of each $\Del_i$. As any element 
      of $H_0$ can be written in the form $h=\smallmat{B&0\\
        \bx&1}a_t$,
      we will be done by
      proving the following:
      \begin{claim}\label{claim: new fix}
      For any $r>0$, $h=\smallmat{B&0\\ \bx&1}\in H_0$, and 
      $\Lam\in \Del_{\cS_r}^\bR$ there exists $r'>r$ such that
      $h\Lam \in \Del_{\cS_{r'}}^\bR$.
    \end{claim}
Let $\cY_\Lam(\cS_r) = \set{t_1<t_2<\dots}$. 
The assumption that $\Lam\in \Del_{\cS_r}^\bR$ is equivalent (see \eqref{def:badset}) to 
saying that there exist $\del>0$ such that for any $\vre>0$ there
exists a subsequence $t_{k_j}$
for which
\begin{equation}\label{eq: two conditions for S le eps}
0<\av{t_{k_j +1}-t_{k_j}}<\vre  \ \ \ \textrm{ and } \ \ \ \limsup_{j
  \to \infty} \frac{j}{t_{k_j}}> \del.
\end{equation}
We will establish a similar statement for $\cY_{\Lam'}(\cS_{r'})$, where $\Lam'\df h\Lam$ and $r'>r$ is large enough.

Let $v_k\in \Lam$ be a sequence of 
vectors satisfying that $a_{t_k}v_k\in D_r$.
We may write 
$$v_k =\mat{\bar{v}_k\\ e^{dt_k}}.$$
Since $a_{t_k}v_k\in D_r$, we have that $\norm{\bar{v}_k}\le
e^{-t_k}r$.
Denoting by $\idist{\cdot, \cdot}$ the standard inner product
on $\R^d$, 
we then have 
$$hv_k = \mat{B\bar{v}_k\\ e^{dt_k}+\idist{\bx,\bar{v}_k}}\in \Lam'.$$
If we let 
$$\tau_k\df  \frac{1}{d}\log(1+\idist{\bx,e^{-dt_k}\bar{v}_k}),$$ then we see that 
$$a_{t_k+\tau_k}hv_k = \mat{e^{t_k+\tau_k}B\bar{v}_k\\ 1}.$$
We claim that $s_k \df t_k+\tau_k\in \cY_{\Lam'}(\cS_{r'})$ for any $r'$ large enough. Indeed, because $B,\bx$
are fixed and by the
bound $\norm{\bar{v}_k}\le e^{-t_k}r$, we have that $\av{\tau_k}\le c$
for some $c>0$. Hence $r'$ could be taken to be $e^c\norm{B}r$. 

 If for some $k$ we have that  
that $s_k= s_{k+1}$, then this means that the lattice
$a_{t_k+\tau_k}\Lam' = a_{t_{k+1}+\tau_{k+1}}\Lam'$ contains the two
distinct vectors $a_{t_k+\tau_k}v_k$ and $a_{t_k+\tau_k}v_{k+1}$ which
are in the disc $D_{r'}$ and so this is a visit time to $\cS_{r'}\sm
\cS_{r'}^\sharp$. From Lemma \ref{lem: new lemma about return times}
we know that this  
can only happen for finitely many values of  $k$.  

Let $\del,\vre, \set{k_j}$ be as in \eqref{eq: two conditions for S le
  eps}. Since $\tau_k$ is bounded (and in fact  
converges to 0),
we deduce that along the same subsequence $k_j$, excluding potentially
finitely many values of $j$,
$$0<\av{s_{k_j+1}-s_{k_j}}<\vre \textrm{ and }\limsup_j \frac{j}{s_{k_j}} = \limsup_j \frac{j}{t_{k_j}}>\del.$$
This shows that $\Lam'\in \Del_{\cS_{r'}}^\bR$.
This establishes Claim \ref{claim: new fix}, concluding the proof of the
Proposition. 
\end{proof}

\section{Concluding the proofs}\label{sec: concluding proofs}
As discussed in \S \ref{sec: lift functionals}, Theorems \ref{thm:
  main best}, \ref{thm: main epsilon} and \ref{thm: number field} all  
follow from Theorem \ref{thm: refinement Horesh bundle}.

\begin{proof}[Proof of Theorem \ref{thm: refinement Horesh bundle}]
  In all cases, the measure $\lambda$ in \eqref{eq: def
    lambda} is finite (see Propositions \ref{prop:description of measures1},
  \ref{prop:description of measures}, \ref{prop:941}) but need not be
  a probability measure. We fix a norm on $\R^d$ and $\vre>0$, and
  define the following probability measures: 
  \[
    \begin{split}
\mu^{(\mathbf{e}_n)} & \text{ is the normalization of } \lambda \text{
  in Case I, best approximations} \\
\mu^{(\mathbf{e}_n, \vec{\alpha})} & \text{ is the normalization of } \lambda \text{
  in Case II, best approximations} \\
\nu^{(\mathbf{e}_n)} & \text{ is the normalization of } \lambda \text{
  in Case I, } \vre\text{-approximations} \\
\nu^{(\mathbf{e}_n, \vec{\alpha})} & \text{ is the normalization of } \lambda \text{
  in Case II, } \vre\text{-approximations}, \\
    \end{split}
  \]
  where in Case II we require that $d \geq 2$ and in Case II for best
  approximations we also require that the norm is either
  the sup-norm, or is $\bar{A}_{\vec{\al}}$-analytic. 

 In Case I, let $\mu = m_{\XXnA}$. Then for Lebesgue
 a.e. $\theta$, by Proposition \ref{prop: weak stable 
   weak}, $\wt{\Lam}_\theta$ is $\left(a_t,
   \mu_{\srotilde}|_{\wt{\cB}}\right)$-generic 
and $\
\left(a_t,
\mu_{\srotilde}|_{\wt{\cS}_\vre} \right)$-generic. 
The map
 $\wt{\psi}$ is continuous but is only defined on
 $\widetilde{\cS}_{r_0}^\sharp$. Since $\widetilde{\cS}_{r_0}^\sharp$
 is open (by Lemma \ref{lem:continuity of vector}) and of full measure
 (by Lemma \ref{lem:Sstar}), $\wt{\psi}$ 
 coincides with a map on $\srotilde$ whose set of discontinuities has
 zero measure with respect to $\mu_{\srotilde} $. Using Lemma
 \ref{lem:convergence-equivalence}, we obtain that  $\widetilde 
 \psi$ maps an equidistributed sequence in
 $\widetilde{\cS}_{r_0}^\sharp $ (with respect to either $
\frac{1}{\mu_{\srotilde}(\wt{\cB})} \mu_{\srotilde}|_{\wt{\cB}} $ or $\frac{1}{\mu_{\srotilde}(\wt{\cS}_\vre)}\mu_{\srotilde}|_{\wt{\cS}_\vre} $) 
 to an equidistributed sequence in 
 $\crly{E}_n \times \R^d \times \widehat{\Z}^n_{\prim}$ (with respect
 to the pushed-forward measure). 
  Using Proposition \ref{prop:
   hitting the subsets}, we
 find that the sequences \eqref{eq:
   seq bundle best}, \eqref{eq: seq bundle eps} are equidistributed
 with respect to  $\mu^{(\mathbf{e}_n)} $ and $\nu^{(\mathbf{e}_n)} $
 respectively. Propositions \ref{prop:
    description of measures I eps}, \ref{prop: splits as a product
    best} and  \ref{prop: properties of the measures} 
 show that the projected measures in 
 the statement of the theorem have the stated properties.

 In Case II, let $\mu = m_{\tilde{L}_{\vec{\al}} \tilde{y}_{\vec{\al}}}$. By Lemma \ref{lem: epsilon zero} any $\vre>\vre_0$ satisfies
 $\mu_{\srotilde}(\wt{\cS}_\vre)>0$.
 By Proposition
 \ref{cor: adelic genericity alpha}, the lattice
 $\wt{\Lam}_{\vec{\alpha}}$ is both $\left(a_t,
 \mu_{\srotilde}|_{\wt{\cB}} \right)$-generic and $\left(a_t,
 \mu_{\srotilde}|_{\wt{\cS}_\vre} \right)$-generic.  
 Using Proposition \ref{prop:
   hitting the subsets}, we
 find that the sequences
  \eqref{eq:
   seq bundle best}, \eqref{eq: seq bundle eps} are equidistributed
 with respect to the pushforwards $\mu^{(\mathbf{e}_n, \vec{\alpha})}
 $ and $\nu^{(\mathbf{e}_n, \vec{\alpha})} $
 respectively. 
The desired properties of these measures are given in Propositions
\ref{prop: description of measures II best} and \ref{prop: description
of measure II eps}. 
  \end{proof}

\begin{remark}\label{remark: for Cantor set}
  It is interesting to find other measures on $\R^d$ which give full
  measure to the set of $\theta$ which 
satisfy the conclusions of Theorems \ref{thm: main best} and \ref{thm:
main epsilon}. In this regard, we note the following.
An examination of the proof of Theorem \ref{thm: refinement Horesh
  bundle} shows that its conclusion holds for every 
$\theta$ for which $\tilde{\Lam}_{\theta}$ is $\left(a_t,
  \mu_{\srotilde}\right)$-generic. Since $\wt{\cB}$
is tempered, by Proposition \ref{prop:JMtempered}, for best approximations
it suffices to show that $\wt{\Lam}_\theta$ is
$(a_t, \mu)$-generic for $\mu = m_{\XXnA}$. It is thus an interesting
question to provide examples of measures on $\R^d$ satisfying that $\wt{\Lam}_\theta$ is
$(a_t, \mu)$-generic for $\mu = m_{\XXnA}$ for a.e.\;$\theta$. Note
that this reduction is not valid in the case of $\vre$-approximations. 

In view of a recent result of Aggarwal and Ghosh (see
  \cite[Cor. 10.8]{Aggarwal_Ghosh}), $\wt{\Lam}_\theta$ is
$(a_t, m_{\XXnA})$-generic if and only if $\Lam_\theta$ is
$(a_t, m_{{\XX}_n})$-generic (we are grateful to an anonymous referee
for pointing this out). Thus it suffices
to find
measures $\mu$ satisfying that $\Lam_\theta$ is $\left(a_t,
  \mu_{\sro}\right)$-generic for $\mu$-a.e.\;$\theta$. Examples of such
measures can be found in \cite{SimmonsWeiss}; they include natural
measures on some self-similar fractals, like the middle-thirds Cantor
set, Sierpinski triangle or Koch snowflake. Furthermore, independently
of \cite[Cor. 10.8]{Aggarwal_Ghosh}, by the strategy employed in this paper, for 
$\theta$ for which $\Lam_\theta$ is $\left(a_t,
  \mu_{\sro}\right)$-generic, the
first two sequences in \eqref{eq: seq best} 
are equidistributed with respect to $\mu^{(\infty)}$.
\end{remark}

\begin{proof}[Proof of Corollary \ref{cor: RR for lines}]
In light of Theorem \ref{thm: refinement Horesh bundle}, it suffices
to compute the measures 
$\lambda^{(\crly{E}_2)}$ in Case I, for both $\vre$-approximations and
best approximations. For $\vre$-approximations we find from
Proposition \ref{prop: description of measures I eps} that
$\lambda^{(\crly{E}_2)} $ is a scalar multiple of $m_{\crly{E}_2}$,
which is just Lebesgue measure on $\left[-\frac12,\frac12\right] \cong
\bT_1$. For best 
approximations we use Propositions \ref{prop: explicit d=1} and
\eqref{eq: RN derivatives Case Ia} to compute the 
density $F$. We leave the details to the reader.  
\end{proof}

In order to prove Corollary~\ref{cor: KL distribution} we isolate the
following lemma. 
\begin{lemma}\label{lem: pushing nu with tau} Let $X$ be lcsc, let $T:
  X \to X$ be a 
  Borel bijection, and let $\nu$ be a $T$-invariant ergodic measure.  
Let $x\in X$ be a $(T,\nu)$-generic point. Let $E\subset X$ be a
$\nu$-JM set with 
$\nu(E)>0$,  let $\tau:E\to \bR$ be  a Borel 
map whose points of discontinuity form a
$\nu$-nullset, and let
$$\nu^{(\tau, E)}\df \frac{1}{\nu(E)}
\tau_*(\nu|_E).$$
Let $\set{k_j}_{j=1}^\infty = \set{k\in \bN: T^kx\in E}$ be the increasing sequence
of visit times to $E$. Then the sequence
$\set{\tau(T^{k_j}x)}_{j\in \bN}$ equidistributes 
with respect to 
$\nu^{(\tau, E)}$. Furthermore, if $\tau$ is integrable on $E$, then the
expectation of $\nu^{(\tau, E)}$ is given by 
$
\frac{1}{\nu(E)}\int_E\tau d\nu. $ 
\end{lemma}
\begin{proof}
Since $E$ is $\nu$-JM, the sequence $T^{k_j}x$ equidistributes with
respect to the normalized restriction  
$\frac{1}{\nu(E)}\nu|_E$. Then, by  \cite[Thm. 2.7]{Billingsley}, the
fact that $\tau$ is $\nu$-almost everywhere continuous gives 
that $\tau(T^{k_j}x)$ equidistributes with respect to $\nu^{(\tau, E)}
= \tau_*(\frac{1}{\nu(E)}\nu|_E)$. For the last assertion we compute the expectation:
$$\int x \, d\nu^{(\tau, E)} 
= \int_E\tau(y) \,
d\left(\frac{1}{\nu(E)}\nu|_E \right)(y) = \frac{1}{\nu(E)}\int_E \tau d\nu.$$
\end{proof}
\begin{proof}[Proof of Corollary \ref{cor: KL distribution}]
Let $\crly{X}_n$ and $\mu = m_{\crly{X}_n}$ be as in 
\S\ref{sec: rhc1},
let $\sro$ be as in \S\ref{sec:main example} and let $\mu_{\sro}$ be the 
cross-section measure. Let $E$ be one of 
$E = \cB$ or $E=\cS_\vre$, for some  
$\vre \in (0, r_0)$. Then $E$ is  
 $\mu_{\sro}$-JM and satisfies
$\mu_{\sro}(E)>0$ (see \S \ref{subsec: 8.3}). Consider the 
dynamical system     
$\left(\sro, T_{\sro},  \frac{1}{\mu_{\sro}(\sro)}\mu_{\sro} \right)$,
where $T_{\sro}$ is the first return map.
With this notation the first return time and first return map, defined
by
$$\tau_E(x) \df \min\set{t > 0:
    a_t(x)\in E}, \ \  \ T_E(x) \df  a_{\tau(x)}x \ \ \ (x \in E)$$ satisfy:
\begin{itemize}
\item For any generic point $x$, $\tau_E(x)$ and $T_E(x)$  
are well-defined (because $E$ is of positive
measure and $\mu_{\sro}$-JM). 
\item $\int_E \tau_E \, d\mu_{\sro} = 1$ (by
  Theorem~\ref{cor:csmeasure}\eqref{item: Kac formula}).
\item
$\tau_E$ is continuous outside a
$\mu_{\sro}$-null set (this can be proved using similar arguments to those in
\S\ref{subsec: 8.3}, we leave the details to the reader). 
  
\end{itemize}

Let 
 $\Lam_\theta$ be the lattice defined in 
 \eqref{eq: def lattice to a vector} and \eqref{eq: def u(v)}. 
It follows from Proposition~\ref{prop: weak stable 3} in conjunction
with Theorem \eqref{thm:Sgenericity}\eqref{item: tempered i} 
 that for Lebesgue almost any $\theta \in \R^d$, the points in $E$
 along the orbit $\{a_t \Lam_\theta\}$ are 
 generic in this dynamical system. We now take such a point 
$x= a_{s_0}\Lam_\theta \in E$, for $\theta$ in this subset of full measure in
$\R^d$, and with $s_0>0$ the smallest possible visit time to $E$.
Define the visit times by 
$$s^{(E)}_k \defi
\tau_E\left(T_E^k(x) \right).$$


Applying Lemma~\ref{lem: pushing nu with tau}, we deduce that 
the sequence 
$\set{s^{(E)}_k}_{k\in \bN}$ equidistributes with respect to 
$$\lam^{(E)} \df \frac{1}{\mu_{\sro}(E)} (\tau_E)_*(\mu_{\sro}|_E).$$
By the correspondence established in Proposition 
 \ref{prop: hitting the subsets}, the visits to $E$
 correspond to an approximation vectors 
 having $n$-th coordinate 
 $$q_k = \exp\left(d t_k\right), \ \ \ \text{ where } t_k =
 \sum_{i=0}^k s^{(E)}_i  = s_0 + T^k_E(x).$$
 From this we immediately obtain \eqref{eq: seq KL}, and see that
 assertion
 \eqref{item: Cor KL 2} of Corollary \ref{cor: KL distribution} holds.
 Lemma~\ref{lem: pushing nu with tau} also gives that
$\lam^{(E)}$ has finite expectation which we denote $\ga^{(E)}$; that
is, 
$$\ga^{(E)}\defi\int_\bR xd\lam^{(E)} =\frac{1}{\mu_{\sro}(E)}
\int_E\tau_E d\mu_{\sro} = \frac{1}{\mu_{\sro}(E)}.$$

We now prove \eqref{item: Cor KL 1} of Corollary \ref{cor: KL distribution}. 
We have shown
in Proposition~\ref{prop: weak stable 3} in conjunction with
Proposition~\ref{prop:JMonallS} that for Lebesgue almost any $\theta$,  
\begin{equation}\label{eq: the correct assymp}
\frac{1}{T}N(\Lam_\theta,E,T)\to \mu_{\sro}(E).
\end{equation}
Restrict further to the set of $\theta$ for which \eqref{eq: the
  correct assymp} holds.
Equation~\eqref{eq: the correct assymp} and 
$\mu_{\sro}(E)>0$ give 
\begin{align*}
\lim_{k \to \infty} \frac{1}{dk}\log q_k = \lim_{k \to \infty}\frac{ t_k }{k} =
  \frac{t_k}{N(\Lam_\theta, E, t_k)}= \frac{1}{\mu_{\sro}(E)}.  
\end{align*}
This immediately gives item \eqref{item: Cor KL 1} of the Corollary.
Finally, 
in the case of $\cS_\vre$, we can compute the expectation explicitly: 
  $$\ga^{(KL)}_{\norm{\cdot},\vre} = \frac{1}{\mu_{\sro}(\cS_\vre)} \stackrel{\mathrm{Prop. }
     \ref{prop:description of measures1}}{=} 
  \frac{\zeta(n)}{d\cdot m_{\XX_n}(\XX_n) \, m_{\R^d}(B_\vre)} =
  \frac{\zeta(n)}{d\cdot V_{d, \|\cdot\|} \, \vre^d}. $$
    \end{proof}

\section{Conclusions regarding equidistributed compact torus orbits}
Until now we have dealt with two types of $a_t$-invariant and ergodic measures:
 $\mu = m_{\XX_n^{\bA}}$ (Case I) and $\mu_{\wt L_{\vec{\al}} \tilde y_{\vec{\al}}}$ (Case II). We proved results regarding the best approximations and $\vre$-approximations of Lebesgue almost any $\theta$ (Case I) and the algebraic vector $\vec{\al}$ (Case II). The two discussions regarding Case I and Case II were carried out simultaneously but were completely independent to one another. We now wish to establish a connection between the two cases which occurs when one varies the algebraic vector
$\vec{\al}$ under some assumptions. First we need the following.
\begin{lemma}\label{lem:lifts equidistribute}
We have that 
$$\set{\mu\in \cP(\XX_n^{\bA}) : \mu \textrm{ is $a_1$-invariant and }\pi_*\mu = m_{\XX_n}} = \set{m_{\XX_n^{\bA}}}.$$ 
As a consequence, 
if $\mu_k\in \cP(\XX_n^{\bA})$ is a sequence of $a_1$-invariant probability measures such that $\pi_*\mu_k\to m_{\XX_n}$, then
$\mu_k\to m_{\XX_n^{\bA}}$.
\end{lemma}
\begin{proof}
The proof relies on two facts about entropy. The probability measures $m_{\XX_n}$ and $m_{\XX_n^{\bA}}$
are the unique $a_1$-invariant measures of maximal entropy for the dynamical systems $(\XX_n,a_1), (\XX_n^{\bA},a_1)$, and in fact
they have equal entropies:
$$\on{h}(a_1, m_{\XX_n}, \XX_n) = \on{h}(a_1, m_{\XX_n^{\bA}}, \XX_n^{\bA}) = n.$$
For the space $\XX_n$ this follows from \cite[Corollary 7.10]{EL}. For the space $\XX_n^{\bA}$ this follows from \cite[Theorem 7.9]{EL} together with the unique ergodicity of the $\SL_n(\bR)$ action on $\XX_n^{\bA}$.

Assume that $\mu\in\cP(\XX_n^{\bA})$ satisfies $\pi_*\mu = m_{\XX_n}$. Then, since $\pi$ is a factor map,
 $$\on{h}(a_1, \mu, \XX_n^{\bA}) \ge \on{h}(a_1, \mu, \XX_n^{\bA}) = n.$$
From the uniqueness of measure of maximal entropy we deduce that $\mu = m_{\XX_n^{\bA}}$ as desired.

Now let $\mu_k\in \cP(\XX_n^{\bA})$ be a sequence of $a_1$-invariant probability measures satisfying that $\pi_*\mu_k\to m_{\XX_n}$. It follows that any weak* accumulation point of $\mu_k$ must be an $a_1$-invariant probability measure that projects under $\pi$ to $m_{\XX_n}$. By the first part of lemma, $\mu_k$ has only one accumulation point - $m_{\XX_n^{\bA}}$. We conclude that $\mu_k\to m_{\XX_n^{\bA}}$ as desired.
\end{proof}

\begin{theorem}\label{thm:for SZ}
Assume $d\ge 2$ and that the norm $\norm{\cdot}$ on $\bR^d$ we use to define the notion of best approximations is either 
the Euclidean norm or the sup norm.
Let $\vec{\al}_k = (\al_{k1},\dots, \al_{kd})^{\on{t}}\in \bR^d$ be a sequence of vectors such that for each $k$,
$\set{1,\al_{k1},\dots,\al_{kd}}$ span a totally real number field of degree $n$ over $\bQ$. 
Let $\bar{h}_{\vec{\al}_k}$ be as in~\eqref{eq: def Bar B} and $x_{\vec{\al}_k}\in \XX_n$ be as in~\eqref{eq:x alpha}. If the sequence of 
periodic probability measures supported on the periodic orbits $\bar{h}_{\vec{\al}_k}^* Ax_{\vec{\al}_k}$ converge weak*
to $m_{\XX_n}$, then 
\begin{equation}\label{eq:com2}
\mu^{(\mb{e}_n,\vec{\al}_k)} \longrightarrow \mu^{(\mb{e}_n)},
\end{equation}
where $\mu^{(\mb{e}_n,\vec{\al}_k)}, \mu^{(\mb{e}_n)} \in \cP(\crly{E}_n\times \bR^d\times \wt{\bZ}^n)$ are  
the probability measure corresponding to best approximations of $\vec{\al}_k$ (resp. Lebesgue almost any $\theta$) by Theorem~\ref{thm: refinement Horesh bundle}. 
\end{theorem}
\begin{proof}
To reduce notational clutter we let $x_k \defi x_{\vec{\al}_k}$ and $\bar{h}_k\defi \bar{h}_{\vec{\al}_k}$. We assume
that the sequence of homogeneous measures of the periodic orbits $\bar{h}_k^*Ax_k$ converge weak* to $m_{\XX_n}$. By applying the involution $x\mapsto x^*$ of $\XX_n$ we deduce that the same holds for the sequence of periodic orbits
$\bar{h}_kAx^*_k = \bar{A}_{\vec{\al}_k} y_{\vec{\al}_k}$, where $\bar{A}_{\vec{\al}_k}$ and $y_{\vec{\al}_k}$ are as defined in~\eqref{eq: def bar A}. That is 
\begin{equation}\label{eq:again measures equidistribute}
m_{ \bar{A}_{\vec{\al}_k} y_{\vec{\al}_k}}\longrightarrow m_{\XX_n}.
\end{equation}
Consider the periodic orbits $\tilde{L}_{\vec{\al}_k}\tilde{y}_{\vec{\al}_k}$ in $\XX_n^\bA$ appearing in 
Proposition \ref{prop: lifted measure case II}. 
Let us denote
\begin{equation}\label{eq:convergence of measures 2}
\mu^{(k)} \defi m_{\tilde{L}_{\vec{\al}_k} \tilde{y}_{\vec{\al}_k}};\;\;\; \mu \defi m_{\XX_n^{\bA}},
\end{equation}
By Proposition \ref{prop: lifted measure case II},  
$\pi_*\mu^{(k)} = m_{\bar{A}_{\vec{\al}_k}y_{\vec{\al}_k}}$. Thus, \eqref{eq:again measures equidistribute} says that 
$\pi_*\mu^{(k)}\to m_{\XX_n}$. Furthermore,
$\mu^{(k)}$ is $a_t$-invariant and so by Lemma~\ref{lem:lifts equidistribute} we deduce that 
\begin{equation}\label{asdfasf}
\mu^{(k)}\longrightarrow \mu
\end{equation}

Let 
$\mu^{(k)}_{\srotilde},  \mu_{\srotilde}$ be the corresponding cross-section measures. 
We apply Proposition~\ref{prop:continuity}, where for the tempered subset  we take $\wt{\cB}\subset \srotilde$.
Regarding the applicability of Proposition~\ref{prop:continuity}, we note that if $\nu$ denotes any of the measures in \eqref{eq:convergence of measures 2}, 
 then the $\nu$-reasonability of $\srotilde$, the $\nu$-Jordan measurability of $\cB$ and its temperedness follow from 
 Theorem \ref{thm: MKf} (note that this is the point where we use the fact that the norm on $\bR^d$ is either the Euclidean norm or the sup norm). 
 
 The conclusion of Proposition~\ref{prop:continuity} is that 
 $\mu^{(k)}_{\srotilde}|_{\wt{\cB}} \longrightarrow  \mu_{\srotilde}|_{\wt{\cB}}$, where the convergence is tight convergence. 
 In particular, the convergence holds after renormalizing the restricted measures to be probability measures:
 $$\frac{1}{\mu^{(k)}_{\srotilde}(\wt{\cB})}\mu^{(k)}_{\srotilde}|_{\wt{\cB}} \longrightarrow \frac{1}{\mu_{\srotilde}(\wt{\cB})} \mu_{\srotilde}|_{\wt{\cB}}. $$
 Let $\wt{\psi}$ be as in \eqref{eq: augmented map}. By definition (see the proof of Theorem~\ref{thm: refinement Horesh bundle} in \S\ref{sec: concluding proofs}), 
 $$\wt{\psi}_*\pa{ \frac{1}{\mu_{\srotilde}(\wt{\cB})} \mu_{\srotilde}|_{\wt{\cB}}} =  \mu^{(\mb{e}_n)};\;\;
 \wt{\psi}_*\pa{\frac{1}{\mu^{(k)}_{\srotilde}(\wt{\cB})}\mu^{(k)}_{\srotilde}|_{\wt{\cB}}} = \mu^{(\mb{e}_n, \vec{\al}_k)}.$$
 By continuity of $\wt{\psi}_*$ at $\mu_{\srotilde}$ (see Lemma~\ref{lem:convergence-equivalence}) we obtain \eqref{eq:com2}.
\end{proof}
\begin{remark}
At the moment we do not have a version of Theorem~\ref{thm:for SZ} for $\vre$-approximations. The reason is that Proposition~\ref{prop:continuity}
requires temperedness and we only know temperedness for $\cB$ and not for $\cS_{\vre}$.
\end{remark}

\bibliographystyle{alpha}
\bibliography{uribib}
\printindex

\end{document}